\documentclass[11pt]{article}
\pdfoutput=1

\usepackage[utf8]{inputenc}
\usepackage[english]{babel}
\usepackage{authblk}
\usepackage{graphicx}        
\usepackage{multicol}        
\usepackage[bottom]{footmisc}
\usepackage{placeins}
\usepackage{changes,todonotes}
\usepackage{pgf}
\usepackage{psfrag}
\usepackage{pgfplots}  
\usepackage{hyperref}
\usepackage{appendix}
\usepackage{amsmath,amssymb}
\usepackage{bm}
\usepackage{mathtools}  
\usepackage{color}
\usepackage{stmaryrd}
\usepackage{amsthm}
\usepackage{blindtext}
\usepackage[labelfont=bf]{caption}
\usepackage{subcaption}
\usepackage[margin=3cm]{geometry}  
\usepackage[ruled,vlined]{algorithm2e}
\usepackage{cases}
\usepackage{stmaryrd}
\usepackage{xparse}
\usepackage{arydshln}
\usepackage[sorting=nty, firstinits=true]{biblatex}
\graphicspath{{./fig/}}

\renewcommand{\vec}[1]{\boldsymbol{#1}}

\NewDocumentCommand{\dgal}{sO{}m}{%
	\IfBooleanTF{#1}
	{\dgalext{#3}}
	{\dgalx[#2]{#3}}%
}
\NewDocumentCommand{\dgalext}{m}{%
	\sbox0{%
		\mathsurround=0pt 
		$\left\{\vphantom{#1}\right.\kern-\nulldelimiterspace$%
	}%
	\sbox2{\{}%
	\ifdim\ht0=\ht2
	\{\kern-.45\wd2 \{#1\}\kern-.45\wd2 \}%
	\else
	\left\{\kern-.5\wd0\left\{#1\right\}\kern-.5\wd0\right\}%
	\fi
}
\NewDocumentCommand{\dgalx}{om}{%
	\sbox0{\mathsurround=0pt$#1\{$}%
	\sbox2{\{}%
	\ifdim\ht0=\ht2
	\{\kern-.45\wd2 \{#2\}\kern-.45\wd2 \}%
	\else
	\mathopen{#1\{\kern-.5\wd0 #1\{}
	#2
	\mathclose{#1\}\kern-.5\wd0 #1\}} 
	\fi  
}
 
\newtheorem{theorem}{Theorem}[section]

\newtheorem{remark}[theorem]{Remark}    

\begin{document}
	
\title{Local Embedded Discrete Fracture Model (LEDFM)}

\author[$\diamond,\ddagger$]{Davide~Losapio}
\author[$\diamond$]{Anna~Scotti}

\affil[$\diamond$]{MOX, Laboratory for Modelling and Scientific Computing, Department of Mathematics, Politecnico di Milano, Via Edoardo Bonardi 9, 20133 Milan (MI), Italy.}
\affil[$\ddagger$]{Corresponding author.}

\affil[ ]{\texttt {davide.losapio@polimi.it, anna.scotti@polimi.it}}

\maketitle
  
\begin{abstract}
\noindent
The study of flow in fractured porous media is a key ingredient for many geoscience applications, such as reservoir management and geothermal energy production. Modelling and simulation of these highly heterogeneous and geometrically complex systems require the adoption of non-standard numerical schemes.
The \emph{Embedded Discrete Fracture Model} (EDFM) is a simple and effective way to account for fractures with coarse and regular grids, but it suffers from some limitations: it assumes a linear pressure distribution around fractures, which holds true only far from the tips and fracture intersections, and it can be employed for highly permeable fractures only.
In this paper we propose an improvement of EDFM which aims at overcoming these limitations computing an improved coupling between fractures and the surrounding porous medium by a) relaxing the linear pressure distribution assumption, b) accounting for impermeable fractures modifying near-fracture transmissibilities.
These results are achieved by solving different types of local problems with a fine conforming grid, and computing new transmissibilities (for connections between fractures and the surrounding porous medium and those through the porous medium itself near to the fractures). Such local problems are inspired from numerical upscaling techniques present in the literature.
The new method is called \emph{Local Embedded Discrete Fracture Model} (LEDFM) and the results obtained from several numerical tests confirm the aforementioned improvements. 
\end{abstract}


\section{Introduction}

The simulation of flow and transport in fractured porous media is an essential ingredient for many geoscience applications, such as reservoir management tasks, that include planning hydrocarbon production or Enhanced Oil Recovery (EOR) operations. 
Other examples of relevant applications in this field include CO$_2$ storage and sequestration, water resources management, nuclear waste disposal and geothermal energy production. 

A typical geological porous medium is heterogeneous on different length scales. The material properties, such as rock permeability, may vary locally by many orders of magnitude, and the geometry is typically complex. An important source of heterogeneity are fractures, which are thin inclusions whose conductivity can be both orders of magnitude higher and lower than the surrounding material, called porous matrix. Hence, they either serve as preferential paths or barriers for the fluid flow. 
The fracture aperture is typically several orders of magnitude smaller than the characteristic size of the porous medium domain, and the fracture characteristic length.
See \cite{Singhal, Dietrich, Adler, Berkowitz2, Sahimi} for basic knowledge on flow and transport in fractured porous media. 

Due to the strong heterogeneities and complicated geometries characterizing reservoirs, it is often challenging to correctly simulate such systems.
In particular, fractures, due to their features, require special treatments in reservoir simulation, and they should be modelled adequately depending on the considered length scale \cite{Tatomir, Assteerawatt}.
Typically, from small to mid length scales they are explicitly described through Discrete Fracture Models (DFM), that often model fractures as lower-dimensional objects, meaning that they are represented by $(n-1)$-dimensional objects, where $n$ refers to the dimension of the porous matrix domain \cite{Alboin2, Martin, Formaggia2, Schwenck}.
On large scales, instead, Continuum Fracture Models (CFM) are often used.
DFM typically provide more accurate results than CFM, but, since fractures are ubiquitous, they cannot be used to model every single fracture present in a certain domain. Hence, there exist Hybrid Fracture Models (HFM), that use DFM to model the most dominant, largest fractures, while the remaining ones are handled with CFM \cite{Lee}. 

Depending on how the grid generation is performed, two different classes of numerical methods for DFM can be mainly distinguished.
On one side we have conforming methods, that honour the fractures geometry exploiting the flexibility of unstructured grids \cite{Reichenberger2, Angot, Karimi-Fard2, Sandve}. However, for complex fracture networks conformity may result in having small matrix cells near the fractures, and eventually large, ill-conditioned linear systems.
On the other side we have embedded methods, that, in order to overcome the limitations associated to conforming methods, usually adopt structured grids where fractures are allowed to cut arbitrarily the matrix grid. However, transfer functions taking into account the fracture-matrix coupling must be added to the model, and these methods are often based on some restrictive assumptions.

Among the latter we mention the Embedded Discrete Fracture Model (EDFM) \cite{Lee, Li}, which is a simple and effective way to account for fractures with coarse and regular grids, but it suffers from some limitations: the expressions for the flux interaction terms between matrix and fracture domains stems from the assumption of linear pressure distribution around fractures, which holds true only far from the tips and fracture intersections. Moreover, it can be applied for highly permeable fractures only.

The method has been extended to the three-dimensional case in \cite{Moinfar} within a compositional reservoir simulator, and it can handle fractures arbitrarily oriented in space, allowing to perform simulations for geometrically complex fractured reservoirs.
Nonplanar fractures are taken into account in \cite{Xu2}, where the issue of very small fracture segments is tackled as well to avoid excessive limitations of the timestep and preconditioning problems.
In \cite{Xu3} EDFM has been extended to corner-point grids, widely used in the industry to better represent geological features, and to the case of full-permeability-tensors.  \cite{Xu} focuses on geometrical preprocessing, accounting for special limit cases, and proposes a geometrical algorithm to find the intersections between a general polyhedron and a general polygon with the aim of determining the intersections between corner-point cells and fractures.

Various improvements have been proposed for the EDFM method.
In \cite{Tene} the Projection-based Embedded Discrete Fracture Model (pEDFM) is introduced to take into account even the case of ``impermeable'' fractures, meaning by that fractures whose normal permeability is lower than that of the porous matrix.
A modified version of pEDFM has been proposed in \cite{Jiang}, which differs slightly in the definition of matrix-fracture transmissibilities.
In \cite{Olorode} pEDFM has been extended to the 3D case for compositional simulation of realistic fractured reservoirs and in \cite{HosseiniMehr} it has been applied to corner-point grids to simulate both flow and heat transfer phenomena.

Other contributions focused on getting a better representation of the matrix-fracture flow with respect to the original version of EDFM, which simply assumes a linear pressure distribution around fractures, to derive improved matrix-fracture transmissibility expressions.
This is particularly relevant to get a more accurate modelling of transient effects, which is a critical issue in the case of very low matrix permeabilities.
To this end, in \cite{Cao} the matrix-fracture flow expressions are derived using the Boundary Element Method (BEM) in a two-dimensional framework rather than using the linear approximation, thus allowing to consider the effect of local grid geometry and boundary pressure.
This work has been later extended in \cite{Rao} to the 3D case and the resulting method was successfully applied to production simulation of a multi-stage fractured horizontal well.
In \cite{Ding} the pressure distribution around fractures is computed using an integral method and assuming steady-state flow. In this way, a suitable and accurate representation of the pressure is also obtained at fracture tips and intersections. New matrix-fracture transmissibilities are then computed based on the steady-state near-fracture pressure solution. 
In \cite{Shao} the Integrally Embedded Discrete Fracture Model (iEDFM) is presented, where the transmissibilities between matrix and fracture cells are computed with a semi-analytic method in which fractures are represented through a collection of point sinks. The near-fracture pressure distribution is obtained superimposing the pressure solutions of the sinks and is then used to compute the new transmissibilities.

The original EDFM employs the classic Two-Point Flux Approximation (TPFA) scheme to approximate the fluxes at porous matrix interfaces. However, it is well known that this scheme is not consistent in the case of non-orthogonal grids and/or anisotropic media characterized by a full permeability tensor. This issue is addressed in \cite{Nikitin2}, where the monotone embedded discrete fracture method (mEDFM) is presented.
This method couples the original EDFM with two different nonlinear schemes: the monotone TPFA scheme \cite{Danilov} and the compact multi-point flux approximation scheme \cite{Chernyshenko, Terekhov}, which satisfies the discrete maximum principle.
In particular, the nonlinear schemes are used to approximate the fluxes only at matrix faces, but this is sufficient to get accurate results also in the cases where TPFA fails.
The method has been later extended in \cite{Yanbarisov} to account also for blocking fractures, and called Projection-based monotone embedded discrete fracture method. Indeed it couples the pEDFM with the nonlinear compact multi-point flux approximation scheme.

This paper focuses on the development of a new method, belonging to the general class of embedded methods, called Local Embedded Discrete Fracture Model (LEDFM), which is capable to overcome both the limitations of the EDFM, here applied to the single-phase fluid flow case and in a two-dimensional framework.
The method adopts local flow based upscaling methods, see \cite{Durlofsky2}, to compute new matrix-fracture and near-fracture matrix-matrix coarse scale transmissibilities. Here the coarse model coincides with an EDFM model, whereas in classical numerical upscaling techniques with fracture networks found in the literature \cite{Karimi-Fard, Gong2, Tatomir2, Fumagalli2} the coarse model is usually a MINC (i.e. a continuum approach), which is not capable of explicitly representing fractures.
The definitions of the local fine scale problems for the new upscaled transmissibilities are inspired from the aforementioned techniques, and conforming methods are used to solve them \cite{Karimi-Fard2, Sandve}.
Unlike previous works fractures having any conductivity contrast with respect to the porous matrix are considered in the local fine scale problems.

Numerical tests, comparing the solutions of different embedded methods, including the newly developed LEDFM, with respect to reference solutions, show that LEDFM overcomes both the limitations of the classic EDFM method: the matrix-fracture coupling is improved by relaxing the linear pressure distribution assumption and it can be successfully applied for both permeable and impermeable fractures in most cases.

As we will see, in some cases a higher accuracy of the description of the near fracture flow is needed, so that the local problems for the computation of matrix-matrix transmissibilities are substituted with a multiscale approach \cite{Jenny, Hesse}.

The paper is organized as follows: in Section~\ref{sec:goveqns} the hybrid dimensional formulation of the Darcy problem in the incompressible single-phase flow case is presented and the EDFM and pEDFM methods are recalled for readers' convenience.
In Section~\ref{sec:LEDFMform} the new LEDFM method is described in detail while in Section~\ref{sec:numtests} several numerical tests are presented to validate the local method. Convergence analyses for different fracture geometries and matrix-fracture permeability contrasts are run for the purpose. Moreover, a tracer transport test is considered to highlight possible differences in the Darcy flow given by the different methods.
Finally, in Section~\ref{sec:conclusions} the main results are summarized and possible future developments and perspectives are discussed.

\section{Governing Equations and Discretization} \label{sec:goveqns}

In this work we focus on the case of incompressible single-phase Darcy flow in a two-dimensional porous medium domain $\Omega$ having a single immersed straight fracture, but the case of a fracture cutting the entire domain can be easily obtained modifying the model described below, as well as the case of multiple, non intersecting fractures. Moreover, gravity effects are neglected.
The work focuses on the two-dimensional case for reasons of simplicity in the presentation of  the new local embedded method, described in Section~\ref{sec:LEDFMform}, and of computational cost, but the method can be in principle extended to the three-dimensional case.

The domain $\Omega$ is divided into two connected subsets, namely $\Omega_m$ for the porous matrix part and $\gamma$ for the fracture, where $\gamma$ is one-dimensional. The fracture is assumed to be lower dimensional with respect to the matrix part. We require $\Omega_m$ and $\gamma$ to be disjoint, i.e.
\begin{equation*}
	\overline{\Omega} = \overline{\Omega}_m \cup \overline{\gamma}
	\quad
	\textnormal{and}
	\quad
	\mathring{\Omega}_m \cap \mathring{\gamma} = \emptyset.
\end{equation*}
We denote the boundary of the porous medium domain as $\Gamma \vcentcolon = \partial \Omega$ and we subdivide it into two disjoint parts, $\Gamma^D$ and $\Gamma^N$, such that
\begin{equation*}
	\overline{\Gamma} = \overline{\Gamma}^D \cup \overline{\Gamma}^N
	\quad
	\textnormal{and}
	\quad
	\mathring{\Gamma}^D \cap \mathring{\Gamma}^N = \emptyset.
\end{equation*}
On $\Gamma^D$ Dirichlet conditions are enforced, while on $\Gamma^N$ we set Neumann conditions.
$\vec{n}_\Gamma$ is the unit outward normal to the boundary of the domain $\Omega$.

\begin{figure}[!t]
	\centering
	\includegraphics[width=.4\textwidth]{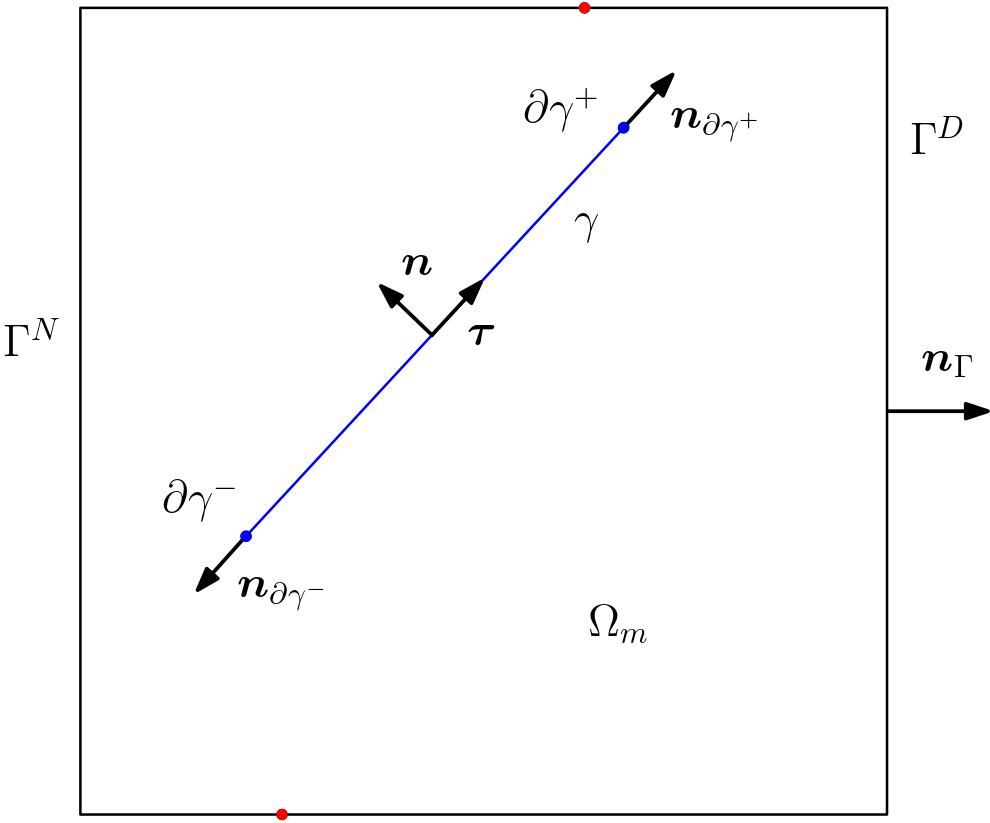}
	\caption{Schematic representation of the domain $\Omega$. The matrix subdomain is denoted by $\Omega_m$ and the fracture by $\gamma$.}
	\label{fig:EDFM_domain}
\end{figure}

The boundary of $\gamma$ is composed by two points $\partial \gamma = \{ \partial \gamma^{+}, \partial \gamma^{-}  \}$, i.e. the tips of the fracture, and we define $\vec{n}_{\partial \gamma^{+}}$ and $\vec{n}_{\partial \gamma^{-}}$ as the unit outwards normals to $\partial \gamma^{+}$ and $\partial \gamma^{-}$, respectively.
We also indicate with $\boldsymbol{\tau}$ the unit tangential vector on $\gamma$ directed from  $\partial \gamma^{-}$ to $ \partial \gamma^{+}$ and with $\vec{n}$ the unit normal vector to the fracture, such that $(\boldsymbol{\tau}, \vec{n})$ are positively oriented.
The domain is depicted in Fig.~\ref{fig:EDFM_domain}.

The hybrid dimensional formulation of the Darcy problem can be written as
\begin{subnumcases}{\label{eq:EDFM}}
	-\nabla\cdot\left(\dfrac{\vec{K}_m}{\mu}\nabla p_m\right) + q_{mf} = Q_m & in $\Omega_m \subset \mathbb{R}^2$, \label{eq:EDFM_m}
	\\
	-\nabla_{\boldsymbol{\tau}} \cdot \left(\dfrac{K_{\gamma,\boldsymbol{\tau}}}{\mu}\nabla_{\boldsymbol{\tau}} p_\gamma\right) + q_{fm}
	=
	Q_{\gamma} & in $\gamma \subset \mathbb{R}$, \label{eq:EDFM_f}
\end{subnumcases}
along with the boundary conditions
\begin{subnumcases}{\label{eq:EDFM_BCs}}
	p_m = p_{D,m} \quad & on  $\Gamma^D$,
	\\
	-\left(\dfrac{\vec{K}_m}{\mu}\nabla p_m\right) \cdot \vec{n}_{\Gamma} = g_{N,m} \quad & on $\Gamma^N$,
	\\
	-\left(\dfrac{K_{\gamma,\boldsymbol{\tau}}}{\mu}\nabla_{\boldsymbol{\tau}} p_\gamma\right) \cdot \vec{n}_{\partial \gamma^{-}} = 0 \quad & on $\partial \gamma^{-}$, \label{eq:EDFM_tipBC1} 
	\\
	-\left(\dfrac{K_{\gamma,\boldsymbol{\tau}}}{\mu}\nabla_{\boldsymbol{\tau}} p_\gamma\right) \cdot \vec{n}_{\partial \gamma^{+}} = 0 \quad & on $\partial \gamma^{+}$, \label{eq:EDFM_tipBC2} 
\end{subnumcases} 
where $p_m$ and $p_\gamma$ are the matrix and fracture pressure, respectively, $\mu$ the fluid dynamic viscosity, $q_{mf}$ and $q_{fm}$ the so-called \emph{flux interactions} between matrix and fracture and $Q_m$ and $Q_\gamma$ the source terms for matrix and fracture, respectively. Moreover, $\nabla_{\boldsymbol{\tau}} \cdot$ and $\nabla_{\boldsymbol{\tau}}$ denote, respectively, the divergence and gradient operators on the line tangent to the fracture.

$\vec{K}_m$ is the matrix permeability tensor, which is symmetric and positive definite, while the fracture permeability is assumed to be orthotropic, meaning that it is possible to identify a normal permeability $K_{f,\vec{n}}$ and a tangential one $K_{f,\boldsymbol{\tau}}$ with respect to the fracture tangential plane. In particular, in~\eqref{eq:EDFM}--\eqref{eq:EDFM_BCs} a key role is played by the equivalent tangential permeability $K_{\gamma, \boldsymbol{\tau}} = K_{f,\boldsymbol{\tau}} d$, where $d$ is the fracture aperture, assumed to be constant. Clearly, the permeability values in the matrix and fracture are likely to differ of several orders of magnitude.

Note that no flow conditions~\eqref{eq:EDFM_tipBC1}--\eqref{eq:EDFM_tipBC2} are adopted at the fracture tips, which is a common choice due to its simplicity from an implementation point of view, but the physical motivation behind it is not obvious.
Finally, note that, to close the problem, we need to provide an expression for the interaction terms $q_{mf}$ and $q_{fm}$.

\begin{figure}[!t]
	\centering
	\includegraphics[width=.4\textwidth]{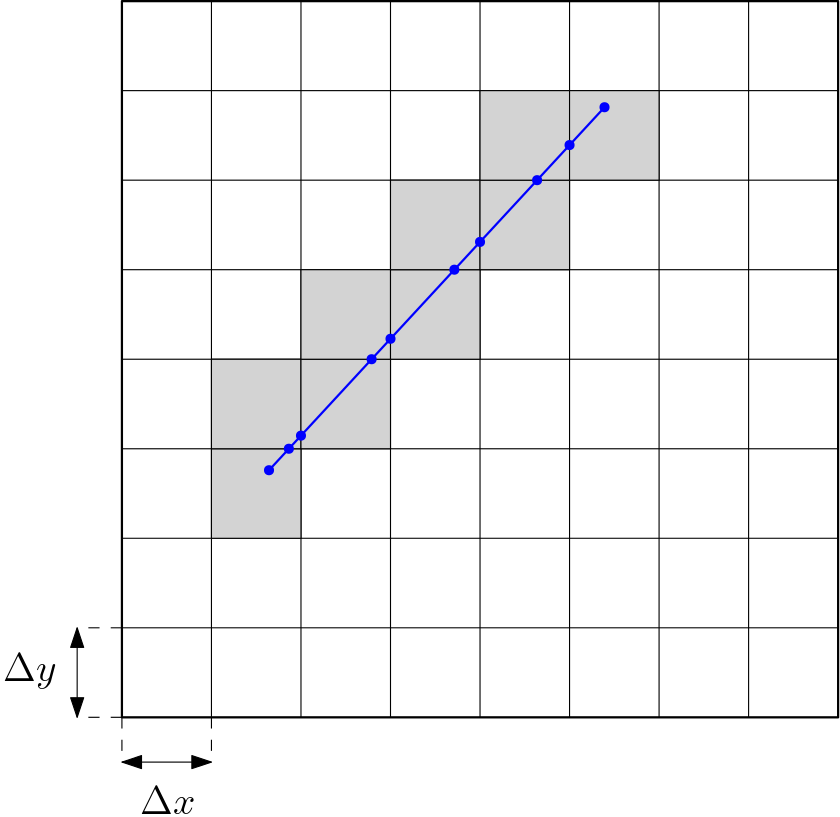}
	\caption{The matrix domain is discretized with a Cartesian grid with spacings $\Delta x$ and $\Delta y$ in the horizontal and vertical directions, respectively. The fracture grid is obtained by the fracture-matrix grid intersections. The matrix cells cut by the fracture are highlighted in grey.}
	\label{fig:EDFM_grid}
\end{figure}

The problem given by~\eqref{eq:EDFM}--\eqref{eq:EDFM_BCs} is discretized using an embedded approach, meaning that the fracture is allowed to cut arbitrarily the matrix cells, and is thus not forced to lie on the faces of the matrix grid.
For simplicity a structured Cartesian grid is considered for the discretization of the matrix domain and the fracture grid is obtained by intersecting the fracture with the background matrix grid, as shown in Fig.~\ref{fig:EDFM_grid}. 

A standard Two-Point Flux Approximation (TPFA) finite volume scheme is used to discretize fluxes for the problem~\eqref{eq:EDFM}--\eqref{eq:EDFM_BCs}, hence we will have one pressure degree of freedom for each matrix and fracture cell.
The \emph{matrix-matrix} flux $F_{m_1m_2}$ between two generic matrix cells $m_1$ and $m_2$ and the \emph{neighbouring fracture-fracture} flux $F_{f_1f_2}$ between two generic fracture cells $f_1$ and $f_2$ belonging to the same fracture are then written as
\begin{equation*}
	F_{m_1m_2} \vcentcolon = \frac{T_{m_1m_2}}{\mu} (p_{m_1} - p_{m_2}),
	\qquad
	F_{f_1f_2} \vcentcolon = \frac{T_{f_1f_2}}{\mu} (p_{f_1} - p_{f_2}),
\end{equation*}
where $T_{m_1m_2}$ is the \emph{total matrix-matrix transmissibility} and $p_{m_1}$ and $p_{m_2}$ are the pressures of the cells $m_1$ and $m_2$, respectively. Similar explanations apply for the fracture flux and, in particular, $T_{f_1f_2}$ is the \emph{total neighbouring fracture-fracture transmissibility}.

Both transmissibilities are computed by taking the half of the harmonic average of the corresponding \emph{half-transmissibilities}, i.e.
\begin{equation*}
	T_{m_1m_2} = \left( \frac{1}{T_{m_1}} + \frac{1}{T_{m_2}} \right)^{-1} = \frac{T_{m_1} T_{m_2}}{T_{m_1} + T_{m_2}},
	\qquad
	T_{f_1f_2} = \left( \frac{1}{T_{f_1}} + \frac{1}{T_{f_2}} \right)^{-1} = \frac{T_{f_1} T_{f_2}}{T_{f_1} + T_{f_2}},
\end{equation*}
where $T_{m_1}$ and $T_{m_2}$ are the half-transmissibilities of the matrix cells $m_1$ and $m_2$, respectively, while $T_{f_1}$ and $T_{f_2}$ are, respectively, the half-transmissibilities of the fracture cells $f_1$ and $f_2$.

The expressions for the half-transmissibilities in the case of a 2D Cartesian grid and constant aperture $d$ for the fracture are given by
\begin{equation*}
	T_{m_i}^V = 2 k_{m_i,x} \frac{\Delta y}{\Delta x},
	\qquad
	T_{m_i}^H = 2 k_{m_i,y} \frac{\Delta x}{\Delta y},
	\qquad
	T_{f_i} = 2 k_{f_i,\boldsymbol{\tau}} \frac{d}{|f_i|},
	\qquad
	i = 1,2,
\end{equation*}
where $T_{m_i}^V$ and $T_{m_i}^H$ are the $i$-th matrix half-transmissibilities for a vertical and horizontal interface, respectively, $k_{m_i,x}$ is the matrix permeability for the cell $m_i$ in the $x$ direction and $k_{m_i,y}$ in the $y$ direction.
$\Delta x$ and $\Delta y$ are the horizontal and vertical uniform spacings of the Cartesian grid, $k_{f_i,\boldsymbol{\tau}}$ is the tangential fracture permeability for the cell $f_i$ and $|f_i|$ denotes the length of the fracture cell $f_i$.

The full transmissibilities are then given by
\begin{equation*}
	T_{m_1m_2}^V = \widetilde{k}_{m_1m_2,x} \frac{\Delta y}{\Delta x},
	\qquad
	T_{m_1m_2}^H = \widetilde{k}_{m_1m_2,y} \frac{\Delta x}{\Delta y},
	\qquad
	T_{f_1f_2} = \frac{ 2 k_{f_1,\boldsymbol{\tau}} k_{f_2,\boldsymbol{\tau}} }{k_{f_1,\boldsymbol{\tau}} |f_2| + k_{f_2,\boldsymbol{\tau}} |f_1|} d,
\end{equation*}
where $\widetilde{k}_{m_1m_2,x}$ denotes the harmonic average of $k_{m_1,x}$ and $k_{m_2,x}$, while $\widetilde{k}_{m_1m_2,y}$ the harmonic average of $k_{m_1,y}$ and $k_{m_2,y}$. Note that assuming a uniform spacing for the grid simplifies the final transmissibility expressions.

In the following section we describe methods that differ in the way they describe the coupling between the matrix and fracture domains with the aim of considering different types of connections and expressions for the transmissibilities.

\subsection{Embedded Discrete Fracture Model (EDFM)}
The first embedded method, described in \cite{Lee, Li}, was conceived for an efficient handling of long fractures, i.e. those whose length is greater than the characteristic grid cell size, and is known as the \emph{Embedded Discrete Fracture Model} (EDFM). This method works only in presence of highly permeable fractures with respect to the matrix: with this assumption the pressures in the matrix cells cut by the fracture can be taken as continuous and can be described by a single average constant value.
Following \cite{Lee} we then model the flux interaction terms $q_{mf}$ and $q_{fm}$ as one typically does in classical well models, such as that of Peaceman \cite{Peaceman}, and assume a linear pressure profile in the direction normal to the fracture, yielding the following expressions:
\begin{equation} \label{eq:flux_interactions}
	q_{mf} \vcentcolon = \frac{2 \cdot CI_{mf}}{|m|}  \frac{\vec{n}^\top \vec{K}_m \vec{n}}{\mu} (p_m - p_\gamma),
	\qquad
	q_{fm} \vcentcolon = \frac{2 \cdot CI_{mf}}{|f|}  \frac{\vec{n}^\top \vec{K}_m \vec{n}}{\mu} (p_\gamma - p_m),
\end{equation}
where $|m|$ and $|f|$ indicate the measure of the matrix and fracture cell, respectively, while $CI_{mf}$ is the \emph{connectivity index}, which is in analogy with the well productivity indices defined by Peaceman, which depends only on grid parameters as
\begin{equation*}
	CI_{mf} \vcentcolon = \frac{|f|}{\langle d \rangle_{mf}},
\end{equation*}
where $\langle d \rangle_{mf}$ is the \emph{average distance} of the matrix cell $m$ from the fracture cell $f$. $\langle d \rangle_{mf}$  can be computed as
\begin{equation*}
	\langle d \rangle_{mf} = \frac{1}{|m|} \int_{m} |(\vec{x} - \vec{x}_f) \cdot \vec{n}|~ dV,
\end{equation*}
where $\vec{x} \in m$ and $\vec{x}_f \in f$ are generic and $|(\vec{x} - \vec{x}_f) \cdot \vec{n}|$ indicates the distance of a point in the matrix cell from the fracture. Typically the average distance is computed numerically using simple quadrature rules, but for very simple cases analytical expressions are available \cite{Tene, Hajibeygi, Pluimers}.

Clearly, $q_{mf}$ is nonnull only in those matrix cells intersected by the fracture. The \emph{matrix-fracture} fluxes $F_{mf}$ and $F_{fm}$ exchanged between a generic matrix cell $m$ and fracture cell $f$ are given by
\begin{equation} \label{eq:mf_fluxes}
	F_{mf} = \int_{m} q_{mf} dV,
	\qquad
	F_{fm} = \int_{f} q_{fm} dS.
\end{equation}
Using the expressions of the flux interaction terms~\eqref{eq:flux_interactions} to compute the matrix-fracture fluxes~\eqref{eq:mf_fluxes} we get the same result, i.e.
\begin{equation*}
	F_{mf} = \int_{m} q_{mf} dV = - \int_{f} q_{fm} dS.
\end{equation*}
This implies that the local conservation of mass between matrix and fracture is satisfied.
In particular, the expression of the matrix-fracture flux is given by
\begin{equation} \label{eq:Fmf}
	F_{mf} = \frac{T_{mf}}{\mu} (p_m - p_f),
\end{equation}
where $p_m$ and $p_f$ are the pressures of the cells $m$ and $f$, respectively, while $T_{mf}$ is the \emph{matrix-fracture transmissibility}, which can be computed as
\begin{equation} \label{eq:mf_transm}
	T_{mf} = 2 \cdot CI_{mf} \cdot \vec{n}^\top \vec{K}_m \vec{n}.
\end{equation}

The above transmissibility expression stands for a fracture cell $f$ immersed in a matrix cell $m$, hence a single matrix-fracture flux exchange takes place. 
However, it may happen that a fracture coincides with a face of the matrix grid: in this case the fracture cell exchanges flux with the two matrix cells sharing the face, and the corresponding transmissibilities are halved with respect to~\eqref{eq:mf_transm}.
Other expressions for the matrix-fracture transmissibility have been proposed in the literature, e.g. in \cite{Tene, Jiang}, but the results are very close to those provided by~\eqref{eq:mf_transm} since the fracture is assumed to be highly permeable with respect to the matrix.

\subsection{Projection-based Embedded Discrete Fracture Model (pEDFM)} \label{ssec:pEDFM}
\begin{figure}[!t]
	\centering
	\includegraphics[width=.75\textwidth]{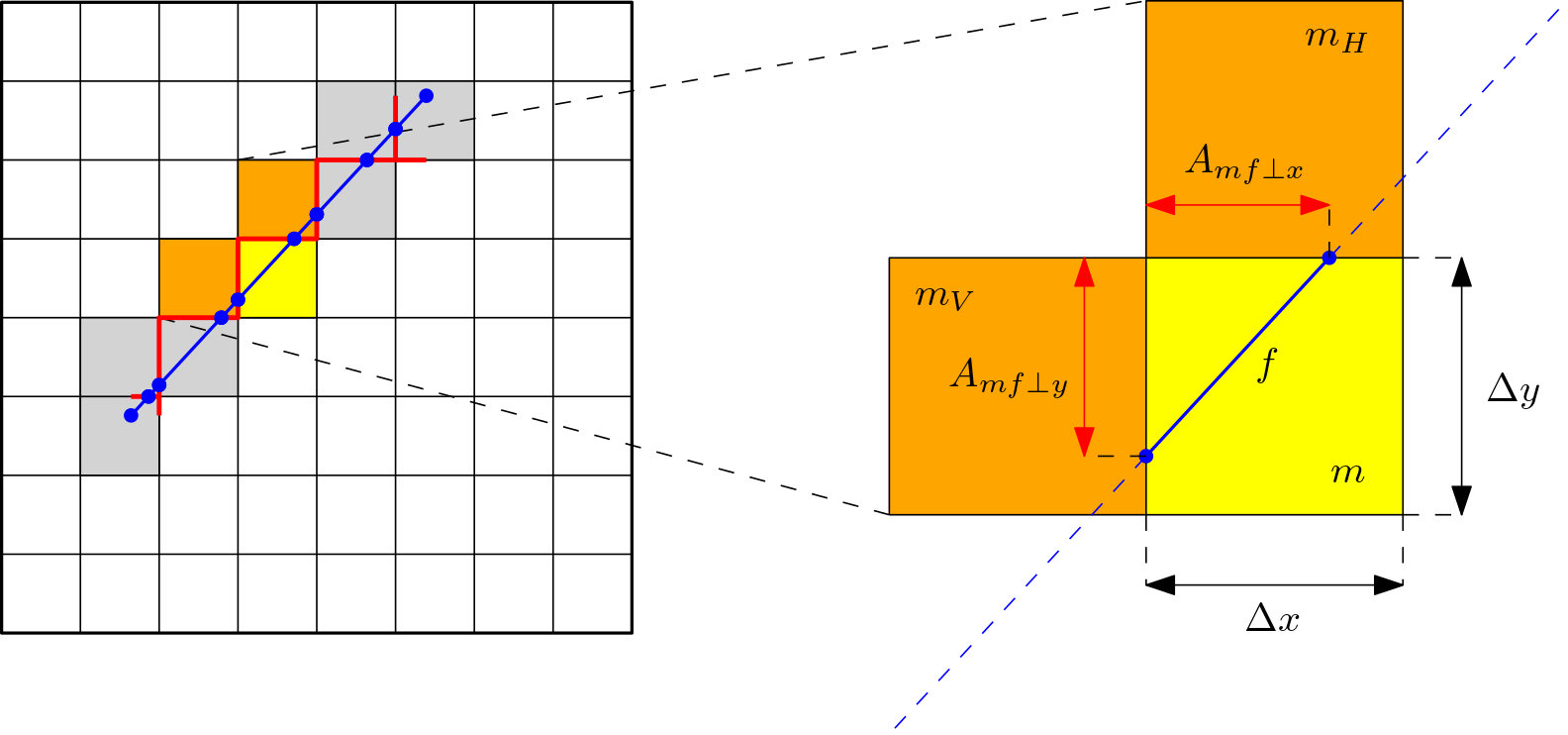}
	\caption{pEDFM for a 2D Cartesian grid. The continuous fracture projection path is highlighted in red. The fracture cell $f$ is directly connected to the yellow matrix cell $m$ and the non-neighbouring matrix cells are highlighted in orange. The fracture cell projections are also put in evidence.}
	\label{fig:pEDFM}
\end{figure} 

The \emph{Projection-based Embedded Discrete Fracture Model} (pEDFM), introduced in \cite{Tene}, is an extension of the EDFM method that allows to take into account, unlike the original embedded formulation, even the case of impermeable fractures.
To this aim, matrix-fracture and matrix-matrix transmissibilities near to the fractures are modified properly.

As a preliminary step, a set of matrix grid faces is selected such that they define a continuous projection path of the fractures on the matrix domain, as highlighted in red on the left side of Fig.~\ref{fig:pEDFM}.
As stated in \cite{Tene}, and as we will see in Section~\ref{ssec:conv_tests} and Appendix~\ref{sec:warn_pEDFM}, it is important to ensure the continuity of the fracture projection paths to obtain a correct representation of the pressure fields, especially for fractures acting as barriers for the fluid flow. This is difficult to achieve in the 3D case and for general, non Cartesian grids.
Then, let us consider a fracture cell $f$ intersecting a matrix cell $m$.
We denote with $m_H$ and $m_V$ the matrix cells sharing with $m$ the faces on which the fracture cell is projected along the $x$ and $y$ directions, respectively.
They are called \emph{non-neighbouring} matrix cells since they are not cut by the fracture, but still they are affected by the modification of the transmissibilities. These cells are depicted in orange in Fig.~\ref{fig:pEDFM}.
We also indicate with $A_{mf \perp x}$ and $A_{mf \perp y}$ the areas of the projections of the fracture cell $f$ along $x$ and $y$, respectively.
Under the assumption of isotropic and diagonal permeability tensors, we define the following transmissibilities
\begin{equation} \label{eq:Tmf_jiang}
	T_{mf} =  \frac{T_{m} T_{f}}{T_{m} + T_{f}},
\end{equation}
\begin{equation} \label{eq:Tmnf_jiang}
	T_{m_Hf} = \frac{T_{m_H} T_{f}}{T_{m_H} + T_{f}},
	\qquad
	T_{m_Vf} =  \frac{T_{m_V} T_{f}}{T_{m_V} + T_{f}}.
\end{equation}
\begin{equation} \label{eq:Tmm_mod}
	T_{mm_H} = \frac{\Delta x - A_{mf \perp x}}{\Delta y} \widetilde{k}_{mm_H}^{DW},
	\qquad
	T_{mm_V} = \frac{\Delta y - A_{mf \perp y}}{\Delta x} \widetilde{k}_{mm_V}^{DW},
\end{equation}
where $T_{m_Hf}$ and $T_{m_Vf}$ are the new \emph{non-neighbouring matrix-fracture transmissibilities}, while $T_{mm_H}$ and $T_{mm_V}$ are modified matrix-matrix transmissibilities. Both of them are needed to take into account the isolating effect of impermeable fractures. 
$T_m$, $T_{m_H}$, $T_{m_V}$, $T_f$ are all half-transmissibilities, whose expressions are given by
\begin{equation} \label{eq:half_transm_jiang}
	T_m = CI_{mf} k_m,
	\qquad
	T_{m_H} = \frac{A_{mf \perp x}}{d_{m_Hf}} k_{m_H},
	\qquad
	T_{m_V} = \frac{A_{mf \perp y}}{d_{m_Vf}} k_{m_V},
	\qquad
	T_f = \frac{2 |f|}{d} k_f,
\end{equation}
where $d_{m_Hf}$ and $d_{m_Vf}$ are the distances from the centroids of $m_H$ and $m_V$, respectively, to the centroid of the fracture cell $f$.
Note that the permeabilities appearing in~\eqref{eq:Tmm_mod} are all distance-weighted harmonically averaged between neighbouring cells.

We recall that the \emph{distance-weighted harmonic average} between two generic cells is defined as follows:
\begin{equation*}
	\widetilde{k}_{12}^{DW} \vcentcolon = \frac{k_1 k_2}{k_1 d_2 + k_2 d_1} (d_1 + d_2), 
\end{equation*}
where $k_1$ and $k_2$ are the permeabilities of the cells, while $d_1$ and $d_2$ the distance weights pertaining to the cells.
In particular, for $\widetilde{k}_{mm_H}^{DW}$ both weights correspond to the distance from the centre of $m$ to the interface shared by $m$ and $m_H$, i.e. $\Delta y / 2$. Following the same reasoning, for $\widetilde{k}_{mm_V}^{DW}$ both distance weights are equal to $\Delta x / 2$.

The transmissibility expressions reported above are taken from~\cite{Jiang} and, apart from the modified matrix-matrix transmissibility formulae~\eqref{eq:Tmm_mod}, they present small differences with respect to those presented in the original pEDFM paper~\cite{Tene}. In particular, the differences become more and more evident as the distances of the non-neighbouring matrix cells from the fracture cells approach the fracture aperture value.

For the special case in which the fracture lies on the matrix grid faces the projection-based method is equivalent to a conforming approach, while this is not the case for the EDFM.
The computation of the fracture cells projections is definitely the most complicated part of the method from the implementation point of view, especially in three-dimensional cases, in presence of grids having cells of generic shape, and when fractures (which could also be curved in the case of faults) are arbitrarily oriented in space.
  
\section{LEDFM Formulation} \label{sec:LEDFMform}

The main idea behind the \emph{Local Embedded Discrete Fracture Model} (LEDFM) is that of using local transmissibility upscaling methods to define new transmissibilities for the embedded formulation of the Darcy problem~\eqref{eq:EDFM}--\eqref{eq:EDFM_BCs}. In particular, two types of local problems are considered:
\begin{itemize}
	\item Matrix-Fracture (M-F) local problems, to define new matrix-fracture transmissibilities;
	\item Matrix-Matrix (M-M) local problems, to define new near-fracture matrix-matrix transmissibilities, i.e. those relative to the faces of the matrix grid having at least one of the neighbour cells cut by a fracture.
\end{itemize}
In this way, no assumptions are made about the near-fracture pressure distribution to compute the transmissibilities between matrix and fracture cells, but they are instead directly computed post-processing the corresponding local problem solution, allowing for improved accuracy at fracture tips and intersections.
Moreover, blocking fractures are accounted for without having to project them on porous matrix interfaces, which is a necessary step in pEDFM that requires a lot of effort to implement in practice, especially for complex grids and the three-dimensional case.

For the sake of simplicity, and also for a faster implementation of the subsequent multiscale version of the method, we focus on the special case of square matrix grid cells.
At the moment the method is constructed so that only one fracture is allowed to cut each local domain. Moreover, the case where the fracture lies on the local domain boundary, i.e. on coarse faces, is not considered.
Clearly, fractures that are completely immersed in one coarse cell are not considered since they should be represented by means of homogenization at the scale of interest.
LEDFM does not adopt a multiple sub-region approach, contrarily to some upscaling methods for fracture networks that can be found in the literature, e.g. \cite{Karimi-Fard, Fumagalli2}, thus only one upscaled transmissibility for each type of connection within a local domain is computed.

In particular, the LEDFM formulation employs the local flow based upscaling of the transmissibility, whose description can be found, e.g., in \cite{Durlofsky2} along with a detailed review of the most important flow based methods for the upscaling of porous media properties.
In local upscaling methods the transmissibility for a coarse grid face is computed using local flow solutions on a fine grid region that includes the coarse cells neighbouring the face. This region could be extended to include more coarse cells (extended local methods, see \cite{Holden, Nielsen, Holden2, Wen2}). These methods are generally easy to implement and their computational cost is low compared to other flow based methods. However, they are heavily dependent on the boundary conditions chosen to solve local flow problems, which is the major drawback of these methods.

This section is organized in the following way: the governing equations of both types of local problems are described in Subsections~\ref{ssec:M-F} and~\ref{ssec:M-M}, then in Subsection~\ref{ssec:discretization} we give further details on the discretization of the local problems and in Subsection~\ref{ssec:workflow} we point out when the local problems should be solved and how to use the computed local solutions to construct the new LEDFM method.
Finally, in Subsection~\ref{ssec:MSFV} the multiscale modification of the local method is described.

\subsection{Matrix-Fracture (M-F) Local Problem} \label{ssec:M-F}

\begin{figure}[!b]
	\centering
	\includegraphics[width=.35\textwidth]{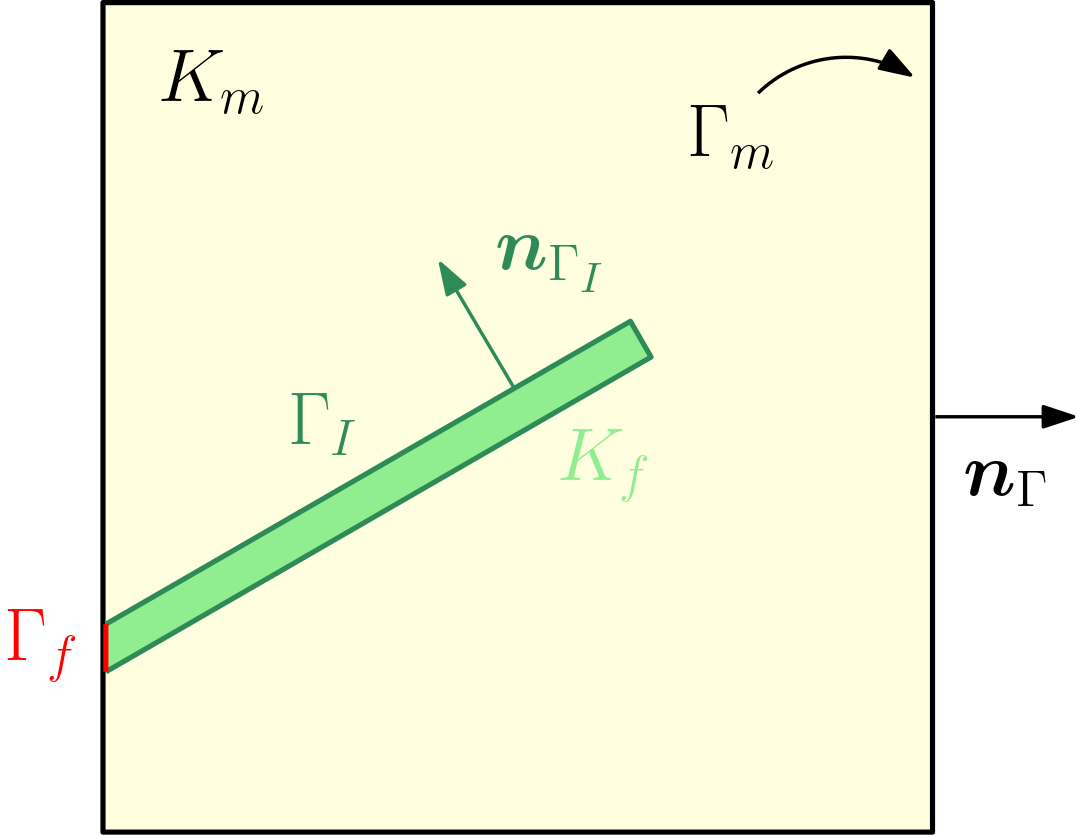}
	\caption{Local domain example for the M-F problem.}
	\label{fig:M_F_LocalProblemDef}
\end{figure} 

Consider the domain depicted in Fig~\ref{fig:M_F_LocalProblemDef}.
Let $\Omega$ be the whole domain and $K \in \Omega$ a grid cell cut by a fracture. The solution of the M-F local problem, taking inspiration from \cite{Karimi-Fard, Fumagalli2}, is chosen to be that of a slightly-compressible single-phase flow equation with impermeable boundaries, in the special case where the porosity $\Phi$ is considered to be constant
\begin{equation} \label{eq:MF_problem}
	\begin{cases}
		c \Phi \dfrac{\partial p}{\partial t} -\nabla\cdot\left(\dfrac{\vec{K}}{\mu}\nabla p\right) = q_f & \text { in }  K  \text { for }  t>0, \\
		-\left(\dfrac{\vec{K}}{\mu}\nabla p\right) \cdot \vec{n}_\Gamma = 0 & \text { on }  \partial K  \text { for }  t>0, \\
		p=0 & \text { in }  K \text { for } t=0.
	\end{cases}
	\quad \ \ \
	q_f = 
	\begin{cases}
		\overline{q}_{f} & \text { in }  K_f \\
		0 & \text { in }  K_m
	\end{cases}
	\quad
	\overline{q}_{f} > 0,
\end{equation}
where $c$ is the fluid compressibility, $q_f$ a piecewise constant source term, strictly positive inside the fracture and zero elsewhere, and $\vec{n}_\Gamma$ the unit outward normal vector to the grid cell boundary $\Gamma \vcentcolon = \partial K$.
Note also that $K = K_m \cup K_f$ and $\Gamma = \Gamma_m \cup \Gamma_f$, where $K_m$ and $K_f$ are the matrix and fracture cell domains, respectively, while $\Gamma_m$ and $\Gamma_f$ are the parts of their boundaries lying on $\partial K$.

Note that the definition of the problem can be easily adapted to the case of a fracture cell cutting the entire local domain from side to side.

The problem describes a continuous injection of fluid in the fracture, so that the overall pressure in the domain increases with time and, after a transient period, a \emph{pseudo steady state} condition is reached, meaning that the time derivative of pressure reaches a constant value, i.e.
\begin{equation*}
	\frac{\partial p}{\partial t} = a \qquad a \in \mathbb{R}.
\end{equation*}
Following the upscaling procedure described in \cite{Karimi-Fard, Fumagalli2}, problem~\eqref{eq:MF_problem} should be solved until the pseudo steady state condition is reached.
It can be shown that the problem is weakly coercive and the numerical errors accumulated in time may largely affect the accuracy of the corresponding numerical solution.

However, finding the pseudo steady state solution of the previous problem can be proved to be equivalent to finding the solution of a stationary problem with a modified source term that depends on the previous source term $\overline{q}_f$, fracture porosity $\Phi_f$, compressibility $c_f$ and pseudo steady state constant $a_f$, as well as the measures of the matrix and fracture domains.
We now give the proof that leads to a slightly different problem from that defined in \cite{Fumagalli2}.

Let us write problem~\eqref{eq:MF_problem} separately for the matrix and fracture domains
\begin{equation*}
	\begin{cases}
		c_m \Phi_m \dfrac{\partial p_m}{\partial t} -\nabla\cdot\left(\dfrac{\vec{K}_m}{\mu}\nabla p_m\right) = 0 & \text { in }  K_m  \text { for }  t>0, \\
		-\left(\dfrac{\vec{K}_m}{\mu}\nabla p_m\right) \cdot \vec{n}_{\Gamma_m} = 0 & \text { on }  \Gamma_m  \text { for }  t>0, \\
		p_m=0 & \text { in }  K_m \text { for } t=0,
	\end{cases}
\end{equation*}
\begin{equation*}
	\begin{cases}
		c_f \Phi_f \dfrac{\partial p_f}{\partial t} -\nabla\cdot\left(\dfrac{\vec{K}_f}{\mu}\nabla p_f\right) = \overline{q}_{f} & \text { in }  K_f  \text { for }  t>0, \\
		-\left(\dfrac{\vec{K}_f}{\mu}\nabla p_f\right) \cdot \vec{n}_{\Gamma_f} = 0 & \text { on }  \Gamma_f  \text { for }  t>0, \\
		p_f=0 & \text { in }  K_f \text { for } t=0,
	\end{cases}
\end{equation*}
where the subscripts $m$ and $f$ refer to porous matrix and fracture quantities, respectively.
Note also that $\vec{n}_{\Gamma_m}$ and $\vec{n}_{\Gamma_f}$ point in the outward direction with respect to $K_m$ and $K_f$, respectively.
We denote with $\Gamma_I$ the matrix-fracture interface and with $\vec{n}_{\Gamma_I}$ the unit normal vector to $\Gamma_I$ directed from $K_f$ to $K_m$.
Pressure and normal flux continuity are enforced on $\Gamma_I$
\begin{equation*}
	\begin{cases}
		p_m = p_f
		\\
		-\left(\dfrac{\vec{K}_m}{\mu}\nabla p_m\right) \cdot \vec{n}_{\Gamma_I} = -\left(\dfrac{\vec{K}_f}{\mu}\nabla p_f\right) \cdot \vec{n}_{\Gamma_I}
	\end{cases}
	\text { on } \Gamma_I \text { for }  t>0.
\end{equation*}

Integrating the matrix and fracture governing equations over their respective domains $K_m$ and $K_f$, and exploiting Gauss's theorem for the divergence term we obtain
\begin{equation} \label{eq:MF_intm}
	\int_{K_m} c_m \Phi_m \dfrac{\partial p_m}{\partial t}
	-
	\int_{\partial K_m} \left(\dfrac{\vec{K}_m}{\mu}\nabla p_m\right) \cdot \vec{n}_{\partial K_m}
	=
	0,
\end{equation}
\begin{equation}  \label{eq:MF_intf}
	\int_{K_f} c_f \Phi_f \dfrac{\partial p_f}{\partial t}
	-
	\int_{\partial K_f} \left(\dfrac{\vec{K}_f}{\mu}\nabla p_f\right) \cdot \vec{n}_{\partial K_f}
	=
	\int_{K_f} \overline{q}_f,
\end{equation}
where $\vec{n}_{\partial K_m}$ and $\vec{n}_{\partial K_f}$ are unit normal vectors pointing outwards with respect to $\partial K_m$ and $\partial K_f$, respectively. 
Since $\partial K_m = \Gamma_m \cup \Gamma_I$ and $\partial K_f = \Gamma_f \cup \Gamma_I$, equations~\eqref{eq:MF_intm} and~\eqref{eq:MF_intf} can be rewritten as 
\begin{equation} \label{eq:MF_intm2}
	\int_{K_m} c_m \Phi_m \dfrac{\partial p_m}{\partial t}
	-
	\int_{\Gamma_m} \left(\dfrac{\vec{K}_m}{\mu}\nabla p_m\right) \cdot \vec{n}_{\Gamma_m}
	+
	\int_{\Gamma_I} \left(\dfrac{\vec{K}_m}{\mu}\nabla p_m\right) \cdot \vec{n}_{\Gamma_I}
	=
	0,
\end{equation}
\begin{equation} \label{eq:MF_intf2}
	\int_{K_f} c_f \Phi_f \dfrac{\partial p_f}{\partial t}
	-
	\int_{\Gamma_f} \left(\dfrac{\vec{K}_f}{\mu}\nabla p_f\right) \cdot \vec{n}_{\Gamma_f}
	-
	\int_{\Gamma_I} \left(\dfrac{\vec{K}_f}{\mu}\nabla p_f\right) \cdot \vec{n}_{\Gamma_I}
	=
	\int_{K_f} \overline{q}_f,
\end{equation}
where the matrix-fracture interface term in equation~\eqref{eq:MF_intm2} changes sign due to the definitions of the unit normal vectors.
Using the no flow boundary conditions, the flux continuity on $\Gamma_I$ and then summing equations~\eqref{eq:MF_intm2} and~\eqref{eq:MF_intf2} the boundary terms vanish and we get
\begin{equation} \label{eq:MF_eq_proof}
	\int_{K_m} c_m \Phi_m \dfrac{\partial p_m}{\partial t}
	+
	\int_{K_f} c_f \Phi_f \dfrac{\partial p_f}{\partial t}
	=
	\int_{K_f} \overline{q}_f.
\end{equation}
When the system reaches a pseudo steady state condition, we have 
\begin{equation*}
	\frac{\partial p_m}{\partial t} = a_m 
	\qquad
	\frac{\partial p_f}{\partial t} = a_f
	\qquad
	a_m, a_f \in \mathbb{R}.
\end{equation*}
Since we are considering slightly-compressible flow, both compressibilities $c_m$ and $c_f$ are constant and small. Moreover, we assume that matrix and fracture porosities are homogeneous in $K_m$ and $K_f$, respectively. Under these hypotheses, and knowing that $\overline{q}_f$ is constant, equation~\eqref{eq:MF_eq_proof} can be rewritten in the following way
\begin{equation*}
	|K_m| c_m \Phi_m a_m + 	|K_f| c_f \Phi_f a_f = 	|K_f| \overline{q}_f.
\end{equation*}
We are now able to reformulate the pseudo steady state as the solution of a modified stationary problem:
\begin{equation} \label{eq:MF_problem_stationary}
	\begin{cases}
		-\nabla\cdot\left(\dfrac{\vec{K}}{\mu}\nabla p\right) = q_{f,s} & \text { in }  K, \\
		-\left(\dfrac{\vec{K}}{\mu}\nabla p\right) \cdot \vec{n}_\Gamma = 0 & \text { on }  \partial K, 
	\end{cases}
	\quad \ \ \
	q_{f,s} = 
	\begin{cases}
		\overline{q}_{f,s} & \text { in }  K_f \\
		-\dfrac{|K_f|}{|K_m|} \overline{q}_{f,s} & \text { in }  K_m
	\end{cases}
	,
\end{equation} 
where $\overline{q}_{f,s} = \overline{q}_f - c_f \Phi_f a_f$. 
Since only homogeneous Neumann boundary conditions are considered for the previous system, we also need to fix a condition to ensure the uniqueness of the pressure, for instance null average:
\begin{equation*}
	\frac{1}{|K|} \int_K p = 0.
\end{equation*}

Clearly, the stationary formulation \eqref{eq:MF_problem_stationary} is preferred over the original for reasons of computational time.

It is important to observe that, when the pseudo steady state condition is reached, the shapes of the isopressure curves do not change any more (even though the associated pressure values change with time), and they depend only on fracture geometry and permeabilities of both domains \cite{Karimi-Fard}.
For this reason the source term $\overline{q}_{f,s}$ is taken to be equal to $\overline{q}_{f,s} = 1$, for simplicity.

\subsection{Matrix-Matrix (M-M) Local Problem} \label{ssec:M-M}

\begin{figure}[!t]
	\centering
	\includegraphics[width=.55\textwidth]{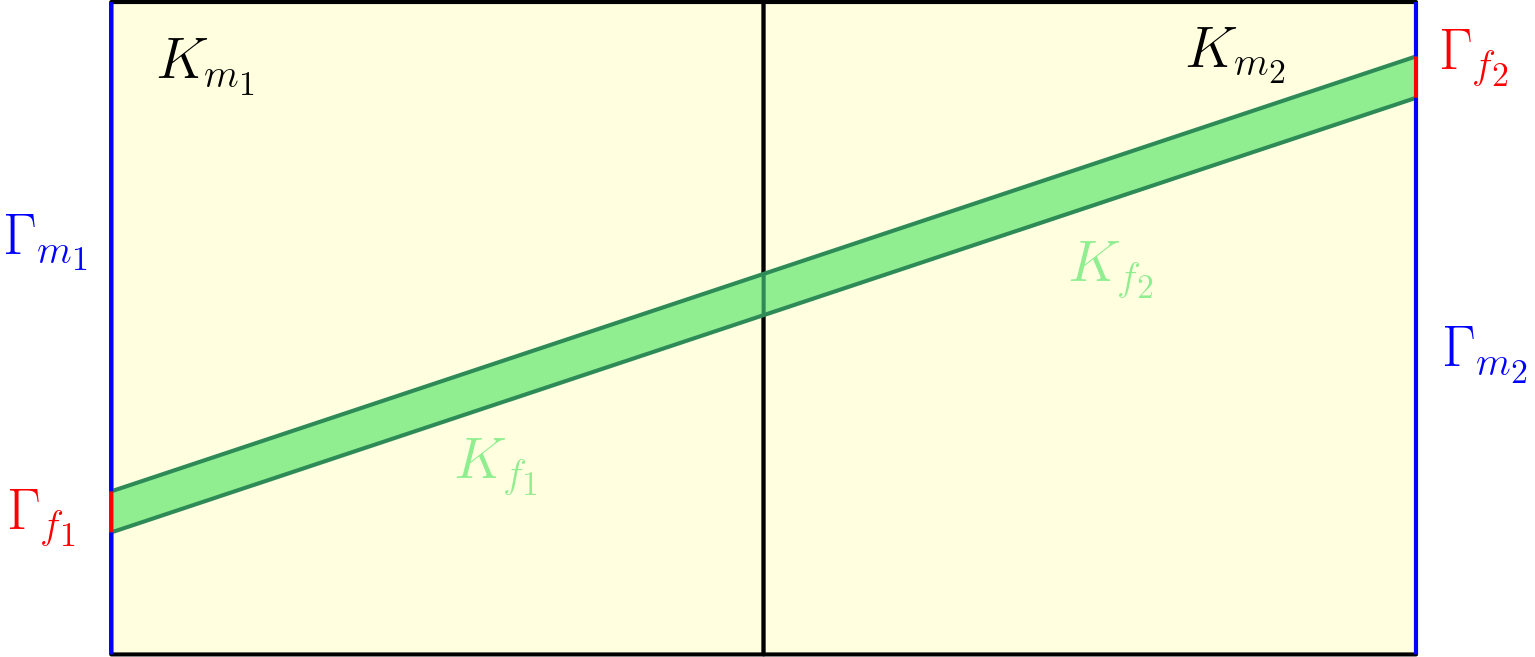}
	\caption{Local domain example for the M-M problem.}
	\label{fig:M_M_LocalProblemDef}
\end{figure} 

Let us now consider the system depicted in Fig.~\ref{fig:M_M_LocalProblemDef}. The solution of the M-M local problem for two generic grid cells $K_1, K_2 \in \Omega$ neighbouring a face and where at least one of the cells is cut by a fracture is chosen to be that of the usual incompressible single-phase Darcy problem.
With respect to traditional transmissibility upscaling, see e.g.~\cite{Durlofsky2}, here a different set of boundary conditions is employed:
\begin{equation} \label{eq:MM_prob}
	\begin{cases}
		\nabla\cdot\left(\dfrac{\vec{K}}{\mu}\nabla p\right) = 0 & \text { in } D, \\
		p=1 & \text { on }  \Gamma_{m_1}, \\
		-\left(\dfrac{\vec{K}}{\mu}\nabla p\right) \cdot \vec{n}=1 & \text { on }  \Gamma_{m_2}, \\
		-\left(\dfrac{\vec{K}}{\mu}\nabla p\right) \cdot \vec{n}= 0 & \text { on } \partial D \setminus (\Gamma_{m_1} \cup \Gamma_{m_2}).
	\end{cases}
\end{equation}
Here, we have that $D = K_1 \cup K_2$, $K_1 = K_{m_1} \cup K_{f_1}$ and $K_2 = K_{m_2} \cup K_{f_2}$.
$K_{m_1}$ and $K_{m_2}$ are the matrix cell domains while $K_{f_1}$ and $K_{f_2}$ the fracture cell ones. $\Gamma_{m_1}$ and $\Gamma_{m_2}$ represent the left and right domain boundaries pertaining to the matrix part, while $\vec{n}$ denotes the unit outward normal vector.

The problem is here described for the case of a fracture cutting the entire local domain $D$ from side to side and for a vertical interface, but its definition can be adapted with ease both to the case of a fracture partially immersed in the local domain and to that of a horizontal interface.

\begin{figure}[!t]
	\centering
	\begin{subfigure}[b]{0.495\textwidth}
		\centering
		\includegraphics[width=\textwidth]{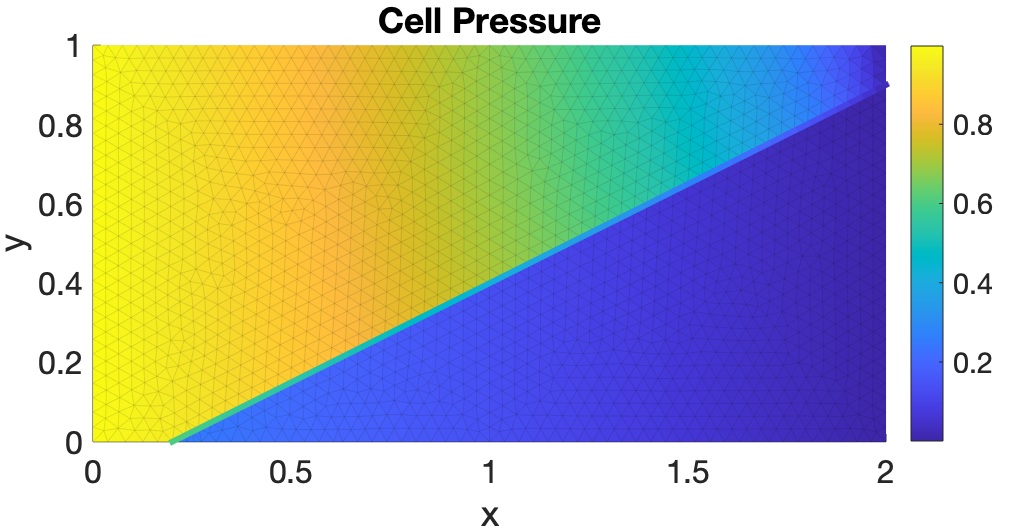}
		\caption{}
		\label{fig:pressure_pressure}
	\end{subfigure}
	\hfill
	\begin{subfigure}[b]{0.495\textwidth}
		\centering
		\includegraphics[width=\textwidth]{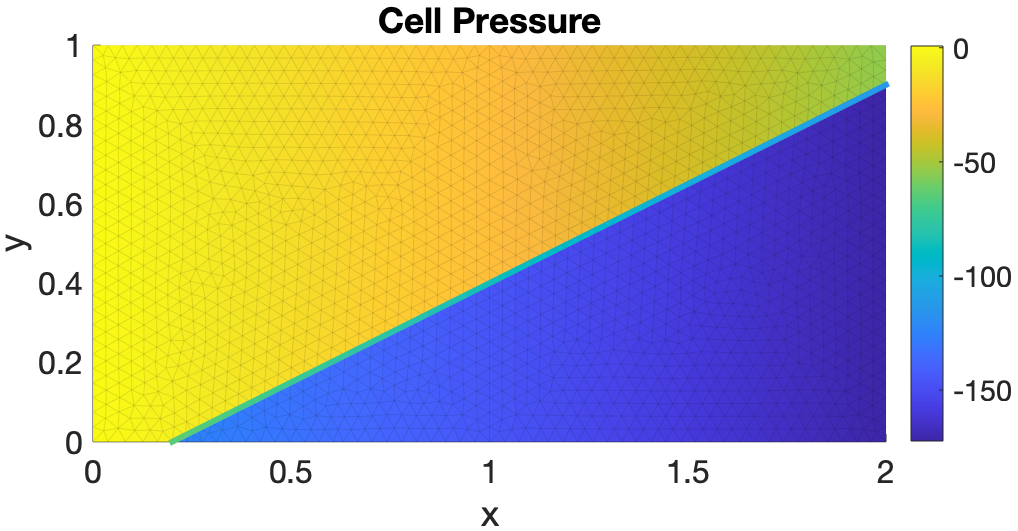}
		\caption{}
		\label{fig:pressure_flux}
	\end{subfigure}
	\caption{M-M local problem solutions for an impermeable fracture $(k_f/k_m = 10^{-4}, d = 10^{-3})$. (a) Pressure-pressure boundary conditions ($p=1$ on $\Gamma_{m_1}$, $p=0$ on $\Gamma_{m_2}$). (b) Pressure-flux boundary conditions as in~\eqref{eq:MM_prob}. Note that for the purpose of transmissibility calculations one should look at pressure differences, not at absolute pressure values.}
	\label{fig:M_M_BCs}
\end{figure}

Pressure-flux conditions are imposed on $\Gamma_{m_1}$ and $\Gamma_{m_2}$, respectively, and the remaining part of the boundary is set to be impermeable: this non-standard choice of the boundary conditions stems from the observation that they are capable of capturing effectively the barrier effect of fractures as opposed to the standard pressure-pressure conditions as shown in Fig.~\ref{fig:M_M_BCs}. Indeed, it can be seen that the pressure of the matrix cells neighbouring the right domain boundary reflects more appropriately the barrier effect of the fracture when pressure-flux conditions are used.

\subsection{Discretization} \label{ssec:discretization}
The local problems~\eqref{eq:MF_problem_stationary} and~\eqref{eq:MM_prob} are solved with a conforming Discrete Fracture-Matrix (DFM) model on triangular matrix grids.
A Control Volume Finite Difference (CVFD) method is adopted with the possibility of choosing between the standard TPFA scheme and the Multi-Point Flux Approximation (MPFA) scheme, the latter being capable of providing a consistent discretization for anisotropic permeabilities.
In particular, a hybrid formulation is adopted, where the fractures are lower-dimensional with respect to the matrix grid, but equi-dimensional for the purpose of transmissibility computations.
We refer to \cite{Karimi-Fard2, Sandve} for more details on the discretization procedure.

Since we are considering triangular matrix grids, the $\vec{K}$-orthogonality condition is not satisfied in general, even with isotropic $\vec{K}$, unless the triangulation is of Delaunay type.
However, even if the TPFA scheme is consistent only for $\vec{K}$-orthogonal grids, we still used it most of the times, since we observed that, in this case, the computed upscaled transmissibilities and the corresponding coarse scale EDFM solution were very similar to those obtained using a MPFA scheme, and this allowed to spare some computational resources. 

Note that in the case of the M-M local problem the grid should also honour the interface shared by the two matrix coarse blocks to allow for an easy computation of the coarse scale flux, which is involved in the definition of the new near-fracture matrix-matrix transmissibility, as it will be evident in Subsection~\ref{ssec:workflow}.

\subsection{Workflow} \label{ssec:workflow}
\begin{figure}[!t]
	\centering
	\begin{subfigure}[b]{0.495\textwidth}
		\centering
		\includegraphics[width=.6\textwidth]{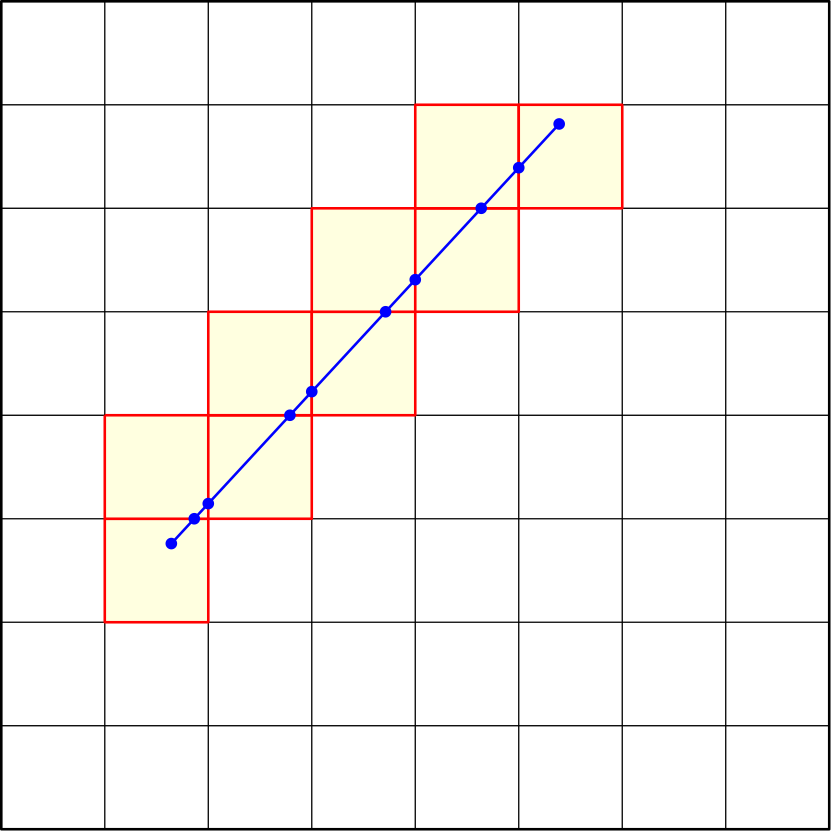}
		\caption{}
		\label{fig:LEDFM_problems1}
	\end{subfigure}
	\hfill
	\begin{subfigure}[b]{0.495\textwidth}
		\centering
		\includegraphics[width=.6\textwidth]{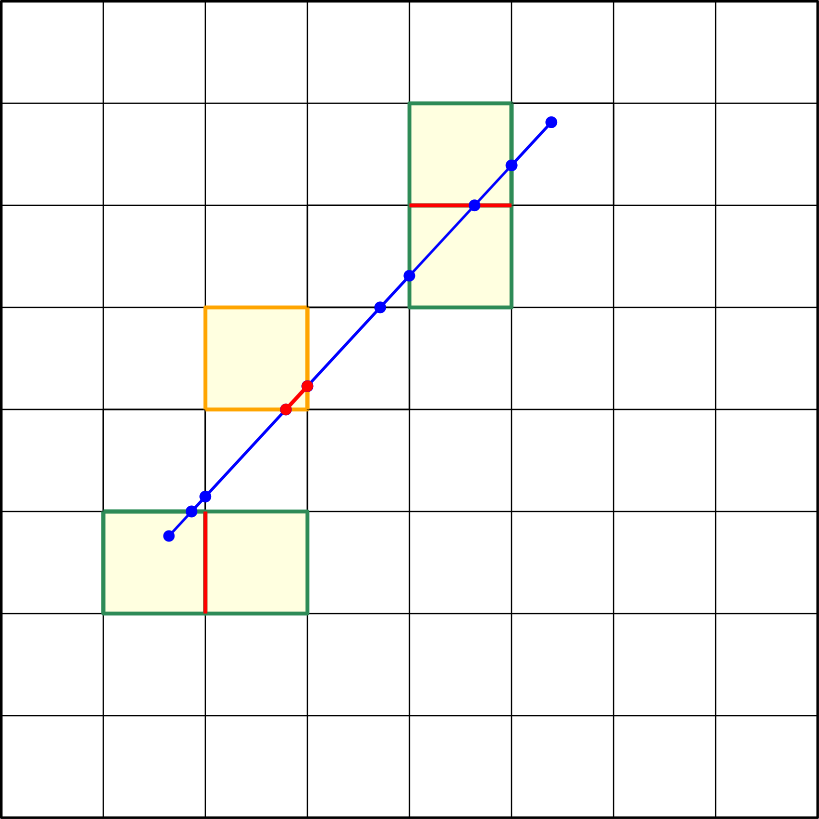}
		\caption{}
		\label{fig:LEDFM_problems2}
	\end{subfigure}
	\caption{LEDFM workflow. (a) It should be solved: a M-F local problem for each yellow cell, a M-M local problem for each red cell edge. (b) Examples of M-F (orange boundary) and M-M (green boundary) local problems domains.}
	\label{fig:LEDFM_problems}
\end{figure}

Let us consider the fractured system illustrated in Fig~\ref{fig:LEDFM_problems}. The LEDFM method requires to solve:
\begin{itemize}
	\item a M-F local problem whenever a matrix cell is cut by a fracture;
	\item a M-M local problem whenever a face of the matrix grid has at least one of the
	neighbour cells cut by a fracture.
\end{itemize}
Note that, before solving the local problems, the local domains are transformed to a unit square for M-F problems and to a $2\times1$ rectangle for M-M problems in order to standardize the way local problems are solved.
Clearly, also the fractures cutting the local domains should be transformed so as to reproduce the shape of the original local domain geometry.
In particular, since the local domain transformations mentioned above include a scaling transformation, also the aperture of the fracture cell within the local domain should be scaled accordingly.
These geometrical transformations are introduced to allow also for a simpler implementation of the routines responsible of the computation of local solutions.

For the M-F problem, let us now denote with $\mathcal{T}_{mf} = \mathcal{T}_{m} \cup \mathcal{T}_{f}$ the local fine grid, where $\mathcal{T}_{m}$ is the fine triangle grid of the matrix coarse cell $K_{m}$ and $\mathcal{T}_{f}$ the fine grid of the coarse fracture cell $K_{f}$.
$\mathcal{T}_{m}$ honours the fracture geometry, and we indicate with $\mathcal{F}_{mf}$ the set of fine scale faces $f$ forming the matrix-fracture interface.

Similarly, for the M-M problem we denote with $\mathcal{T}_{mm} =  \mathcal{T}_{m_1} \cup \mathcal{T}_{m_2} \cup \mathcal{T}_{f_1} \cup  \mathcal{T}_{f_2}$ the corresponding local fine grid, where $\mathcal{T}_{m_1}$ and $\mathcal{T}_{m_2}$ are the fine triangle grids of the matrix coarse cells $K_{m_1}$ and $K_{m_2}$, respectively, while $\mathcal{T}_{f_1}$ and $\mathcal{T}_{f_2}$ the fine grids of the coarse fracture cells $K_{f_1}$ and $K_{f_2}$, respectively.
Both the coarse face $e_{m_1|m_2}$ shared by the matrix coarse blocks and the fracture should be honoured by the fine discretization, and we indicate with $\mathcal{F}_{m_1m_2}$ the set of fine scale faces $f$ forming $e_{m_1|m_2}$.

Once that all the M-F and M-M local problems have been solved, the corresponding new upscaled transmissibilities are computed in the following way
\begin{equation} \label{eq:LEDFM_transm}
	T_{mf} =  \mu \left| \frac{F_{mf}}{\overline{p}_{m} - \overline{p}_{f} } \right|
	\qquad
	T_{m_1m_2} = \mu \left| \frac{F_{m_1m_2}}{\overline{p}_{m_1} - \overline{p}_{m_2} } \right|,
\end{equation}
where $F_{mf}$ is the total flux exchanged through the matrix-fracture interface, while $F_{m_1m_2}$ is the total flux through the face $e_{m_1|m_2}$ of the coarse grid, and they are computed as the sum of fine scale fluxes
\begin{equation*}
	F_{mf} = \sum_{f \in \mathcal{F}_{mf}} T_{l_ml_f} (p_{l_m} - p_{l_f})
	\qquad
	F_{m_1m_2} = \sum_{f \in \mathcal{F}_{m_1m_2}} T_{l_1l_2} (p_{l_1} - p_{l_2}),
\end{equation*}
where $l_m \in \mathcal{T}_{m}$, $l_f \in \mathcal{T}_{f}$ are fine scale cells and $p_{l_m}$, $p_{l_f}$ are their corresponding pressure values.
$T_{l_ml_f}$, instead, denotes the fine scale transmissibility of the face $f$ shared by the fine scale cells $l_m$ and $l_f$.

Finally, $\overline{p}_{m}$, $\overline{p}_{f}$, $\overline{p}_{m_1}$ and $\overline{p}_{m_2}$ are the volume weighted averages of the pressure over the coarse cells $K_{m}$, $K_{f}$,  $K_{m_1}$ and $K_{m_2}$, respectively, defined as
\begin{equation*}
	\overline{p}_{m} = \frac{\sum_{l \in \mathcal{T}_{m}} |l| p_l  }{|K_{m}|}
	\qquad
	\overline{p}_{f} = \frac{\sum_{l \in \mathcal{T}_{f}} |l| p_l  }{|K_{f}|}
	\qquad
	\overline{p}_{m_1} = \frac{\sum_{l \in \mathcal{T}_{m_1}} |l| p_l  }{|K_{m_1}|}
	\qquad
	\overline{p}_{m_2} = \frac{\sum_{l \in \mathcal{T}_{m_2}} |l| p_l  }{|K_{m_2}|},
\end{equation*}
where $p_l$ is the pressure of the fine scale cell $l$.

\subsection{Multiscale Modification} \label{ssec:MSFV}
The expression of the upscaled matrix-matrix transmissibility $T_{m_1m_2}$ in~\eqref{eq:LEDFM_transm} comes from a two-point flux approximation assumed at the coarse grid level.
As mentioned before, TPFA is accurate only for $\vec{K}$-orthogonal grids, for instance for Cartesian grids with isotropic $\vec{K}$.
However, fractures can be seen as objects that introduce full-tensor effects since they give preferential paths to the flow.
If one wants to accurately capture these effects, the transmissibility upscaling procedure for the near-fracture matrix-matrix connections needs to be reconsidered in a multi-point flux approximation framework.
To do that, one of the possible approaches is that of replacing the M-M local problems with the Multiscale Finite Volume Method (MSFV).
For Cartesian grids, this method constructs 9-point stencil modifications of the coarse grid cells whose coefficients are computed post-processing local numerical fine scale solutions of elliptic problems in the interaction regions of the dual grid.
A more detailed explanation of the method is given in \cite{Jenny, Hesse}. 

In our case, the 9-point stencil modification of the MSFV method should be applied only on a subset of the matrix cells of the coarse grid, i.e. those near to the fracture, so that the increase in computational cost is limited.

\begin{figure}[!t]
	\centering
	\includegraphics[width=.35\textwidth]{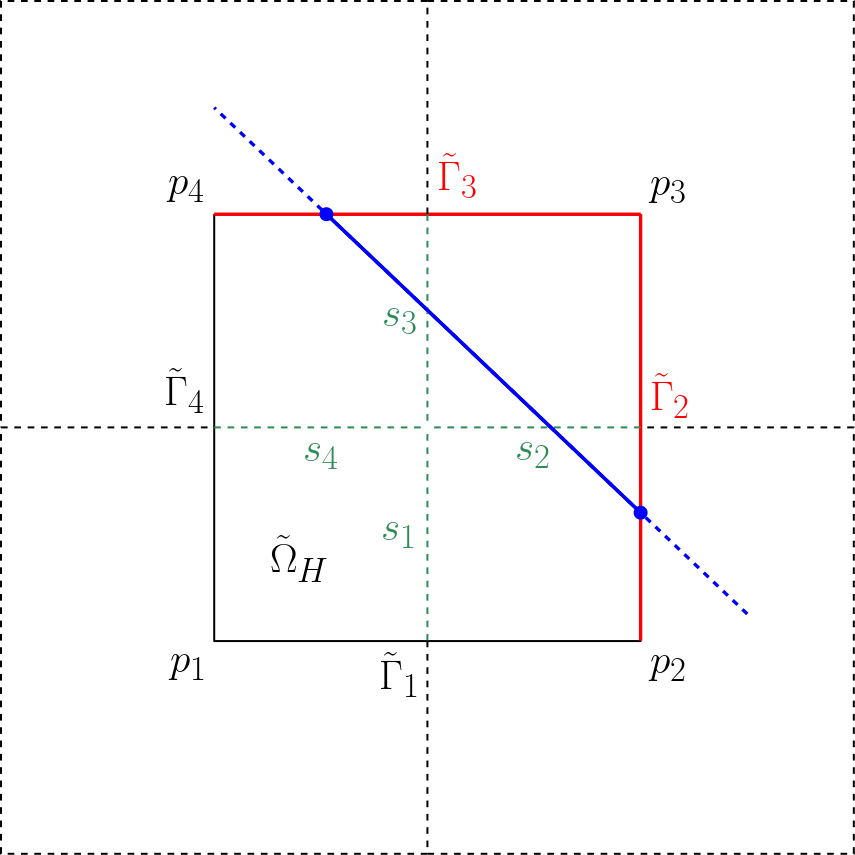}
	\caption{Example of an interaction region $\tilde{\Omega}_{H}$ (black and red solid lines) cut by a fracture (blue line). $p_1$ to $p_4$ are the pressures of the coarse cells (dashed lines) associated to the interaction region $\tilde{\Omega}_{H}$. The subinterfaces $s_1$ to $s_4$ are also highlighted (green dashed lines). In this case, the 1D problems on the boundary~\cite[equation~(2.6)]{Hesse} in $\tilde{\Gamma}_{2}$ and $\tilde{\Gamma}_{3}$ should be replaced with a 1D Lower-Dimensional Discrete Fracture-Matrix Model, explained in Appendix~\ref{sec:1D_DFMML}.}
	\label{fig:MSFV_red}
\end{figure}

Hence, the workflow of the multiscale modification of the LEDFM method is the following:
\begin{itemize}
	\item construct the coarse dual grid;
	\item find the interaction regions $\tilde{\Omega}_{H}$ of the dual grid that are cut by the fracture;
	\item for each intersected interaction region solve the four local fine scale problems~\cite[equation~(2.4)]{Hesse} and build the corresponding transmissibility matrices $\vec{T} \in \mathbb{R}^{4x4}$~\cite[equation~ (1.8)]{Hesse}, that relate the fluxes $f_1$ to $f_4$ across the four subinterfaces $s_1$ to $s_4$ of the interaction region $\tilde{\Omega}_{H}$ to the coarse cell pressures $p_1$ to $p_4$ in the corners of $\tilde{\Omega}_{H}$, see Fig.~\ref{fig:MSFV_red};
	\item find the subset of near-fracture coarse matrix cells, i.e. those for which at least one of the associated interaction regions is cut by the fracture;
	\item compute, using the transmissibility matrices $\vec{T}$, MSFV 9-point stencil modifications~\cite[equations~(1.6),~(1.9)--(1.16)]{Hesse} for the coarse matrix cells belonging to the subset of cells defined in the previous point, hence getting a near-fracture multi-point approximation.
\end{itemize}
In our case, differently from the original version of the MSFV method, the local fine scale problems~\cite[equation~(2.4)]{Hesse} for the interaction regions should take into account the presence of a fracture.

As already done for the M-F and M-M local problems, the conforming method mentioned in Subsection~\ref{ssec:discretization} has been adopted to solve the local problems with triangle matrix grids. Once again, the TPFA version of the method has been employed at the fine local scale for the same reasons explained in Subsection~\ref{ssec:discretization}.
Clearly, the fine grids used to solve the local problems should honour the fracture geometry as well as the four subinterfaces $s_1$ to $s_4$ of the interaction region $\tilde{\Omega}_{H}$, see Fig.~\ref{fig:MSFV_red}.
Moreover, before solving the local problems in the interaction regions, the latter are transformed in unit squares, as done for the M-F problems. 
The fracture cells within the interaction regions are also rescaled in the way previously described in Subsection~\ref{ssec:workflow}.

Since the effect of fractures on the solution is non-negligible, \emph{reduced} boundary conditions~\cite[equation~(2.6)]{Hesse} are the best choice in this case.
Here, the word \emph{reduced} is used to point out that the lower dimensional version of the problem considered in the interaction region should be solved on its boundary to obtain the conditions to be imposed on it.
Hence, the 1D problems that should be solved in order to get the reduced boundary conditions should take into account the presence of the fracture when it intersects any of the boundary segments $\tilde{\Gamma}_{n}$, with $n=1,2,3,4$, see Fig.~\ref{fig:MSFV_red} for an example.
This is done by replacing the 1D elliptic problems on the boundary~\cite[equation~(2.6)]{Hesse} with a 1D Lower-Dimensional Discrete Fracture-Matrix Model, described in Appendix~\ref{sec:1D_DFMML}, where the intersection of the fracture with the boundary is assumed to be a point in a mixed-dimensional setting.
For simplicity, we assume that the matrix fine scale permeability along the boundary segments is constant, meaning that the interaction region is associated to four matrix coarse cells characterized by the same permeability tensor, but the 1D model that will be introduced in Appendix~\ref{sec:1D_DFMML} can be extended to the heterogeneous case too.
It is worth mentioning that, thanks to their simplicity, the equations relative to the 1D Lower-Dimensional Discrete Fracture-Matrix Model can be solved analytically, so that no further source of errors nor computational cost are introduced in the overall procedure.     
  
\section{Numerical Tests} \label{sec:numtests}

In this section some numerical tests in two-dimensional fractured porous media are presented to test the performance of the newly developed LEDFM.
First, convergence tests are carried out in the case of incompressible single-phase Darcy flow. 
Then, a tracer transport test is presented, where the advective velocity field is obtained from the single-phase flow problem.

\subsection{Convergence Tests} \label{ssec:conv_tests}
 
The local method is compared to two other numerical schemes, namely EDFM and pEDFM, in terms of convergence results.
The MATLAB Reservoir Simulation Toolbox (MRST), see \cite{MRST}, is used to perform the simulations: EDFM and pEDFM are already implemented in the toolbox, while routines related to local problems solutions have been written and integrated in MRST to obtain the implementation of LEDFM.
In particular, the MESH2D toolbox \cite{Engwirda} has been used to generate conforming triangular grids for the local problems.
It should be pointed out that some modifications have been made to the pEDFM code of MRST to correctly simulate the fractured systems presented in this section. The most important changes, along with some general warnings about pEDFM implementation, are detailed in Appendix~\ref{sec:warn_pEDFM}.
In the results reported in the next subsections we will refer to the original implementation of the pEDFM code of MRST as ``MRST pEDFM'' and to the updated version of the code as ``updated pEDFM''.

To compare results in a quantitative way, we compute a fine scale equidimensional reference solution $p_{ref}$ and the error indicator is the $L^2$ pressure error in the matrix domain, computed with the following formula, borrowed from \cite{Flemisch2}:
\begin{equation*}
	e_{m}^{2}=\frac{1}{|\Omega|\left(\Delta p_{ref}\right)^{2}} \sum_{I=K_{ref, m} \cap K_{m}}|I|\left(\left.p\right|_{K_{m}}-\left.p_{ref}\right|_{K_{ref, m}}\right)^{2},
\end{equation*}
where $|\Omega|$ indicates the size of the matrix domain, $\Delta p_{ref}$ the maximum pressure difference in the reference solution for normalization, $I$ a generic intersection between a reference matrix element $K_{ref, m}$ and a coarse one $K_m$, $|I|$ the size of the intersection and $p$ the coarse pressure field. 
We stress that for the computation of the pressure matrix error defined above, only cells belonging to the matrix part of the equidimensional grid are considered.

In Subsection~\ref{ssec:perm_oblique_frac} we show the convergence results obtained in the case of an oblique highly permeable fracture for the different embedded methods, then in Subsection~\ref{ssec:imperm_horiz_frac} we do the same in the case of a horizontal impermeable fracture and lastly in Subsection~\ref{ssec:imperm_oblique_frac} we run the same analyses for the oblique fracture presented initially, now considered to be impermeable with respect to the surrounding porous matrix.

\subsubsection{Test 1 -- Permeable Oblique Fracture} \label{ssec:perm_oblique_frac}

\begin{figure}[!t]
	\centering
	\includegraphics[width=.45\textwidth]{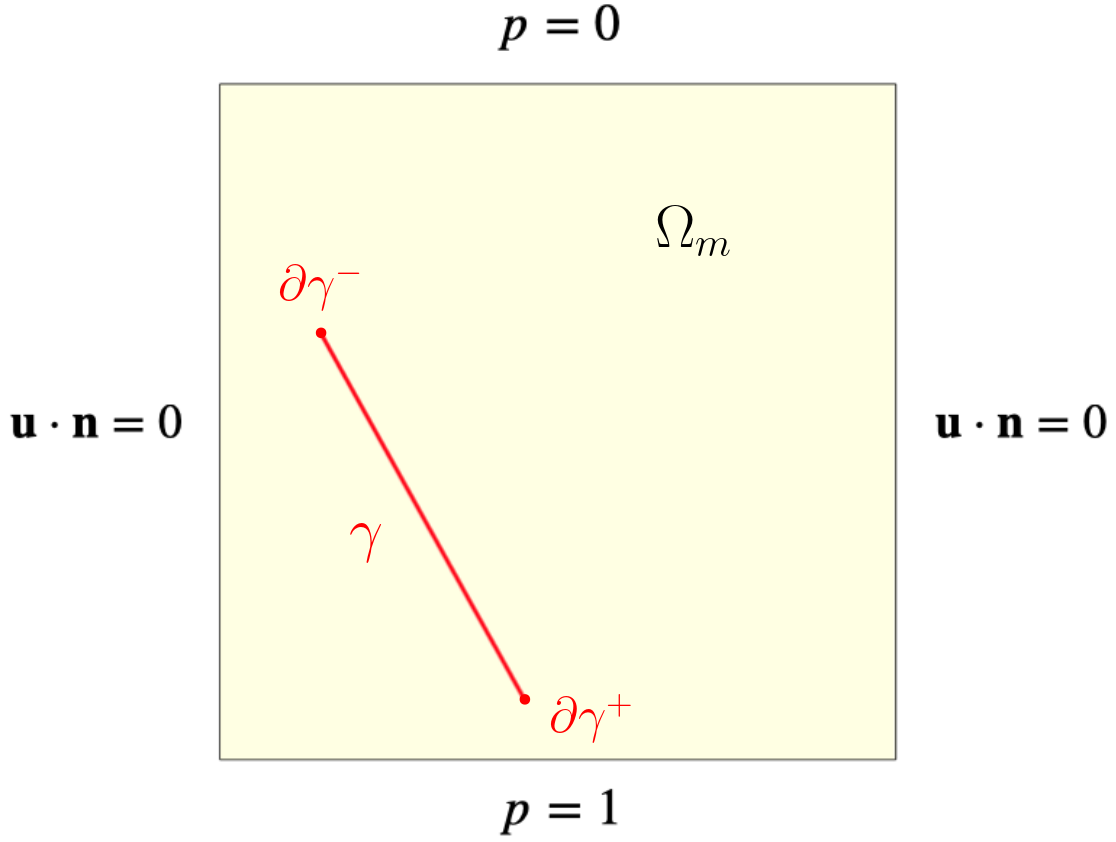}
	\caption{Oblique fracture problem domain. Boundary conditions are put in evidence.}
	\label{fig:oblique_problem}
\end{figure}

The first example considers a two-dimensional square domain $\Omega = [0,1]^2$ where an immersed permeable oblique fracture $\Omega_f$ is present, whose aperture is constant and equal to $d = 10^{-4}$. The fracture domain is then defined as
\begin{equation*}
	\Omega_f = \left\{ \vec{x} \in \Omega : \vec{x} = \vec{s} + r\vec{n}, \ \vec{s} \in \gamma, \ r \in \left( -d/2, d/2 \right) \right\},
\end{equation*}
where $\gamma$ is the line segment connecting the points $\vec{s}_1 = [0.15, 0.63]^\top$ and $\vec{s}_2 = [0.45, 0.09]^\top$. Using the definitions introduced in Section~\ref{sec:goveqns}, we identify $\vec{s}_1$ with $\partial \gamma^{-}$ and $\vec{s}_2$ with $\partial \gamma^{+}$, so that the unit tangential vector $\boldsymbol{\tau}$ and unit normal vector $\vec{n}$ to $\gamma$ can be defined as explained in Section~\ref{sec:goveqns}.
The porous matrix domain is then given by $\Omega_{m} = \Omega \setminus \Omega_f$, or, in the case of hybrid dimensional methods, as $\Omega_{m} = \Omega \setminus \gamma$.
Fig.~\ref{fig:oblique_problem} shows the considered domain along with boundary conditions, where the fracture is identified with $\gamma$.
The flow is driven in the upward direction by imposing pressure values of $p = 1$ and $p = 0$ on the bottom and top boundaries of the domain, respectively, while no-flow conditions are set on the left and right sides. 
No source terms are present in the domain.
Viscosity is set to $\mu = 1$ in the whole domain.
The permeability of the matrix region is considered to be homogeneous, isotropic and equal to $\vec{K}_m = k_m \vec{I}$, with $k_m = 1$.
The same stands for the fracture permeability tensor, which is equal to $\vec{K}_f = k_f \vec{I}$.

\begin{figure}[p]
	\centering
	\begin{subfigure}[b]{0.475\textwidth}
		\centering
		\includegraphics[width=.9\textwidth]{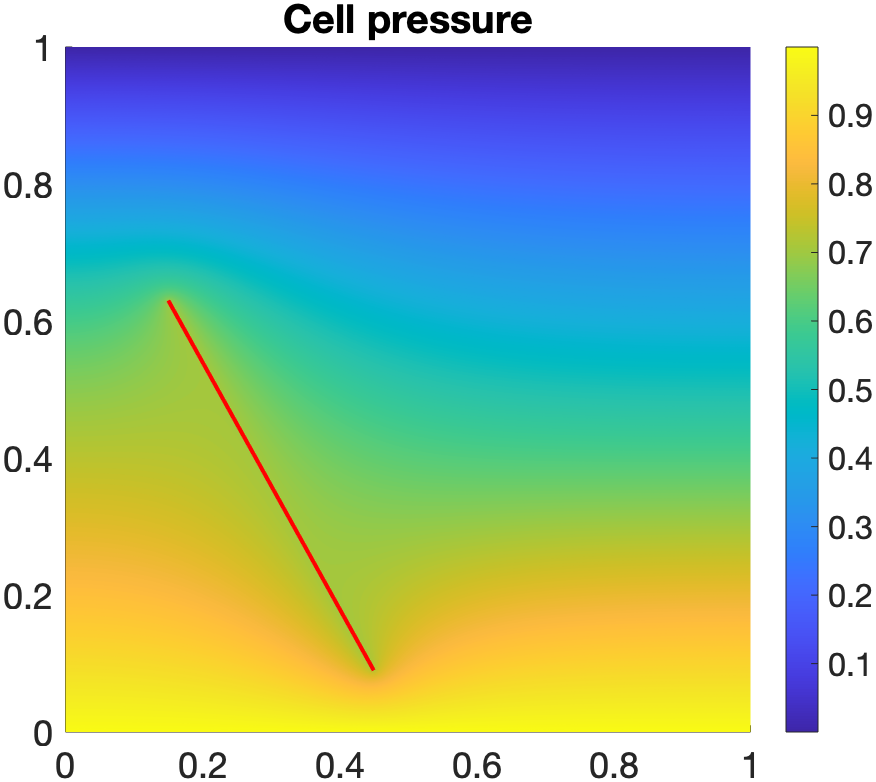}
		\caption{}
		\label{fig:oblique_perm_1e+8_press}
	\end{subfigure}
	\hfill
	\begin{subfigure}[b]{0.515\textwidth}
		\centering
		\includegraphics[width=\textwidth]{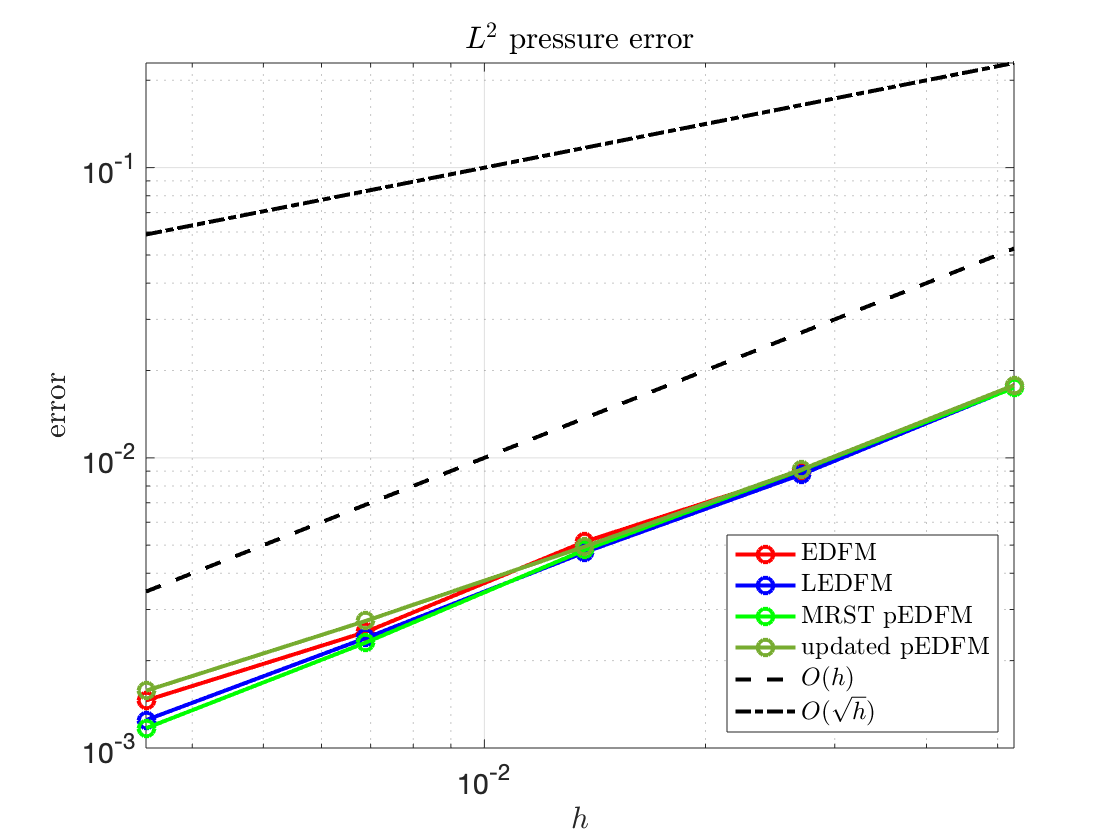}
		\caption{}
		\label{fig:oblique_perm_1e+8_conv_corr}
	\end{subfigure}
	\caption{Permeable oblique fracture. Scenario $k_f/k_m = 10^{8}$. (a) Reference solution. (b) Convergence plot.}
	\label{fig:oblique_perm_1e+8}
\end{figure}

\begin{figure}[p]
	\centering
	\begin{subfigure}[b]{0.475\textwidth}
		\centering
		\includegraphics[width=.9\textwidth]{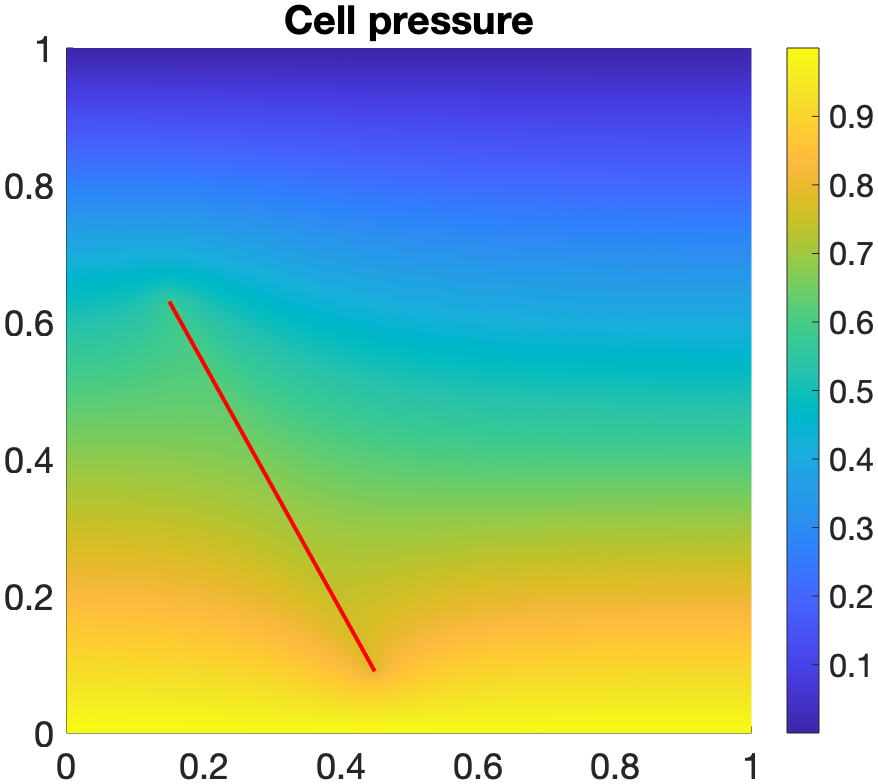}
		\caption{}
		\label{fig:oblique_perm_1e+4_press}
	\end{subfigure}
	\hfill
	\begin{subfigure}[b]{0.515\textwidth}
		\centering
		\includegraphics[width=\textwidth]{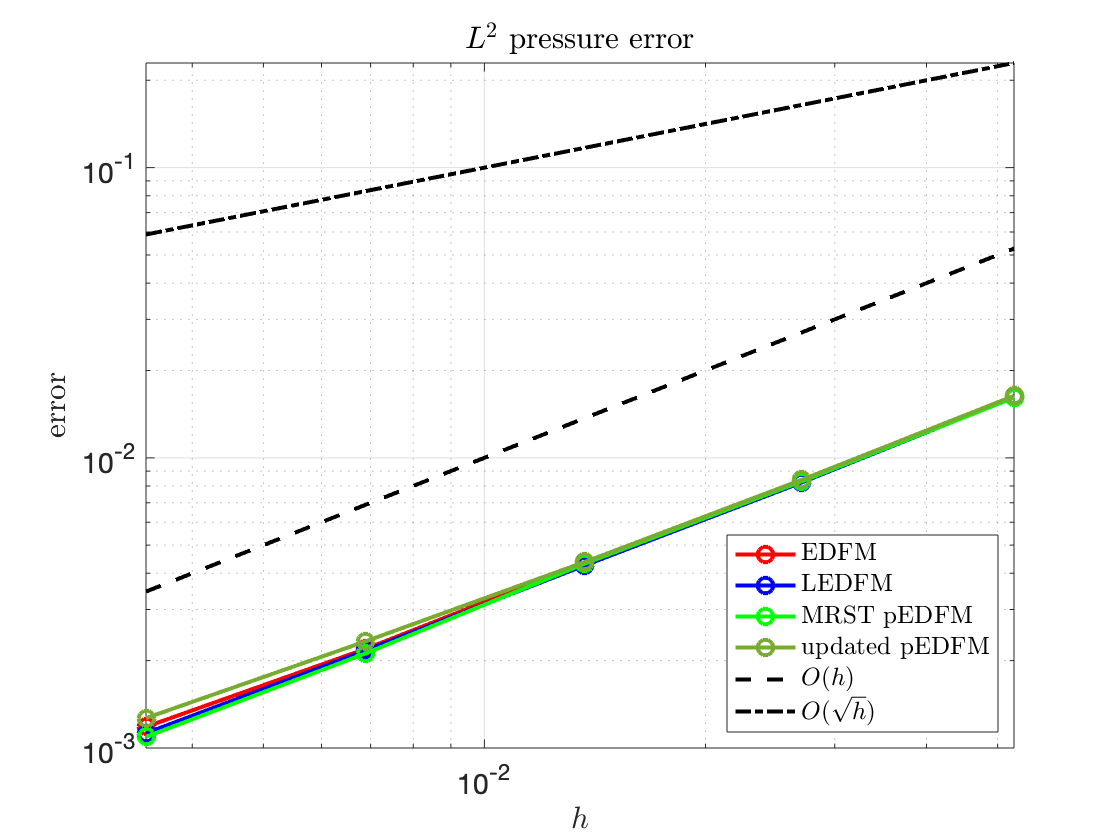}
		\caption{}
		\label{fig:oblique_perm_1e+4_conv_corr}
	\end{subfigure}
	\caption{Permeable oblique fracture. Scenario $k_f/k_m = 10^{4}$. (a) Reference solution. (b) Convergence plot.}
	\label{fig:oblique_perm_1e+4}
\end{figure}

Two different scenarios are considered: one in which the fracture is 8 orders of magnitude more conductive than the matrix ($k_f/k_m = 10^{8}$), and another with a lower conductivity contrast, where the fracture is 4 orders of magnitude more conductive than the matrix ($k_f/k_m = 10^{4}$).

The equidimensional reference solution is computed on a very fine grid that discretizes the matrix in $1,512,320$ triangular elements and the fracture in $5,376$ quadrilateral elements.
The software Gmsh has been used to create the grid \cite{Geuzaine}.
The Virtual Element Method of order $k = 2$, already implemented in MRST and based on \cite{Klemetsdal}, is used to discretize the problem.

Fig.~\ref{fig:oblique_perm_1e+8_press} shows the reference solution in the scenario in which $k_f/k_m = 10^{8}$, while in Fig.~\ref{fig:oblique_perm_1e+8_conv_corr} the corresponding convergence plot for all the competing embedded methods is reported. We employ $N \times N$ square Cartesian grids, where $N$ is gradually increased. To be specific, grids with $N = 19, 37, 73, 145, 289$ cells over each axis are considered.
We observe that all methods show approximately the same behaviour, i.e. linear convergence. 

Similar results are obtained for the second scenario in which $k_f/k_m = 10^{4}$. In particular, Fig.~\ref{fig:oblique_perm_1e+4_press} shows the corresponding reference solution, while Fig.~\ref{fig:oblique_perm_1e+4_conv_corr} the convergence plots, that show once again linear behaviours for all the methods.

\subsubsection{Test 2 -- Impermeable Horizontal Fracture} \label{ssec:imperm_horiz_frac}
\begin{figure}[!t]
	\centering
	\includegraphics[width=.45\textwidth]{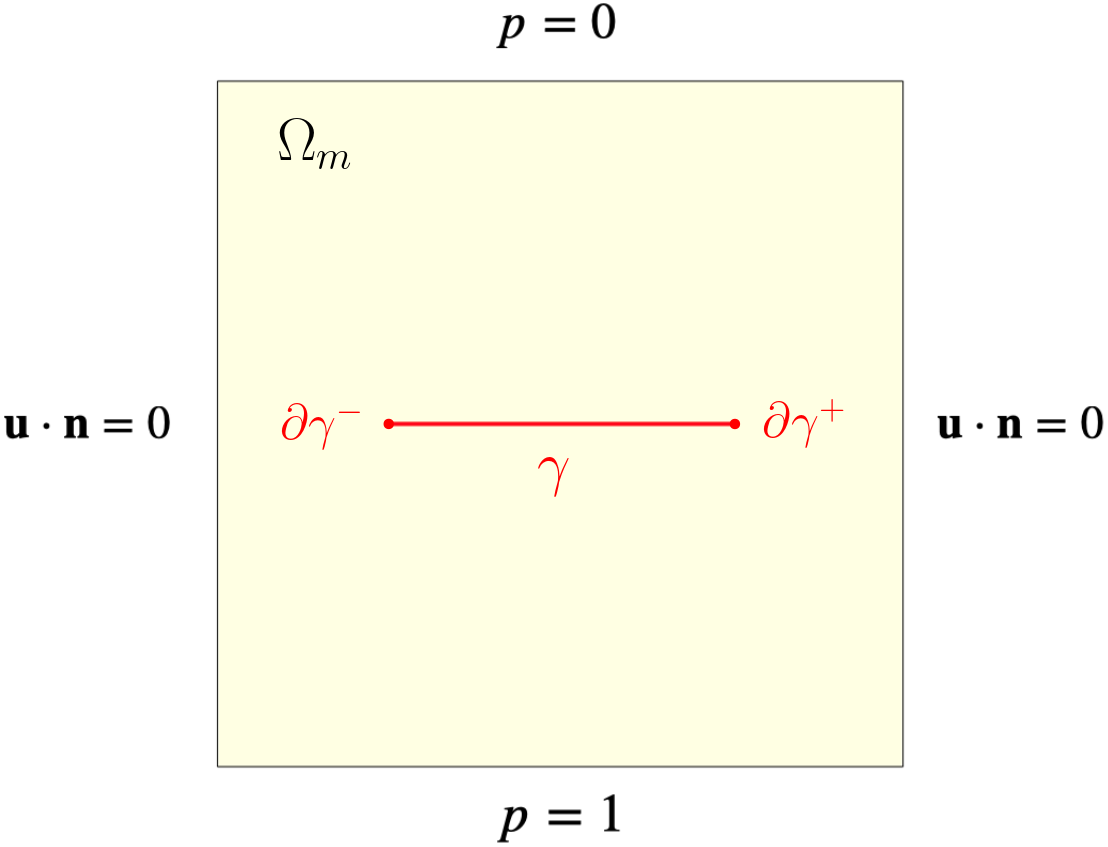}
	\caption{Horizontal fracture problem domain. Boundary conditions are put in evidence.}
	\label{fig:horizontal_problem}
\end{figure}

The second example considers the same two-dimensional square domain $\Omega = [0,1]^2$, but this time an immersed impermeable horizontal fracture is present. In particular $\gamma$ now is the line segment joining the points $\vec{s}_1 = [0.25, 0.5]^\top$ and $\vec{s}_2 = [0.75, 0.5]^\top$. The same aperture of the previous case is considered. 
Fig.~\ref{fig:horizontal_problem} depicts the considered domain along with boundary conditions, which are the same of the previous example.
As in Test 1, no source terms are present. 

\begin{figure}[p]
	\centering
	\begin{subfigure}[b]{0.475\textwidth}
		\centering
		\includegraphics[width=.9\textwidth]{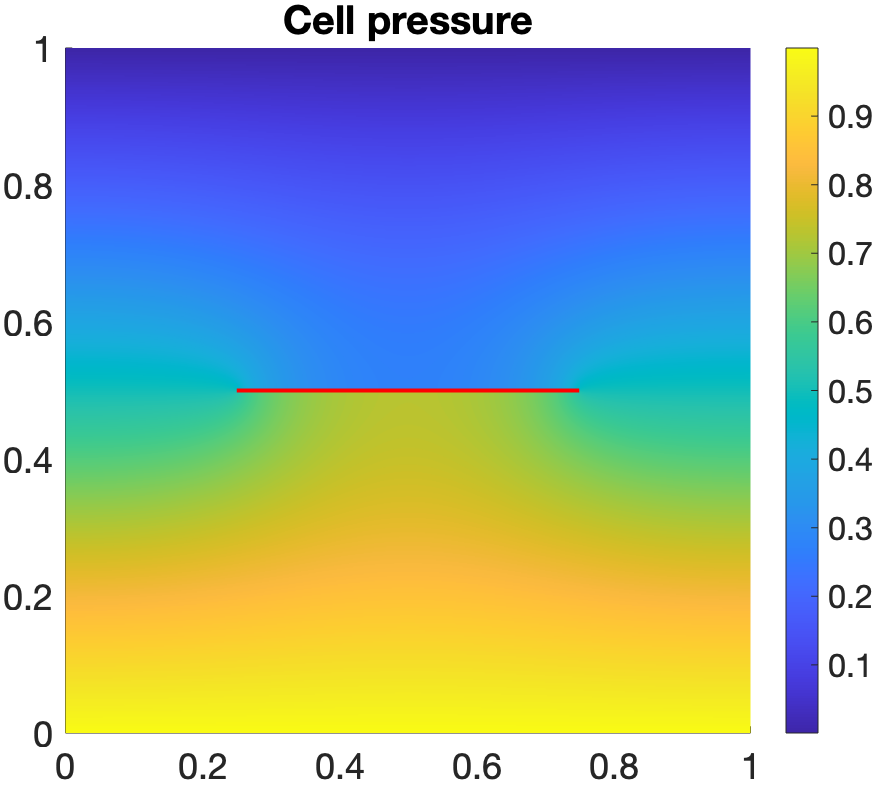}
		\caption{}
		\label{fig:horiz_imperm_1e-8_press}
	\end{subfigure}
	\hfill
	\begin{subfigure}[b]{0.515\textwidth}
		\centering
		\includegraphics[width=\textwidth]{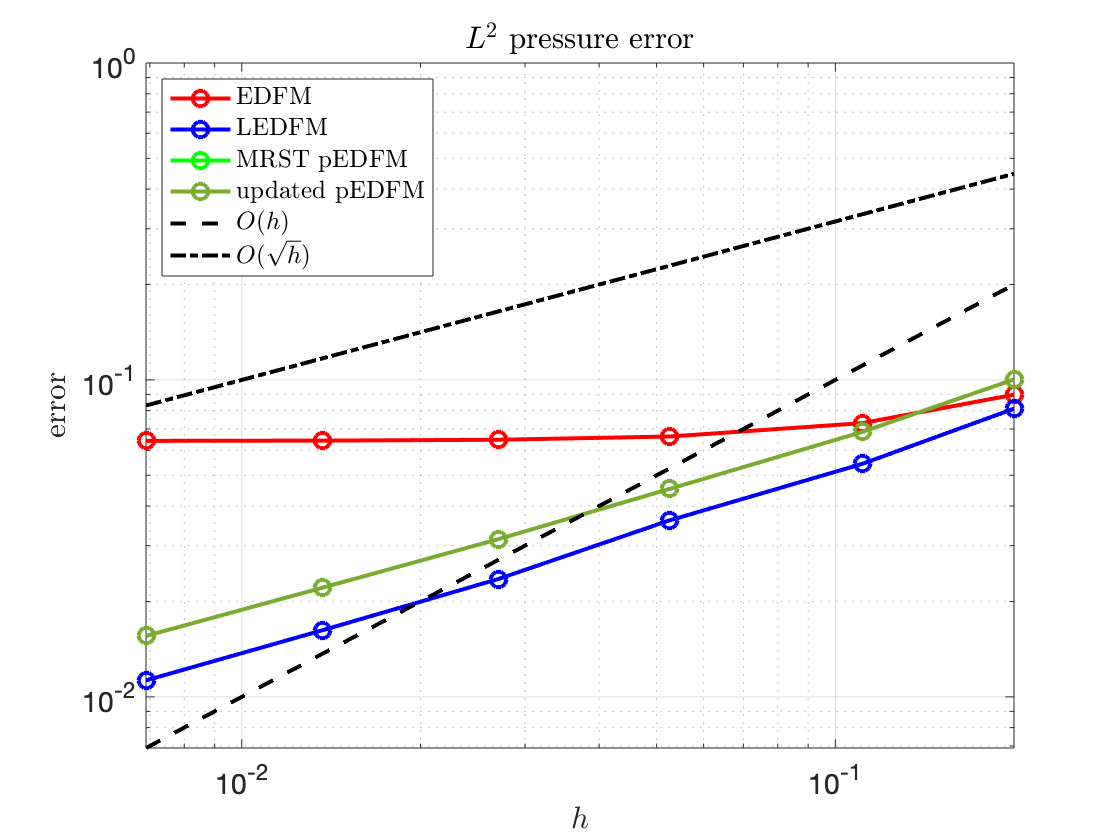}
		\caption{}
		\label{fig:horiz_imperm_1e-8_conv_i}
	\end{subfigure}
	\caption{Impermeable horizontal fracture. Scenario $k_f/k_m = 10^{-8}$. (a) Reference solution. (b) Convergence plot.}
	\label{fig:horiz_imperm_1e-8}
\end{figure} 

\begin{figure}[p]
	\centering
	\begin{subfigure}[b]{0.495\textwidth}
		\centering
		\includegraphics[width=.9\textwidth]{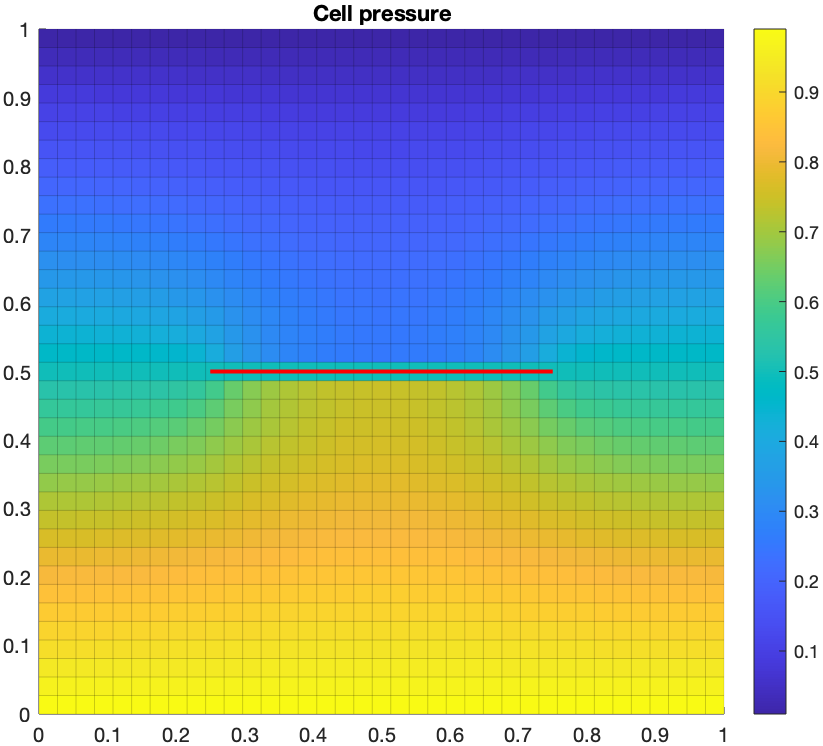}
		\caption{}
		\label{fig:horiz_frac_LEDFM_37x37_Kratio1e-8_p}
	\end{subfigure}
	\hfill
	\begin{subfigure}[b]{0.495\textwidth}
		\centering
		\includegraphics[width=.9\textwidth]{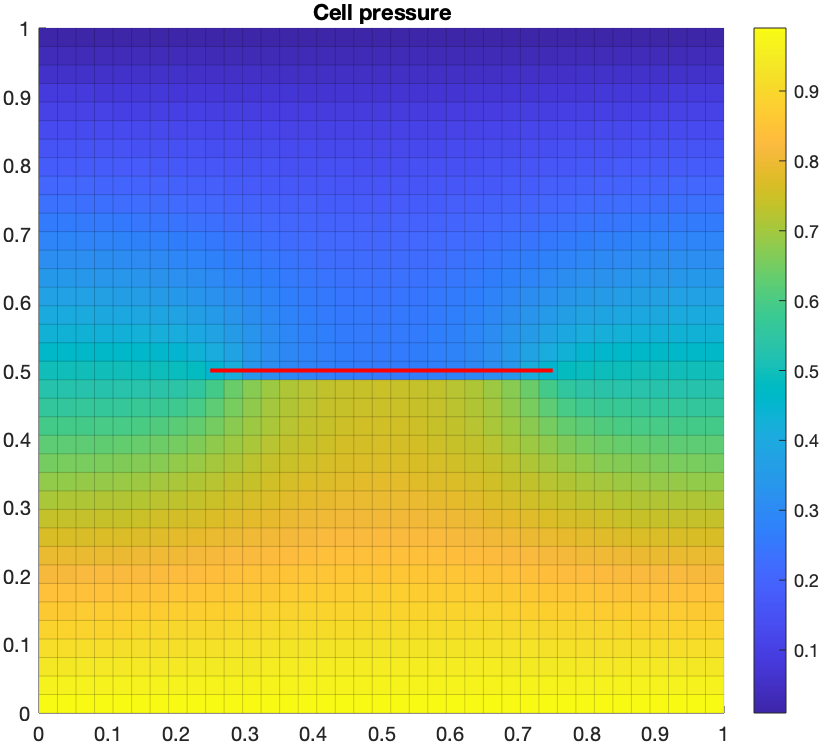}
		\caption{}
		\label{fig:horiz_frac_pEDFM_i_corr_37x37_Kratio1e-8_p}
	\end{subfigure}
	\caption{Impermeable horizontal fracture. Scenario $k_f/k_m = 10^{-8}$. (a) LEDFM solution on $37\times37$ grid. (b) pEDFM solution on $37\times37$ grid.}
	\label{fig:horiz_frac_37x37_Kratio1e-8_p}
\end{figure}

Viscosity and permeability tensors are also taken to be equal to Test 1, but, since the fracture considered now is impermeable, the two scenarios that will be examined are $k_f/k_m = 10^{-8}$ and $k_f/k_m = 10^{-4}$.

The equidimensional reference solution is computed on a fine rectilinear grid that discretizes the matrix in $1,000,000$ rectangular elements and the fracture in $10,000$ rectangular elements as well. 
The grid is progressively refined in the $y$ direction as we approach the fracture.
MRST has been used for generating the grid and solving the problem using a simple TPFA scheme, since the grid is $\vec{K}$-orthogonal.

Fig.~\ref{fig:horiz_imperm_1e-8_press} shows the reference solution in the scenario in which $k_f/k_m = 10^{-8}$, while in Fig.~\ref{fig:horiz_imperm_1e-8_conv_i} the corresponding convergence plots for all the competing embedded methods are reported. 
The square Cartesian grids have, for the different refinement levels, $N = 5, 9, 19, 37, 73, 145$ cells over each axis.
With this choice we avoid that the horizontal fracture coincides with the matrix grid interfaces, which is a limit case not yet handled by the local method. 
We notice that the classical EDFM is not able to provide the correct solution since the fracture is impermeable, while pEDFM and LEDFM can do it.
In particular, we note that the MRST and updated pEDFM curves coincide since the fracture is highly impermeable and the fracture projection path is determined correctly for both cases.
However, the order of convergence is approximately equal to $O(\sqrt{h})$ and not linear as in the permeable case for both pEDFM and LEDFM. This is due to the fact that, see Fig.~\ref{fig:horiz_frac_37x37_Kratio1e-8_p}, the fracture is cutting matrix grid cells. There, in embedded FV methods, pressure is represented by a single value so an error of $O(1)$ is observed with respect to the true, discontinuous solution. If the number of the cut cells is $O(N)$ this leads to the observed $\sqrt{h}$ convergence. Indeed, if we exclude the cut cells from error computation, linear convergence is restored.
We also notice that lower errors are obtained with the local method compared to pEDFM.

\begin{figure}[p]
	\centering
	\begin{subfigure}[b]{0.475\textwidth}
		\centering
		\includegraphics[width=.9\textwidth]{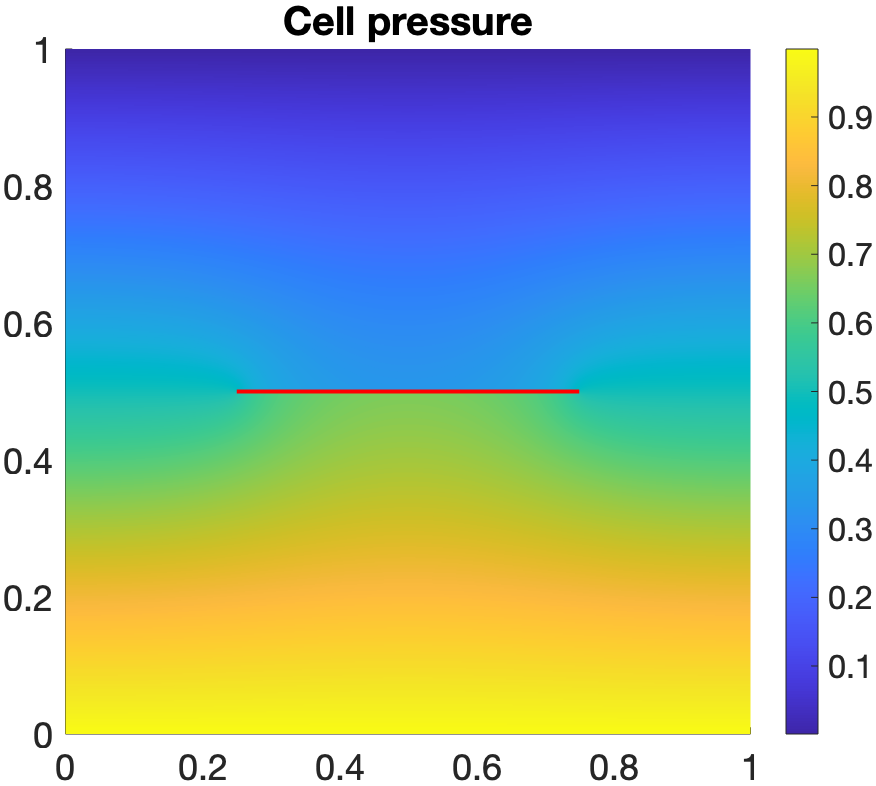}
		\caption{}
		\label{fig:horiz_imperm_1e-4_press}
	\end{subfigure}
	\hfill
	\begin{subfigure}[b]{0.515\textwidth}
		\centering
		\includegraphics[width=\textwidth]{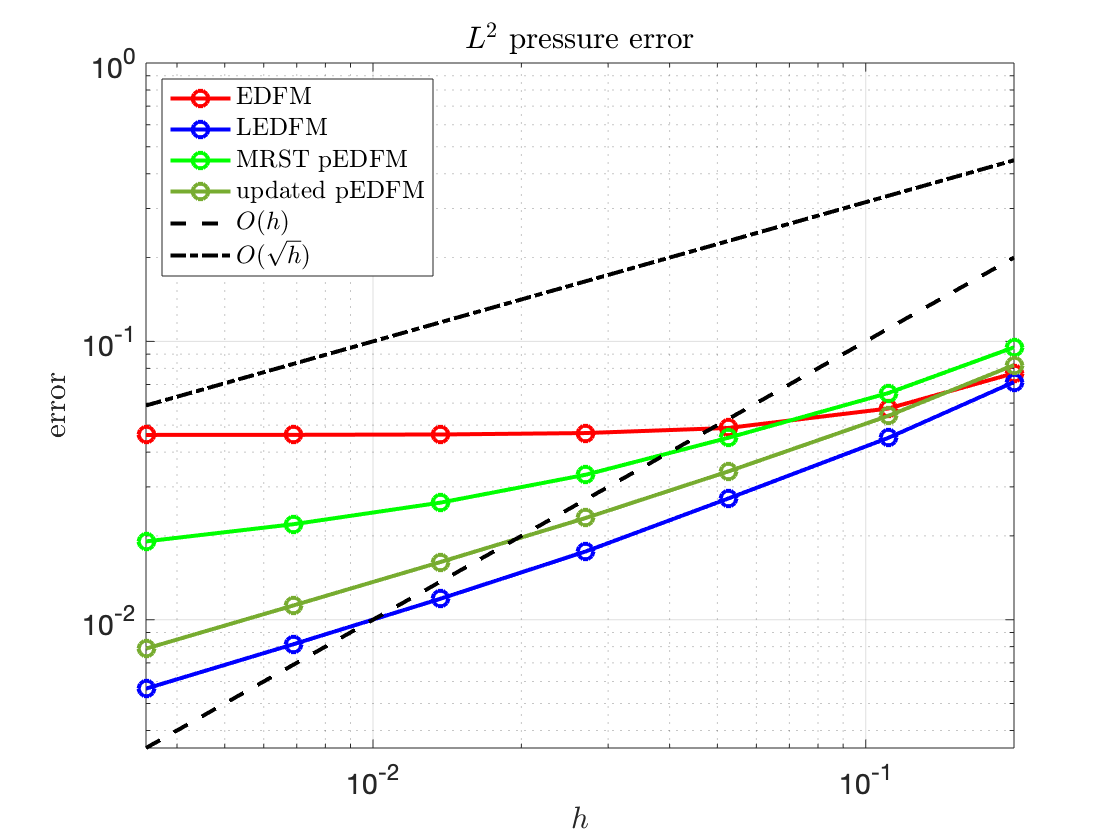}
		\caption{}
		\label{fig:horiz_imperm_1e-4_conv_i}
	\end{subfigure}
	\caption{Impermeable horizontal fracture. Scenario $k_f/k_m = 10^{-4}$. (a) Reference solution. (b) Convergence plot.}
	\label{fig:horiz_imperm_1e-4}
\end{figure} 

\begin{figure}[p]
	\centering
	\begin{subfigure}[b]{0.327\textwidth}
		\centering
		\includegraphics[width=\textwidth]{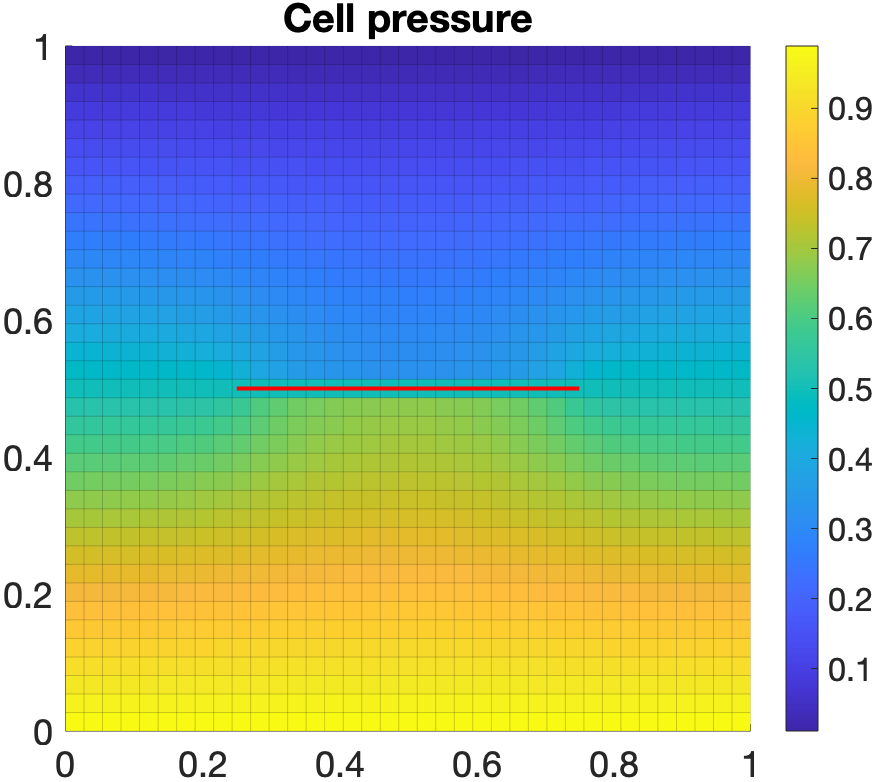}
		\caption{}
		\label{fig:horiz_imperm_1e-4_LEDFM_37x37}
	\end{subfigure}
	\begin{subfigure}[b]{0.327\textwidth}
		\centering
		\includegraphics[width=\textwidth]{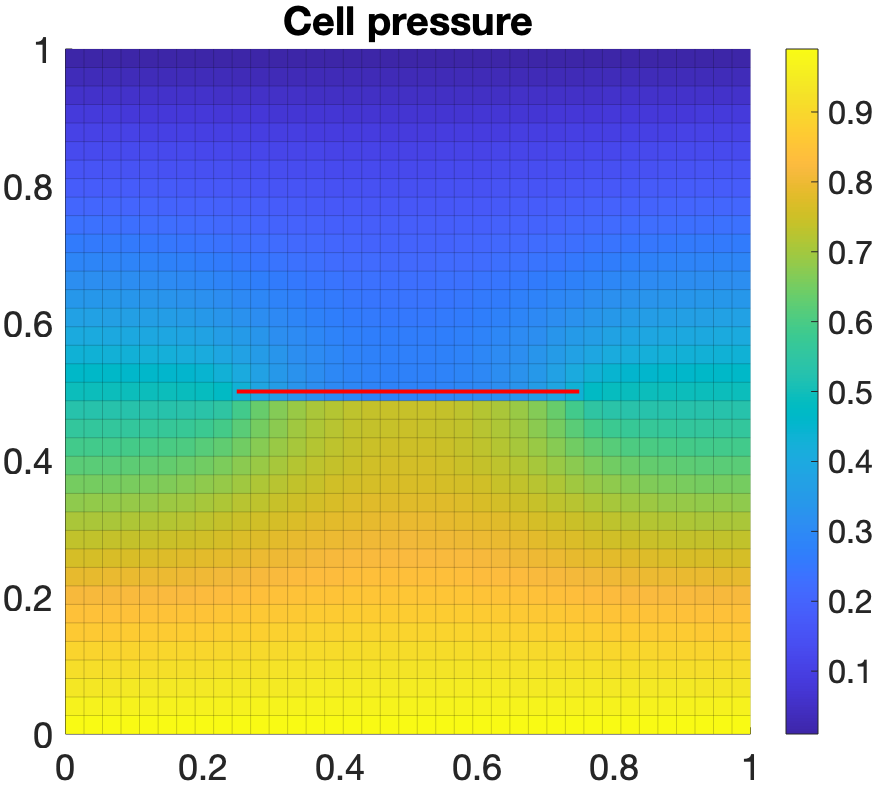}
		\caption{}
		\label{fig:horiz_imperm_1e-4_pEDFM_37x37}
	\end{subfigure}
	\begin{subfigure}[b]{0.327\textwidth}
		\centering
		\includegraphics[width=\textwidth]{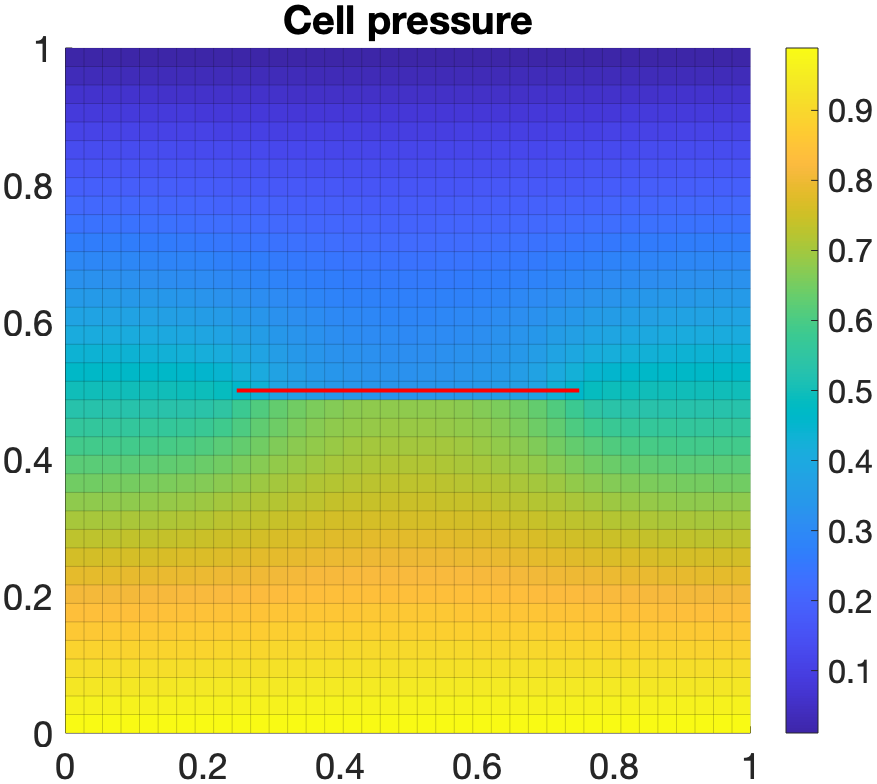}
		\caption{}
		\label{fig:horiz_imperm_1e-4_pEDFM_i_c_37x37}
	\end{subfigure}
	\caption{Impermeable horizontal fracture. Scenario $k_f/k_m = 10^{-4}$. (a) LEDFM solution on $37\times37$ grid. (b) MRST pEDFM solution on $37\times37$ grid. (c) updated pEDFM solution on $37\times37$ grid.}
	\label{fig:horiz_imperm_1e-4_press_37x37}
\end{figure}  

Moreover, from Fig.~\ref{fig:horiz_frac_37x37_Kratio1e-8_p} we observe that pressure maps for the pEDFM and LEDFM method are quite different, especially in correspondence of the matrix cells cut by the fracture: in pEDFM the fracture is projected to the lower matrix grid interfaces and the pressure is discontinuous across them, while in LEDFM the pressure jump is redistributed on the upper and lower interfaces.

Moving to the second scenario in which $k_f/k_m = 10^{-4}$, in Fig.~\ref{fig:horiz_imperm_1e-4_press} and Fig.~\ref{fig:horiz_imperm_1e-4_conv_i} the reference solution and the convergence plots are shown, respectively.
Here, an additional refinement level for the square Cartesian grids corresponding to $N = 289$ has been considered.
The main difference in performance with respect to the previous case is observed for the MRST pEDFM curve.
Indeed, the convergence rate of the MRST pEDFM degrades as the grid is refined, while this does not occur for the updated version and LEDFM.
This occurs since the transmissibility formulae implemented in MRST apply a simple harmonic averaging of the permeabilities, as explained in detail in Appendix~\ref{sec:1DpEDFM}.

Fig.~\ref{fig:horiz_imperm_1e-4_press_37x37} depicts the pressure maps obtained for the second scenario in the case of a $37\times37$ grid for the LEDFM, MRST pEDFM and updated pEDFM method, respectively.
Observations similar to those made for Fig.~\ref{fig:horiz_frac_37x37_Kratio1e-8_p} apply to the LEDFM and updated pEDFM pressure maps. 
The pressure map of the MRST pEDFM, instead, shows an higher pressure jump at the fracture with respect to the other methods. This is, once again, the consequence of using the wrong pEDFM transmissibility formulae with intermediate (impermeable) fracture-matrix permeability contrast values, as can be seen in Fig.~\ref{fig:1D_pEDFM_1e-4} for the 1D case.

\subsubsection{Test 3 -- Impermeable Oblique Fracture} \label{ssec:imperm_oblique_frac}
The third and last example considers the same configuration described in the first test, see Fig.~\ref{fig:oblique_problem}, but this time the oblique fracture is considered to be impermeable. 
In particular, two scenarios in which $k_f/k_m = 10^{-8}$ and $k_f/k_m = 10^{-4}$ will be examined.

\begin{figure}
	\centering
	\begin{subfigure}[b]{0.475\textwidth}
		\centering
		\includegraphics[width=.9\textwidth]{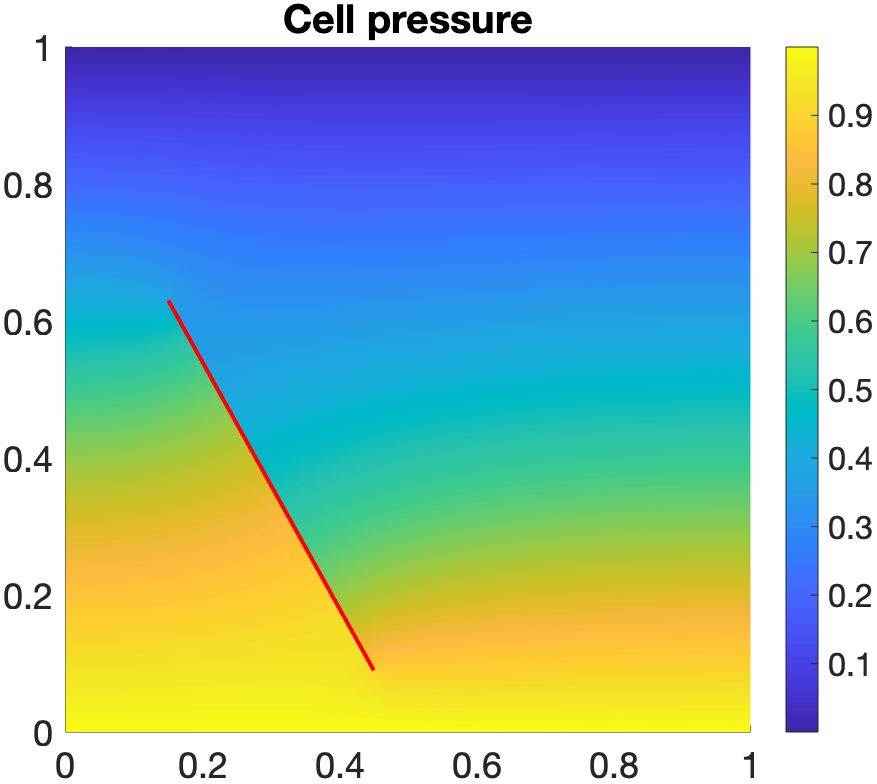}
		\caption{}
		\label{fig:oblique_imperm_1e-8_press}
	\end{subfigure}
	\hfill
	\begin{subfigure}[b]{0.515\textwidth}
		\centering
		\includegraphics[width=\textwidth]{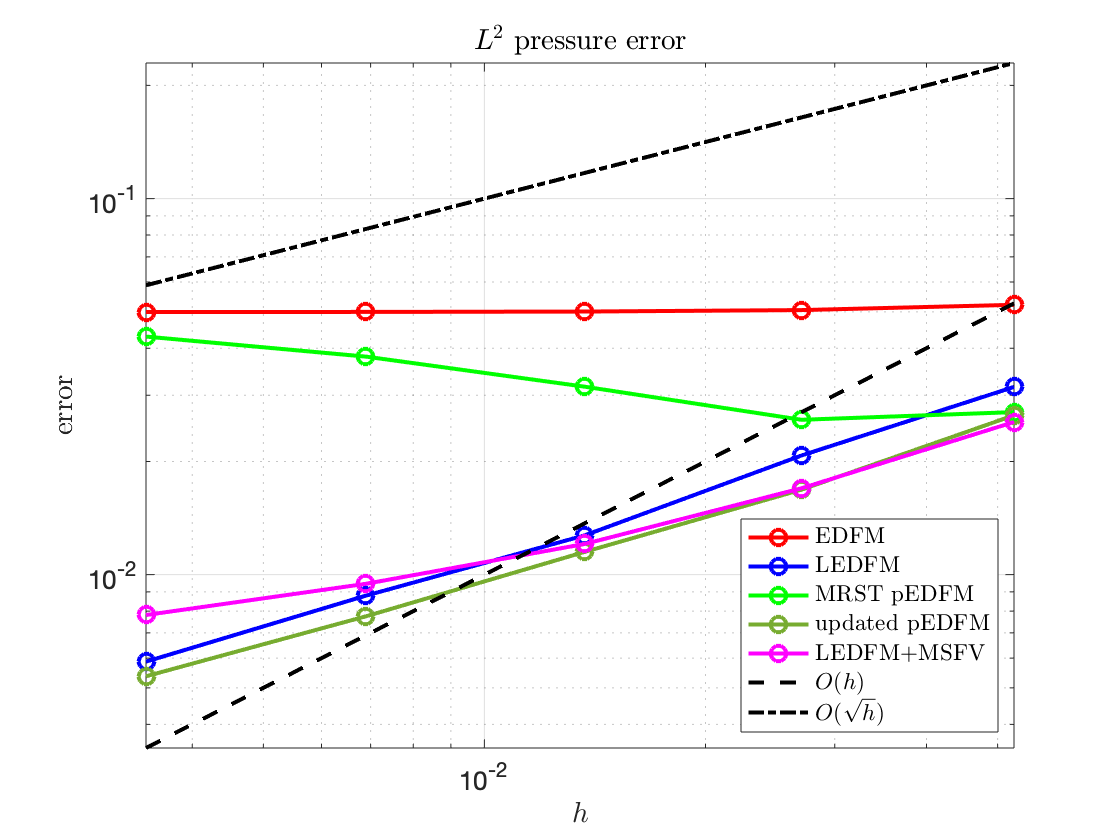}
		\caption{}
		\label{fig:oblique_imperm_1e-8_conv_i}
	\end{subfigure}
	\caption{Impermeable oblique fracture. Scenario $k_f/k_m = 10^{-8}$. (a) Reference solution. (b) Convergence plot.}
	\label{fig:oblique_imperm_1e-8}
\end{figure}

Fig.~\ref{fig:oblique_imperm_1e-8_press} shows the reference solution in the scenario in which $k_f/k_m = 10^{-8}$, while in Fig.~\ref{fig:oblique_imperm_1e-8_conv_i} the corresponding convergence plot for all the competing embedded methods is reported, including the multiscale modification of the LEDFM method (LEDFM+MSFV). 
The square Cartesian grids have, for the different refinement levels, $N = 19, 37, 73, 145, 289$ cells over each axis, as in the permeable case.
As already mentioned in the second example, the classical EDFM is not capable of representing this scenario due to the fact that the fracture is impermeable, regardless of the grid resolution.
The LEDFM method, instead, still shows a good convergence behaviour, corresponding to a convergence rate of approximately $O(\sqrt{h})$, as in the horizontal fracture case. 
Similar results are obtained for the multiscale modification of the local method. In particular, the latter performs slightly better than the standard version for coarser grids and slightly worse for finer grids.
The updated pEDFM method also performs well, with slightly lower errors compared to both LEDFM versions.
In the MRST pEDFM case, instead, as the grid is refined the error rises.
This is because, due to the MRST implementation, the projection of the fracture may not create a continuous path and, as the grid is refined, more and more holes are created in the projection, allowing the fluid to flow across the fracture through multiple spots. These ``leaks'' greatly lower the barrier effect of the fracture, so that very inaccurate pressure solutions are obtained.

\begin{figure}
	\centering
	\begin{subfigure}[b]{0.475\textwidth}
		\centering
		\includegraphics[width=.9\textwidth]{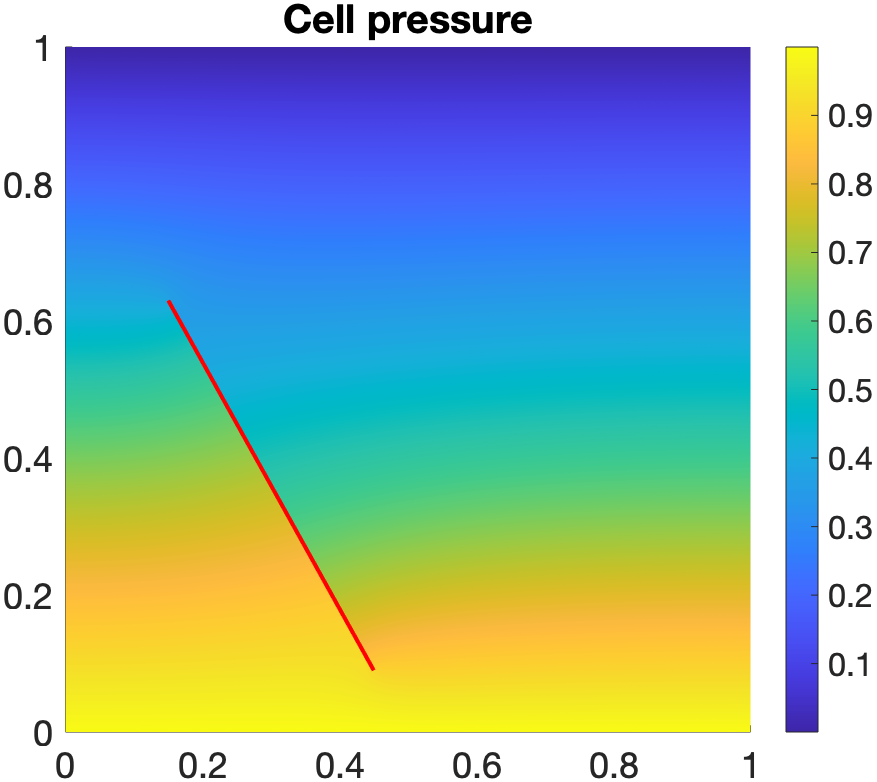}
		\caption{}
		\label{fig:oblique_imperm_1e-4_press}
	\end{subfigure}
	\hfill
	\begin{subfigure}[b]{0.515\textwidth}
		\centering
		\includegraphics[width=\textwidth]{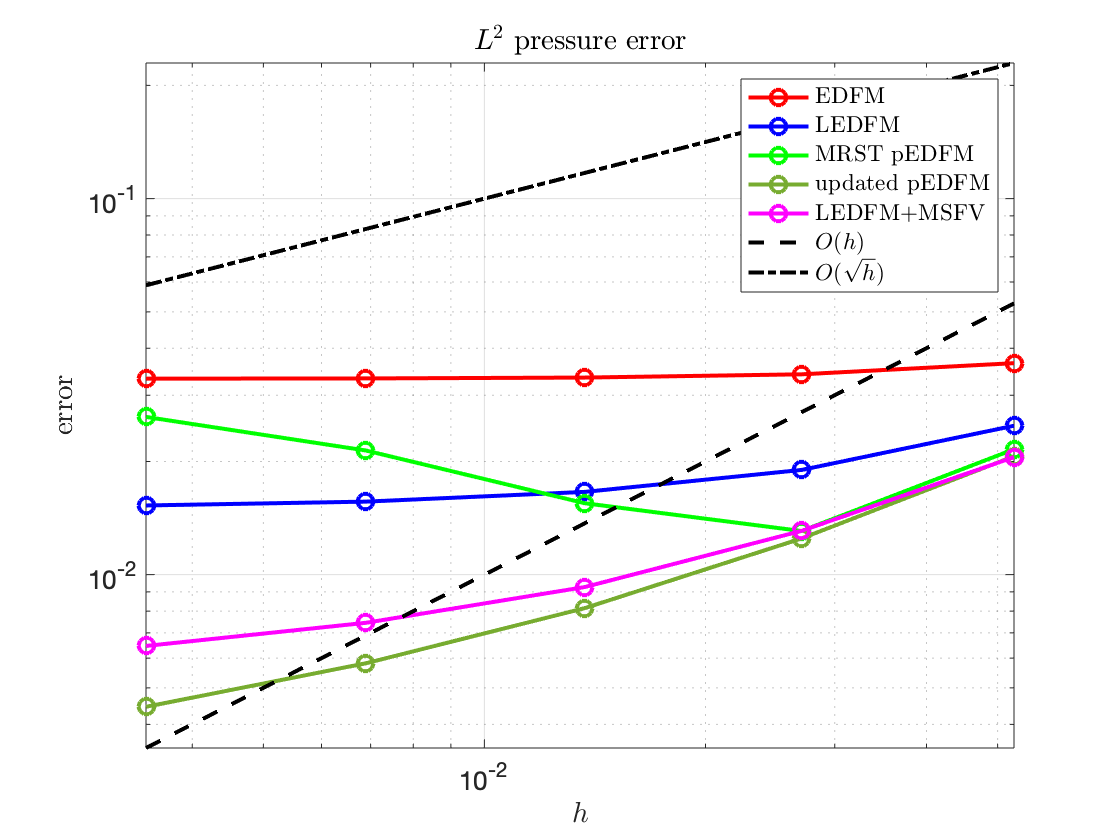}
		\caption{}
		\label{fig:oblique_imperm_1e-4_conv_i}
	\end{subfigure}
	\caption{Impermeable oblique fracture. Scenario $k_f/k_m = 10^{-4}$. (a) Reference solution. (b) Convergence plot.}
	\label{fig:oblique_imperm_1e-4}
\end{figure}

\begin{figure}
	\centering
	\begin{subfigure}[b]{0.495\textwidth}
		\centering
		\includegraphics[width=\textwidth]{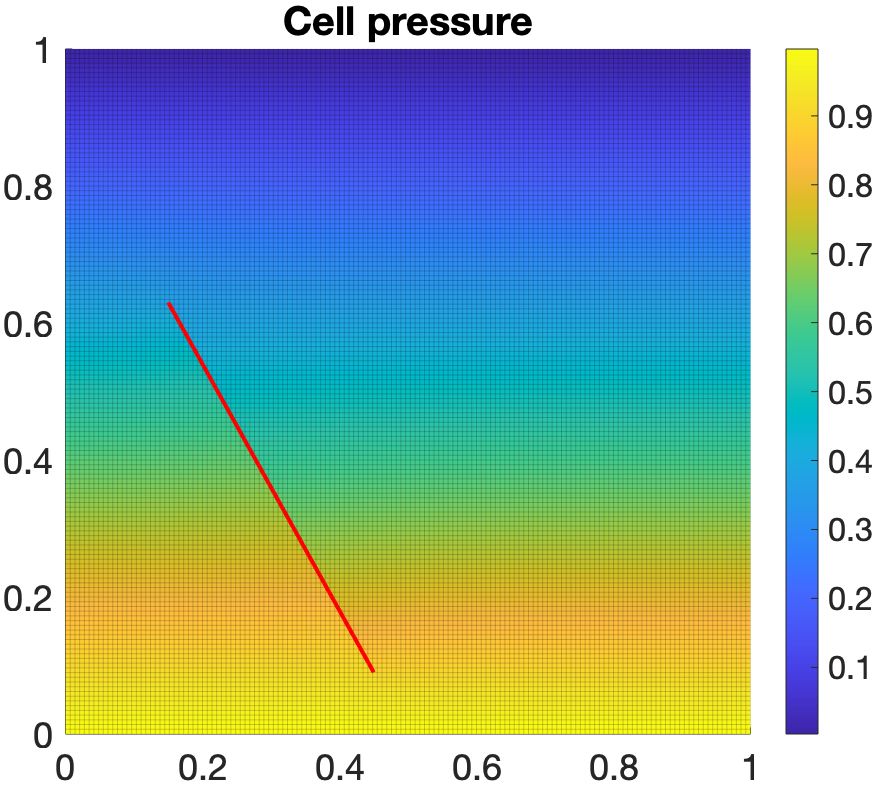}
		\caption{}
		\label{fig:oblique_imperm_1e-8_press_pEDFM_145x145}
	\end{subfigure}
	\hfill
	\begin{subfigure}[b]{0.495\textwidth}
		\centering
		\includegraphics[width=\textwidth]{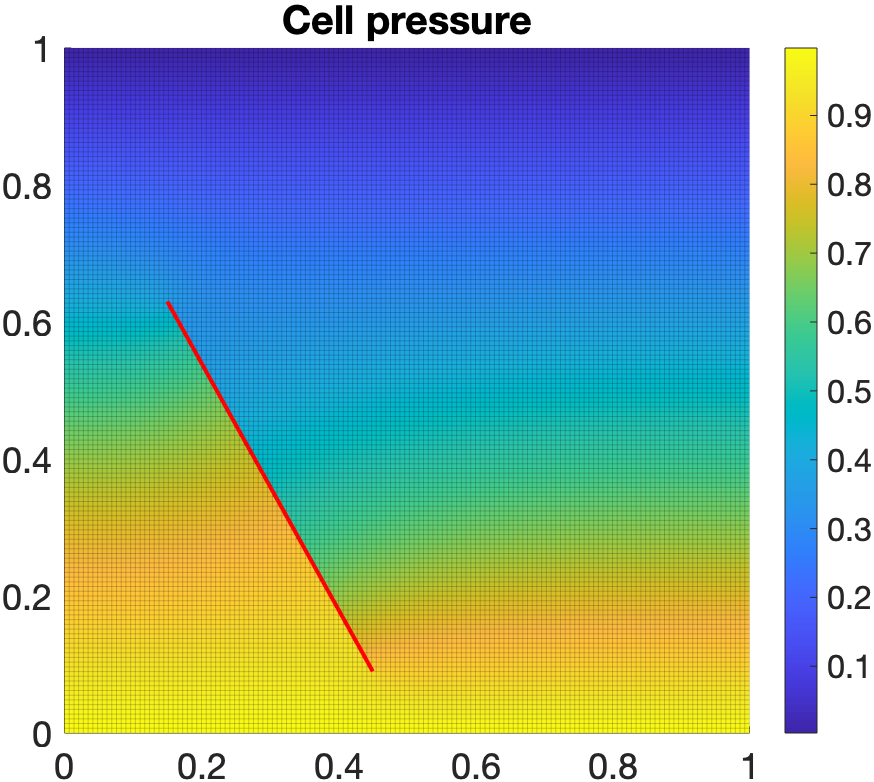}
		\caption{}
		\label{fig:oblique_imperm_1e-8_press_pEDFM_i_145x145}
	\end{subfigure}
	\vskip\baselineskip
	\begin{subfigure}[b]{0.495\textwidth}
		\centering
		\includegraphics[width=\textwidth]{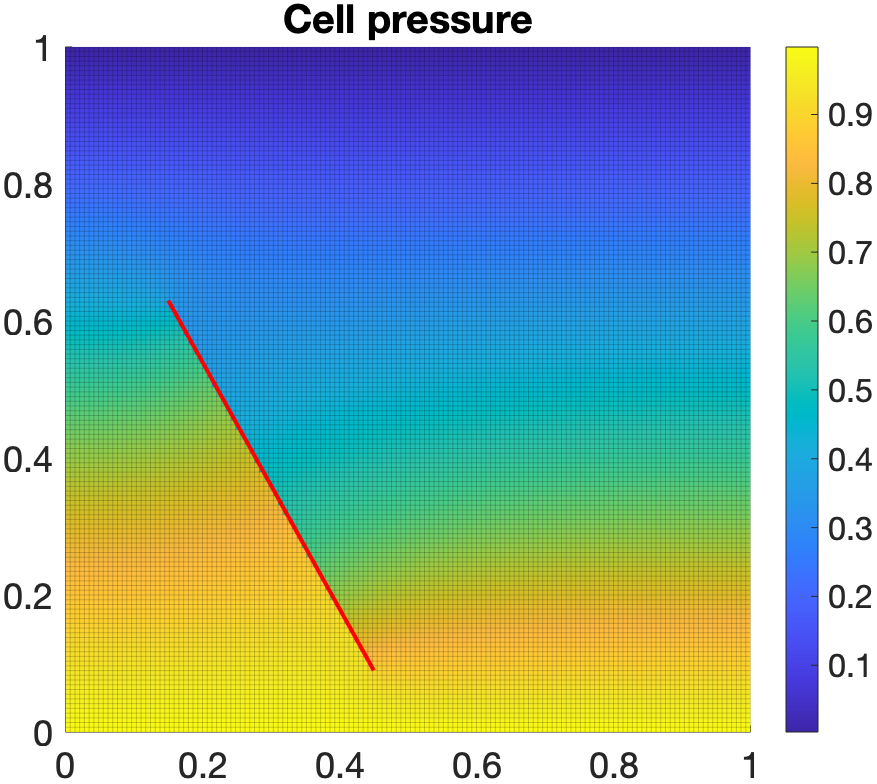}
		\caption{}
		\label{fig:oblique_imperm_1e-8_press_LEDFM_145x145}
	\end{subfigure}
	\hfill
	\begin{subfigure}[b]{0.495\textwidth}
		\centering
		\includegraphics[width=\textwidth]{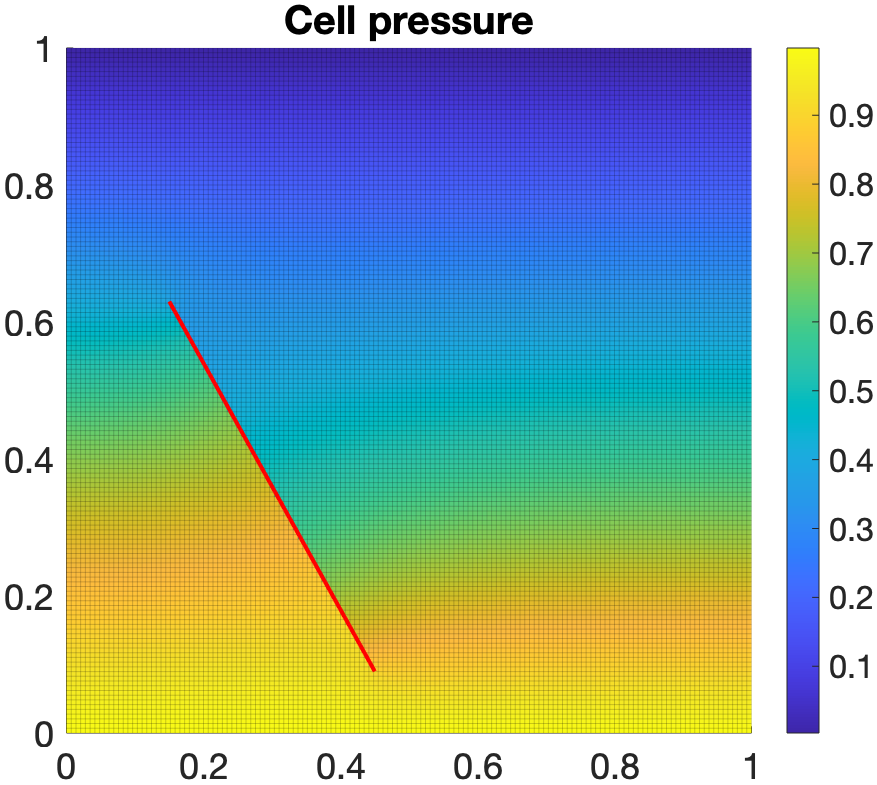}
		\caption{}
		\label{fig:oblique_imperm_1e-8_press_LEDFM_MSFV_145x145}
	\end{subfigure}
	\caption{Impermeable oblique fracture. Scenario $k_f/k_m = 10^{-8}$. (a) MRST pEDFM solution on $145\times145$ grid. (b) updated pEDFM solution on $145\times145$ grid. (c) LEDFM solution on $145\times145$ grid. (d) LEDFM+MSFV solution on $145\times145$ grid.}
	\label{fig:oblique_imperm_1e-8_press_145x145}
\end{figure}

\begin{figure}
	\centering
	\begin{subfigure}[b]{0.495\textwidth}
		\centering
		\includegraphics[width=\textwidth]{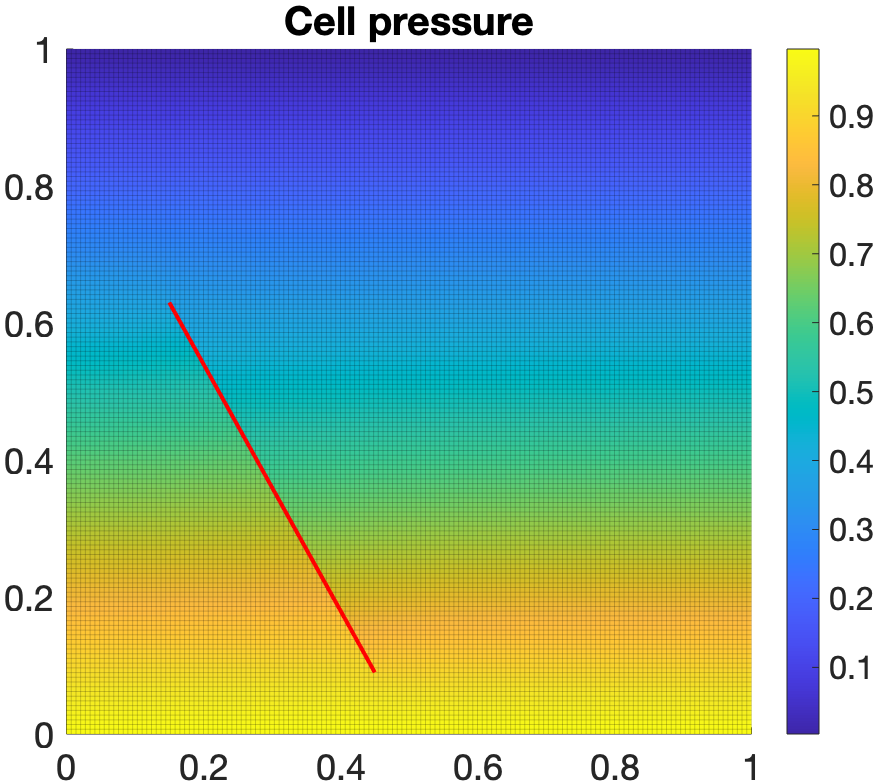}
		\caption{}
		\label{fig:oblique_imperm_1e-4_press_pEDFM_145x145}
	\end{subfigure}
	\hfill
	\begin{subfigure}[b]{0.495\textwidth}
		\centering
		\includegraphics[width=\textwidth]{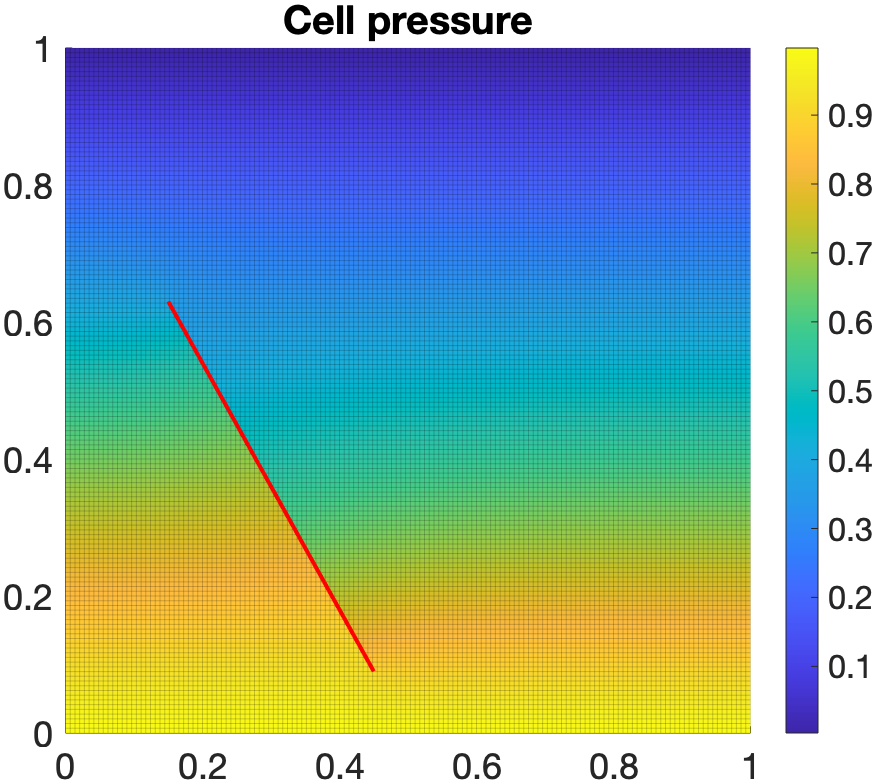}
		\caption{}
		\label{fig:oblique_imperm_1e-4_press_pEDFM_i_145x145}
	\end{subfigure}
	\vskip\baselineskip
	\begin{subfigure}[b]{0.495\textwidth}
		\centering
		\includegraphics[width=\textwidth]{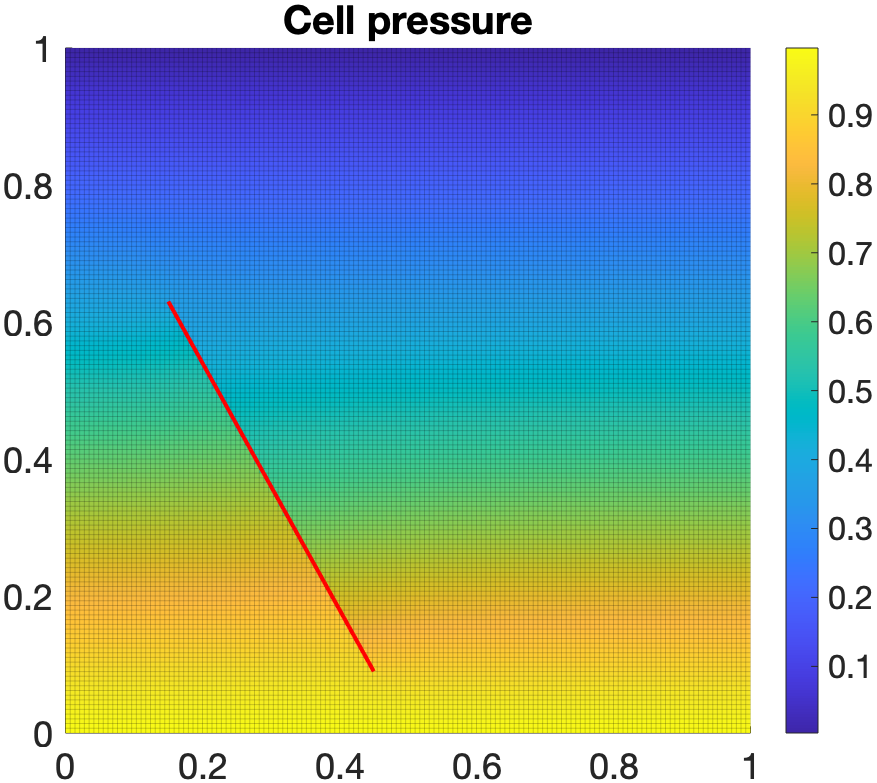}
		\caption{}
		\label{fig:oblique_imperm_1e-4_press_LEDFM_145x145}
	\end{subfigure}
	\hfill
	\begin{subfigure}[b]{0.495\textwidth}
		\centering
		\includegraphics[width=\textwidth]{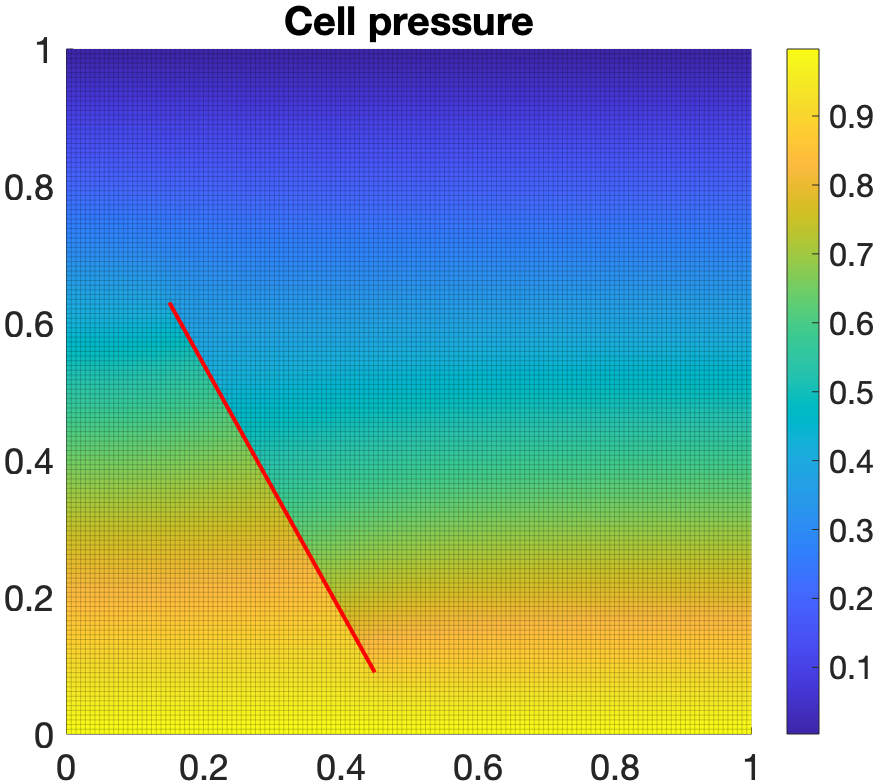}
		\caption{}
		\label{fig:oblique_imperm_1e-4_press_LEDFM_MSFV_145x145}
	\end{subfigure}
	\caption{Impermeable oblique fracture. Scenario $k_f/k_m = 10^{-4}$. (a) MRST pEDFM solution on $145\times145$ grid. (b) updated pEDFM solution on $145\times145$ grid. (c) LEDFM solution on $145\times145$ grid. (d) LEDFM+MSFV solution on $145\times145$ grid.}
	\label{fig:oblique_imperm_1e-4_press_145x145}
\end{figure}

Fig.~\ref{fig:oblique_imperm_1e-8_press_145x145} depicts the pressure for a $145\times145$ grid for all the competing methods: note that the pressure jump across the fracture is not correctly captured with the MRST pEDFM, contrarily to the other methods, as if the presence of the impermeable fracture is not accounted for, especially with fine grids, as it happens in the classical EDFM.

Fig.~\ref{fig:oblique_imperm_1e-4_press} shows the reference solution in the second scenario in which $k_f/k_m = 10^{-4}$, while in Fig.~\ref{fig:oblique_imperm_1e-4_conv_i} the corresponding convergence plot for all the methods is reported.
The same grids used in the first scenario are employed.
This time the standard version of LEDFM shows a poor convergence rate with respect to the other methods, apart of course from the classical EDFM and the MRST pEDFM, the latter showing the same behaviour observed in the high contrast case. 
We believe that, in this moderate contrast case, full-tensor effects due to the presence of the fracture in conjunction with non negligible flow velocities in the vicinity of it could lead to inconsistencies with a TPFA approximation.
Indeed, the Cartesian background grid is $\vec{K}$-orthogonal only if $\vec{K}$ is diagonal whereas fractures create different preferential directions. 
We thus applied the MSFV version of the method obtaining an improvement of the convergence rate with respect to the standard version of the method. In particular, for coarser grids it provides very similar results to those of the pEDFM, whose updated version overall performs slightly better than the MSFV version of the local method. However, the trend is still not optimal for fine grids: as the grid is refined the error tends to saturate reaching an almost constant value.

Fig.~\ref{fig:oblique_imperm_1e-4_press_145x145} depicts the pressure maps obtained for the second scenario in the case of a $145\times145$ grid for all the competing methods.
The standard LEDFM map, as well as the MRST pEDFM one, are substantially different from the reference solution, where the pressure jump across the fracture is more visible.
In the case of the multiscale modification of LEDFM and updated pEDFM, instead, the pressure distribution is represented more accurately than with the previous methods.

\subsection{Tracer Transport Test}

After verifying the convergence of the local method in the incompressible single-phase flow case, we want to use the corresponding velocity field to transport a tracer through a fractured porous medium. This allows us to test the accuracy of the velocity field.

Denoting the tracer concentration with the scalar quantity $c$, its advective transport through a porous medium $\Omega$ is described by the conservation equation:
\begin{equation} \label{eq:tracer_tr_eqn}
	\Phi \frac{\partial c}{\partial t} + \nabla \cdot (c \vec{u}) = 0 \qquad \textnormal{in } \Omega, 
\end{equation}
where $\Phi$ is the porosity and $\vec{u} = - (\vec{K}/ \mu) \nabla p$ the velocity field obtained from the single-phase flow problem.
We consider Dirichlet boundary conditions on the inflow boundary:
\begin{equation} \label{eq:tracer_tr_BCs}
	c = c_D \quad \textnormal{on } \Gamma^D_c, 
	\qquad \quad
	\Gamma^D_c = \{ \vec{x} \in \partial \Omega \colon \vec{u} \cdot \vec{n} < 0 \},
\end{equation}
where $c_D$ is the prescribed concentration value on the inflow boundary $\Gamma^D_c$ and $\vec{n}$ the unit outward normal vector to it.

The transport term is approximated with the standard first-order upwind scheme and the backward Euler method with a fixed time step is used for time discretization. The flow and transport problems are solved both with the classic and local versions of the embedded method on a $51 \times 51$ Cartesian grid and the results are compared to those obtained with a conforming method, already mentioned in Subsection~\ref{ssec:discretization}, which uses triangular grids of  $2594$ elements and the MPFA scheme.
In-house codes have been written to implement the embedded methods to solve both the single-phase flow and tracer transport problem. MRST was used, instead, to simulate both problems with the conforming method.

Let us now consider a square domain $\Omega = (0,100)^2$ crossed by a fracture $\gamma$, represented by the line segment joining the points $\vec{s}_1 = [25, 24]^\top$ and $\vec{s}_2 = [75, 74]^\top$ and characterized by a constant aperture $d = 10^{-2}$. The flow is driven upwards by imposing $p = 1.1 \cdot 10^6$ on the bottom boundary and $p = 10^6$ on the top boundary, while the remaining part of the boundary is impermeable. Hence, the inflow boundary is the bottom one and a tracer concentration of $c_D = 5.8 \cdot 10^{-5}$ is imposed here.
$\mu = 2.8 \cdot 10^{-4}$, $\Phi = 0.15$ in the whole domain and the matrix and fracture permeability tensors are homogeneous and isotropic with values $k_m = 10^{-13}$ and $k_f = 10^{-9}$, respectively, so that the fracture is more conductive than the porous matrix. The final simulation time is $5 \cdot 10^7$ discretized with $100$ time steps.

\begin{figure}[!h]
	\centering
	\begin{subfigure}[b]{0.3\textwidth}
		\centering
		\includegraphics[width=\textwidth]{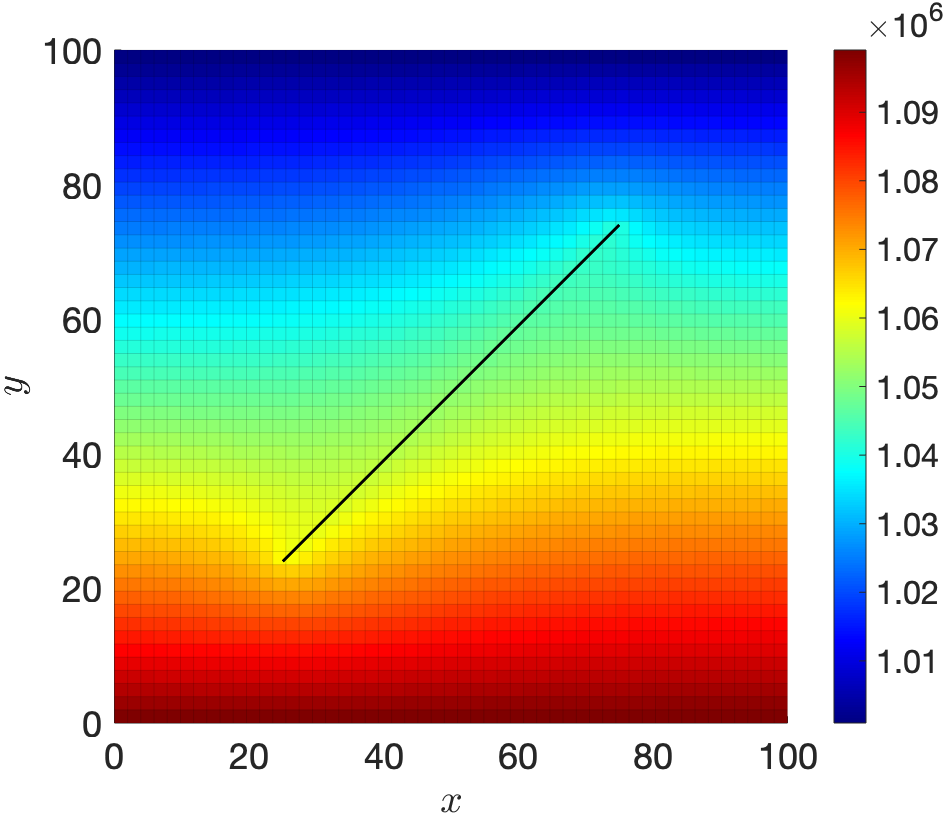}
		\caption{Classic EDFM.}
		\label{fig:tracer_tr_p_EDFM_51x51_classic}
	\end{subfigure}
	\begin{subfigure}[b]{0.3\textwidth}
		\centering
		\includegraphics[width=\textwidth]{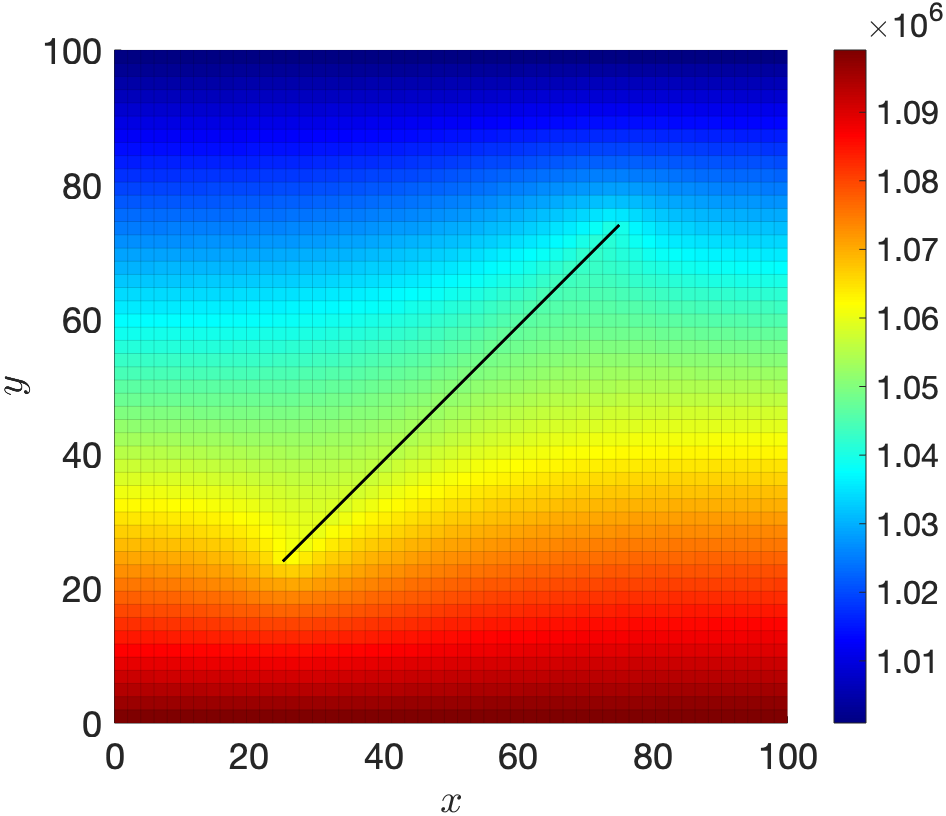}
		\caption{Local EDFM.}
		\label{fig:tracer_tr_p_EDFM_51x51_local}
	\end{subfigure}
	\begin{subfigure}[b]{0.3\textwidth}
		\centering
		\includegraphics[width=\textwidth]{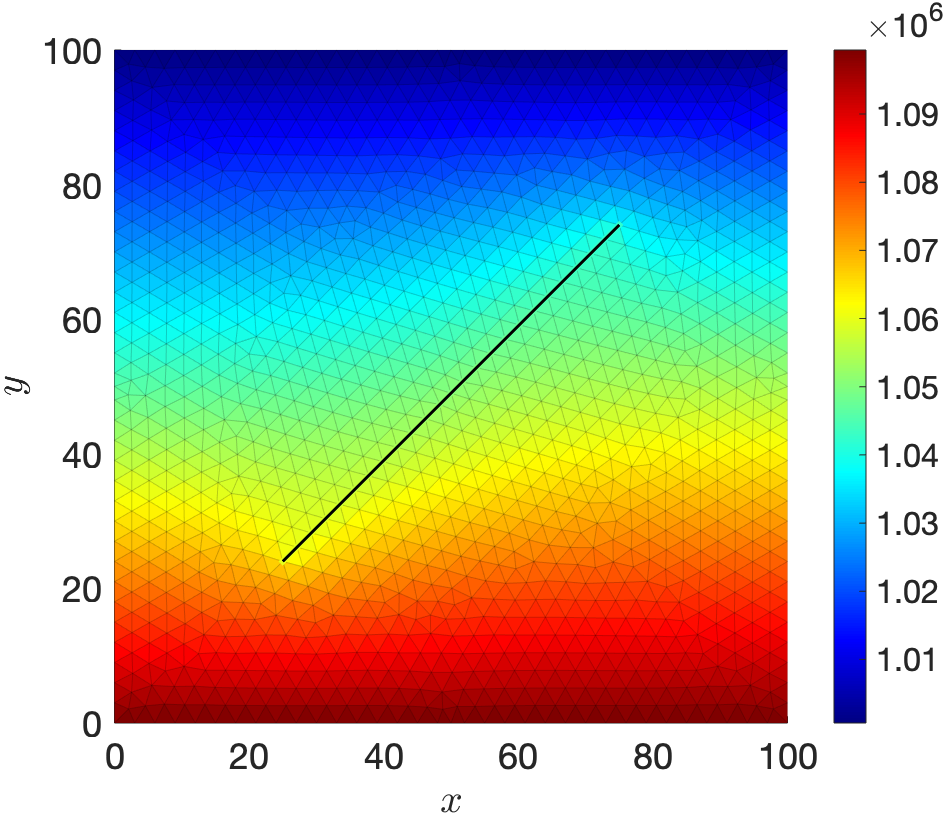}
		\caption{Conforming DFM.}
		\label{fig:tracer_tr_p_DFM_h3}
	\end{subfigure}
	\caption{Tracer transport test. Pressure field for all the considered methods.}
	\label{fig:tracer_tr_p}
\end{figure}

In Fig.~\ref{fig:tracer_tr_p} we observe that the pressure fields computed with all the considered methods are very similar to each other.

\begin{figure}[!h]
	\centering
	\begin{subfigure}[b]{0.3\textwidth}
		\centering
		\includegraphics[width=\textwidth]{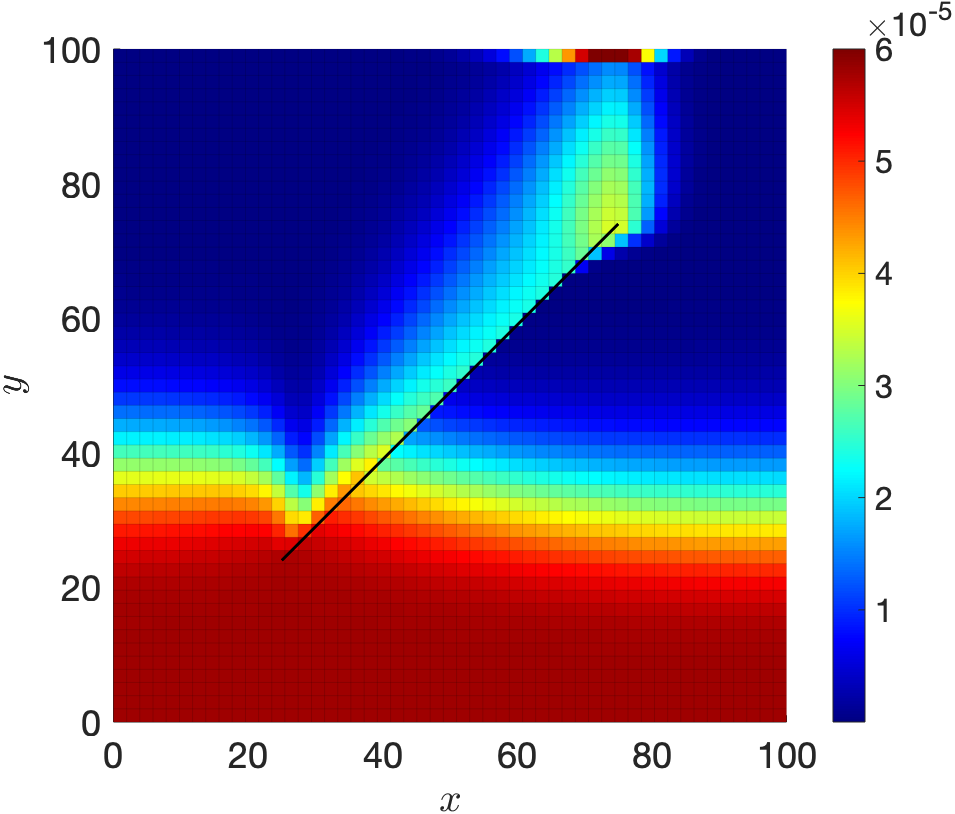}
		\caption{Classic EDFM.}
		\label{fig:tracer_tr_c_EDFM_51x51_classic}
	\end{subfigure}
	\begin{subfigure}[b]{0.3\textwidth}
		\centering
		\includegraphics[width=\textwidth]{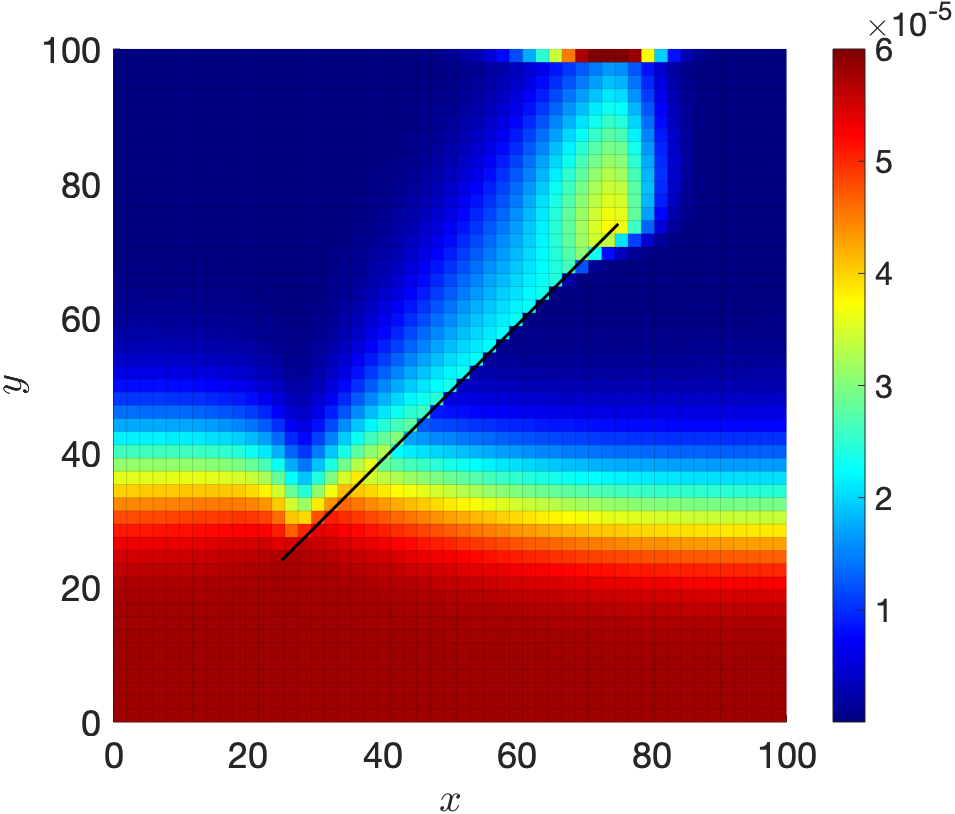}
		\caption{Local EDFM.}
		\label{fig:tracer_tr_c_EDFM_51x51_local}
	\end{subfigure}
	\begin{subfigure}[b]{0.3\textwidth}
		\centering
		\includegraphics[width=\textwidth]{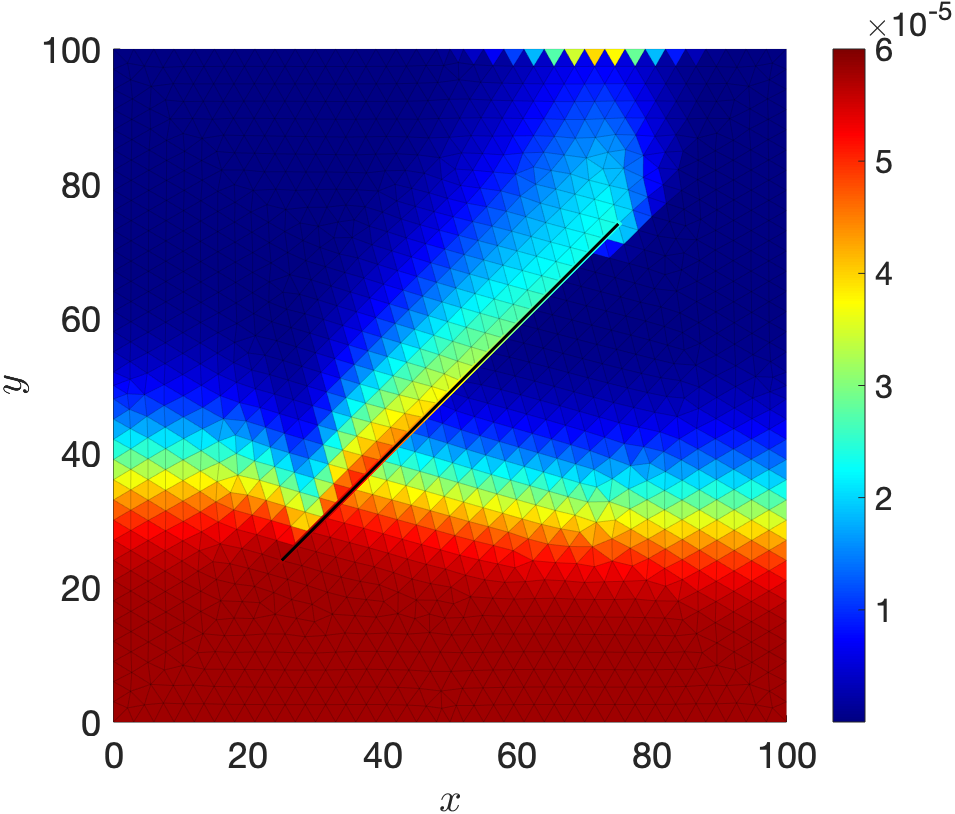}
		\caption{Conforming DFM.}
		\label{fig:tracer_tr_c_DFM_h3}
	\end{subfigure}
	\caption{Tracer transport test. Tracer concentration field for all the considered methods at time $t = 1.5 \cdot 10^7$.}
	\label{fig:tracer_tr_c}
\end{figure} 

\begin{figure}[!h]
	\centering
	\includegraphics[width=.4\textwidth]{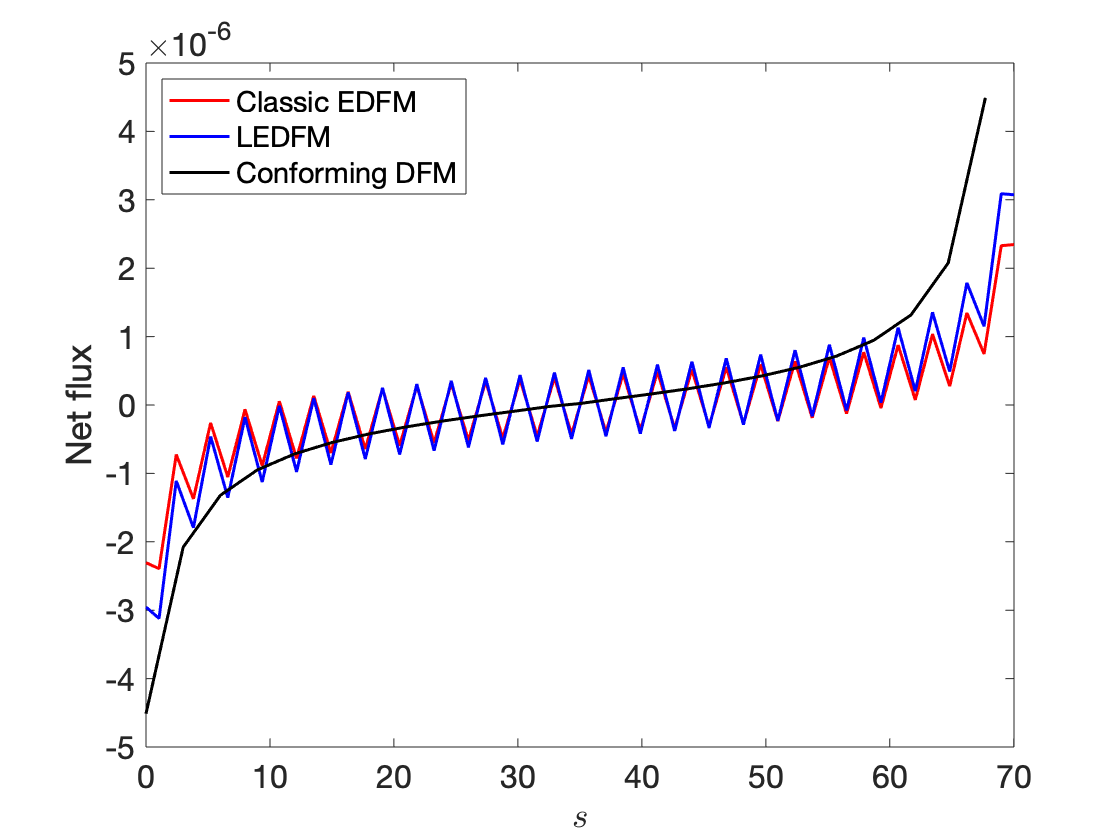}
	\caption{Tracer transport test. Net flux exiting the fracture for all the considered methods. $s$ is the arc length variable going from the lower to the upper fracture tip.}
	\label{fig:tracer_tr_netFlux}
\end{figure}

In Fig.~\ref{fig:tracer_tr_c} the tracer concentration distribution is depicted  for both the embedded and conforming methods at time $t = 1.5 \cdot 10^7$.
Similar results are obtained for the two embedded methods, but small differences are observed near to the upper fracture tip, where higher concentration values are attained with LEDFM compared to the classic version. This is believed to be the consequence of the slightly higher flux values obtained with LEDFM near fracture tips. 
The conforming DFM concentration distribution, instead, is quite different from the previous ones, especially in the near-fracture region: this is believed to happen due to differences in the gridding strategy and matrix-fracture flow modelling between conforming and embedded methods. Indeed, in 
Fig.~\ref{fig:tracer_tr_netFlux} we can observe the net flux exiting the fracture in the embedded and the conforming methods, computed as the difference between the flux exiting the fracture and directed towards the upper cells and the flux entering the fracture from the lower cells for the conforming method, while in the case of embedded methods it is represented by the flow directed from the fracture to the corresponding cut cell, computed from the matrix-fracture flow formula~\eqref{eq:Fmf}.
The average trend  (excluding the oscillations due to grid effects) is similar for the two methods and it shows that the flux is entering the fracture from the lower tip to approximately half of the fracture length and exiting the fracture from the middle of it to the upper tip. 
However net flux values are greater in magnitude for the conforming method near to the fracture tips with respect to embedded methods, justifying the different concentration maps.

\section{Conclusions} \label{sec:conclusions}

This paper presents a novel embedded method for flow simulation in fractured porous media, namely the \emph{Local Embedded Discrete Fracture Model} (LEDFM). Its formulation requires the solution of local problems, whose definitions are inspired from flow-based upscaling methods, for the computation of improved transmissibility coefficients.

The performance of the new method was compared to other well-known embedded methods, namely the \emph{Embedded Discrete Fracture Model} (EDFM) and the \emph{Projection-based Embedded Discrete Fracture Model} (pEDFM).
The results showed that LEDFM is able to capture the effect of the presence of both permeable and impermeable fractures with any orientation, even though, when dealing with impermeable fractures having intermediate conductivity contrasts, and not parallel to the underlying grid, a degradation of the convergence was observed for fine grids.

To improve the performance of the local method  a modified version of LEDFM replacing the M-M local problems with the \emph{Multiscale Finite Volume Method} (MSFV) was proposed.
The suggestion for the adoption of the MSFV method stems from the observation that full-tensor effects are accurately represented with MPFA methods, and fractures can be seen as objects that introduce full-tensor effects since they give preferential paths to the flow.
Indeed, using a MSFV method is equivalent to adopting a MPFA scheme at the coarse grid level, while the solution of M-M local problems on cells pairs corresponds to the adoption of a TPFA scheme.
Better results are obtained with the multiscale modification of the local method for tilted impermeable fractures having intermediate conductivity contrasts. However, some loss of convergence is still experienced in these cases.

Results similar to those obtained with the local method were observed with pEDFM, provided that the correct transmissibility expressions are used and that fracture cells projections are computed properly, as pointed out in Appendix~\ref{sec:warn_pEDFM}.

The local embedded method applicability was also verified for a tracer transport problem, where the velocity field computed in the single-phase flow case was used to transport the tracer through a porous medium cut by a conductive fracture.
The results were compared both with the classic EDFM and a conforming method.
The net flux exiting the fracture showed oscillations for the embedded methods due to grid effects, but the average trends were similar to that of the conforming method.
Moreover, LEDFM net flux values were closer to those of the conforming method compared to the classic EDFM, hence showing an improved matrix-fracture coupling.

In conclusion, LEDFM is an accurate, although expensive, model. Indeed,  it requires the solution of many local problems. In order to reduce the computational burden of the method, neural networks could be used in an offline case and trained to compute transmissibilities for different fractures geometries and permeability contrasts. This extension will be the subject of a forthcoming paper.

Possible directions for future work include the possibility of allowing the presence of more than one fracture in the local domain (possibly intersecting), the handling of the limit case of fractures lying on the local domain boundary and the extension of the method to the three-dimensional case.
Furthermore, the performance of the method in the case of impermeable intermediate conductivity contrasts should be improved. To do that, first the performances of the multiscale modification of the local method need to be further investigated, as well as the reasons behind the convergence degradation observed in some cases.  

\appendix
\section{1D Lower-Dimensional Discrete Fracture-Matrix Model} \label{sec:1D_DFMML}

This section provides the details of the 1D Lower-Dimensional Discrete Fracture-Matrix Model, that can be obtained, for the most part, as a particularization of the corresponding two- and three-dimensional model introduced, e.g., in \cite{Martin, Fumagalli3}. Such problems provide the boundary conditions for the local problem in the multiscale version of the local method.

Let us consider the case of a single ``fracture'' in the 1D porous medium domain, and incompressible single-phase flow. A ``fracture'' in this context should be intended as a narrow 1D region with respect to the porous domain size, characterized by a very different permeability value compared to the neighbouring matrix regions.

\begin{figure}[!t]
	\centering
	\includegraphics[width=.4\textwidth]{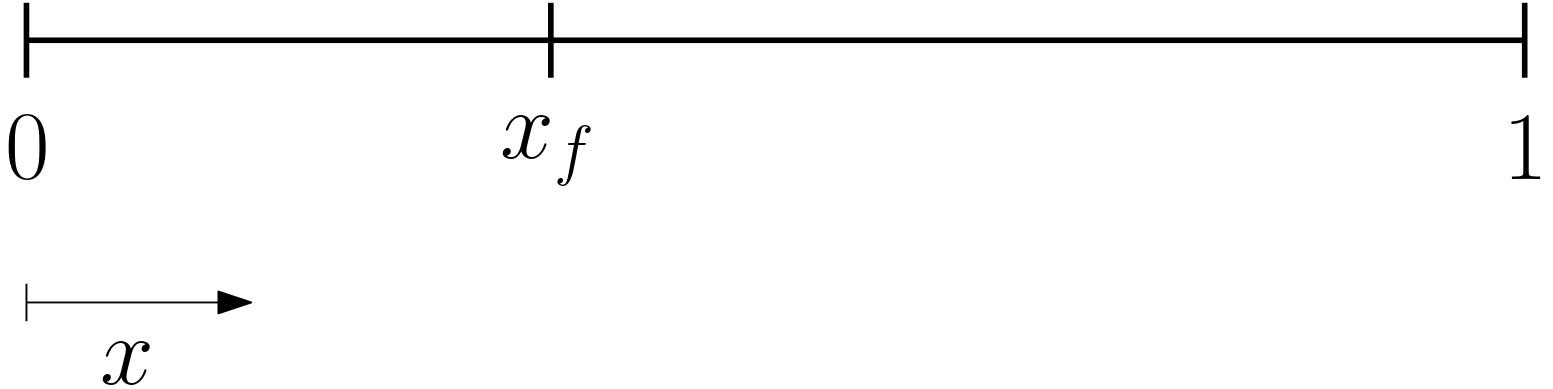}
	\caption{Domain $\Omega$ for the 1D Lower-Dimensional Discrete Fracture-Matrix Model. The fracture is positioned at $x = x_f$ and has a thickness $d \ll 1$.}
	\label{fig:1D_DFMM_L}
\end{figure}

The domain is $\Omega = (0, 1) \subset \mathbb{R}$, where the fracture is originally equidimensional, and corresponds to the region $\Omega_f = (x_f - d/2, x_f + d/2) \subset \Omega$. In the equi-dimensional case we have a Darcy problem with a different permeability in $\Omega_f$ and in the remaining part of the domain, $\Omega_m=\Omega \setminus \Omega_f$, and pressure and flux continuity at the interface.
We want to reduce the dimensionality of the fracture by collapsing it to a point to obtain a reduced formulation of the Darcy problem.
Assuming that the thickness $d$ of the fracture region is much smaller than the characteristic size of the porous medium domain we collapse the fracture region $\Omega_f$ to a point $\gamma = \{ x_f \}$.
The matrix region is now given by $\Omega_m = \Omega \setminus \gamma$, which can be split in two disjoint parts: $\Omega_{m_1} = (0, x_f)$ on the left and $\Omega_{m_2} = (x_f, 1)$ on the right of the fracture $\gamma$, as shown in Fig.~\ref{fig:1D_DFMM_L}.

Let us consider the domain $\Omega_{m_1}$. The 1D incompressible single-phase flow equation simply becomes 
\begin{equation} \label{eq:1D_left_eq}
	\frac{{\rm d}^2 p}{{\rm d} x^2} = 0 \qquad x \in \Omega_{m_1},
\end{equation}
since the permeability is constant in this region.
Solving equation~\eqref{eq:1D_left_eq} with the boundary conditions 
\begin{equation*} 
	\begin{cases}
		p(x=0) = p_0
		\\
		p(x = x_{f_-}) = p_L
	\end{cases}
	,
\end{equation*}
where $p_L$ is the pressure on the left side of the fracture, gives
\begin{equation} \label{eq:1D_p_matrix1}
	p(x) = \frac{p_L - p_0}{x_f} x + p_0
	\qquad
	x \in \Omega_{m_1}.
\end{equation}
Repeating the same steps for the domain $\Omega_{m_2}$, but with the boundary conditions
\begin{equation*} 
	\begin{cases}
		p(x= x_{f_+}) = p_R
		\\
		p(x = 1) = p_1
	\end{cases}
	,
\end{equation*}
where $p_R$ is the pressure on the right side of the fracture, gives
\begin{equation} \label{eq:1D_p_matrix2}
	p(x) = \frac{p_1 - p_R}{1 - x_f} x + \frac{p_R - p_1 x_f}{1 - x_f}
	\qquad
	x \in \Omega_{m_2}.
\end{equation}

$p_L$ and $p_R$ can be obtained by particularizing to the 1D case the interface conditions described in~\cite[Problem~3.2]{Fumagalli3} for the multi-dimensional case.
Since the fracture is 0D the net incoming flux should be null in the absence of source terms, therefore
\begin{equation} \label{eq:1D_red_frac_flow}
	\left( \frac{{\rm d} p}{{\rm d} x} \right)_L = \left( \frac{{\rm d} p}{{\rm d} x} \right)_R.
\end{equation}
Thus the pressure slopes will be identical to the left and right of the fracture. 
The second interface condition, instead, becomes
\begin{equation} \label{eq:1D_interf_avg_cond}
	-\frac{d}{2} \frac{k_m}{k_f} \left[  \left( \frac{{\rm d} p}{{\rm d} x} \right)_L + \left( \frac{{\rm d} p}{{\rm d} x} \right)_R \right] = p_L - p_R,
\end{equation}
where
\begin{equation*}
	\left( \frac{{\rm d} p}{{\rm d} x} \right)_L = \frac{p_L - p_0}{x_f},
	\qquad
	\left( \frac{{\rm d} p}{{\rm d} x} \right)_R = \frac{p_1 - p_R}{1 - x_f}.
\end{equation*}
Using these relations in equations~\eqref{eq:1D_red_frac_flow} and~\eqref{eq:1D_interf_avg_cond} we get
\begin{equation*}
	\begin{cases}
		\dfrac{p_L - p_0}{x_f} = \dfrac{p_1 - p_R}{1 - x_f}
		\\
		-\dfrac{d}{2} \dfrac{k_m}{k_f} \left( \dfrac{p_L - p_0}{x_f} + \dfrac{p_1 - p_R}{1 - x_f}  \right) = p_L - p_R
	\end{cases}
	,
\end{equation*}
which is a linear algebraic system giving the following expressions for $p_L$ and $p_R$:
\begin{align}
	p_L &= \frac{(1-x_f) k_f + d k_m}{k_f + d k_m} p_0 + \frac{x_f k_f}{k_f + d k_m} p_1, \label{eq:pL}
	\\
	p_R &= \frac{(1-x_f) k_f}{k_f + d k_m} p_0 + \frac{x_f k_f + d k_m}{k_f + d k_m} p_1. \label{eq:pR}
\end{align}

\begin{figure}[!t]
	\centering
	\includegraphics[width=.45\textwidth]{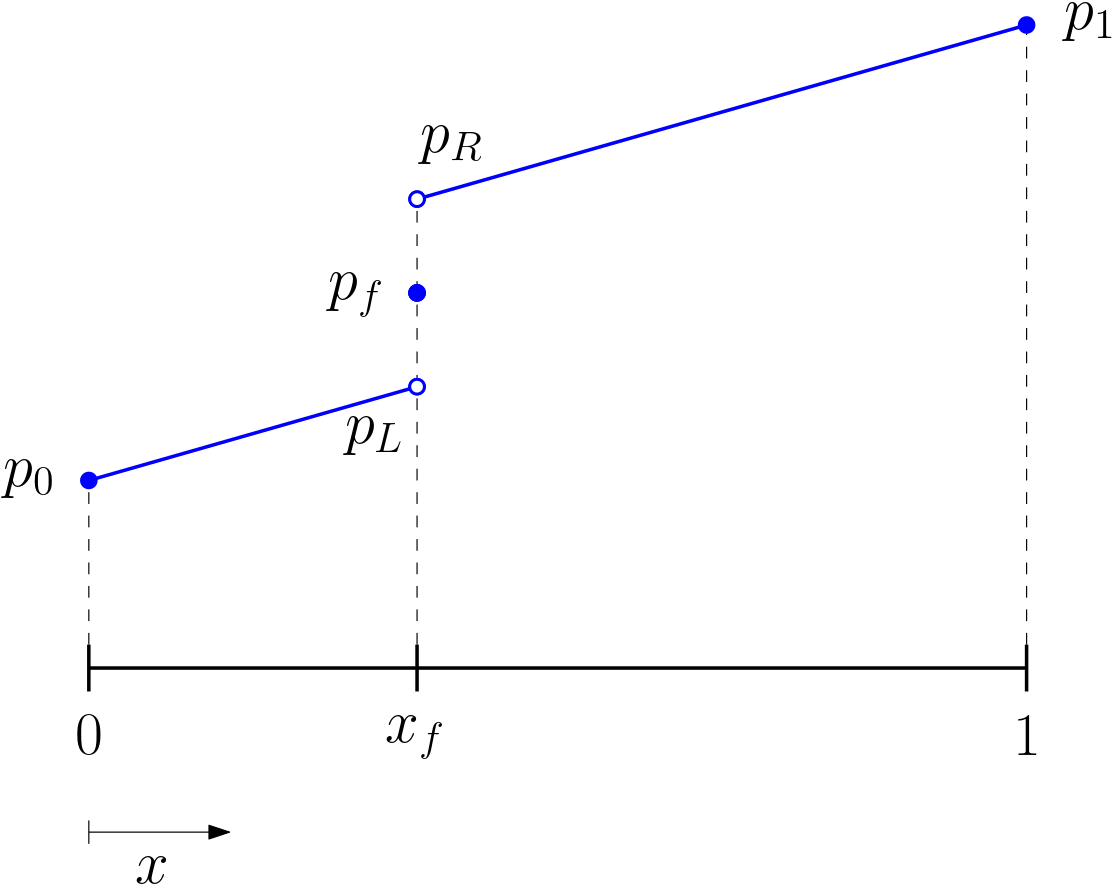}
	\caption{Example of a pressure profile $p(x)$ for the 1D Lower-Dimensional Discrete Fracture-Matrix Model. Case $p_1 > p_0$.}
	\label{fig:1D_DFMM_L_p}
\end{figure}

Finally, the pressure profile for the entire 1D porous medium domain is given by
\begin{align} \label{eq:1D_DFMM_L_p}
	p(x)
	=
	\begin{cases}
		\dfrac{k_f}{k_f + d k_m} (p_1 - p_0) x + p_0 &\qquad x \in \Omega_{m_1} = (0, x_f),
		\\
		\dfrac{2 (1-x_f) k_f + d k_m}{2 (k_f + d k_m)} p_0 + \dfrac{2 x_f k_f + d k_m}{2 (k_f + d k_m)} p_1 &\qquad x \in \gamma = \{ x_f \},
		\\
		\dfrac{k_f}{k_f + d k_m} (p_1 - p_0) x + \dfrac{k_f}{k_f + d k_m} p_0 + \dfrac{d k_m}{k_f + d k_m} &\qquad x \in \Omega_{m_2} = (x_f, 1).
	\end{cases}
\end{align}
A qualitative trend of the pressure is shown in Fig.~\ref{fig:1D_DFMM_L_p} for a case in which $p_1 > p_0$.

\section{Notes on the pEDFM implementation} \label{sec:warn_pEDFM}

Here we present some of the most important modifications made to the pEDFM code of MRST, described in \cite{Olorode, MRST2}, along with some general warnings about the implementation of the method itself. The changes were necessary to correctly simulate the examples proposed in Section~\ref{ssec:conv_tests}.

First of all, the code implemented in MRST employs the transmissibility expressions proposed in the original pEDFM paper~\cite{Tene}. To be precise, and using the notation introduced in Section~\ref{ssec:pEDFM}, the distances $d_{m_Hf}$ and $d_{m_Vf}$ are used in the non-neighbouring transmissibility formulae instead of the average ones $\langle d \rangle_{m_Hf}$ and $\langle d \rangle_{m_Vf}$.
However, a simple harmonic averaging of the permeabilities is used to compute the transmissibilities instead of a distance-weighted one: this has no effect on the modified matrix-matrix transmissibility expressions, but leads to differences in the matrix-fracture and non-neighbouring matrix-fracture transmissibility formulae with respect to the original version of the method. The latter can have a negative impact on the accuracy of the numerical solution if impermeable fractures are to be simulated, especially if they are thin, as explained in detail in Appendix~\ref{sec:1DpEDFM}.
For this reason, the transmissibility formulae~\eqref{eq:Tmf_jiang}--\eqref{eq:half_transm_jiang} proposed in~\cite{Jiang} have been implemented in MRST, which give accurate results in every situation, as also pointed out in Appendix~\ref{sec:1DpEDFM}.

Another issue detected in the code is the incomplete identification of the non-neighbouring matrix-fracture connections. Indeed, it may happen that some of them are missing and this negatively affects the quality of the numerical solution. This problem has also been solved. 

\begin{figure}[!b]
	\centering
	\begin{subfigure}[b]{0.495\textwidth}
		\centering
		\includegraphics[width=.55\textwidth]{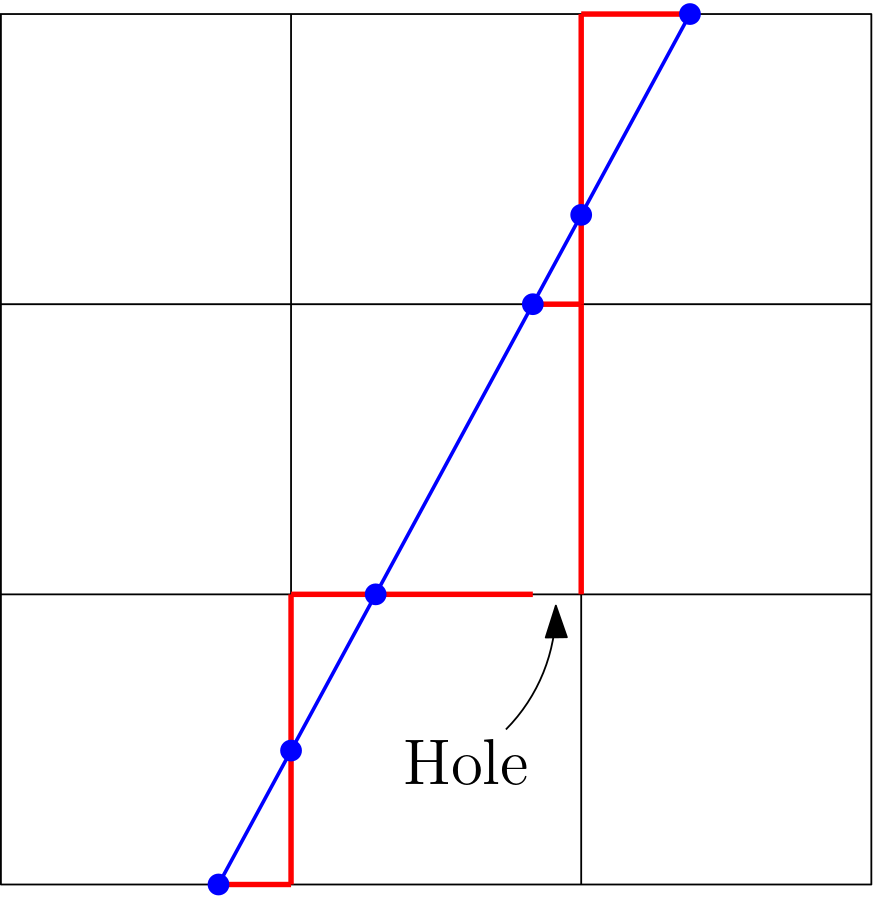}
		\caption{}
		\label{fig:pEDFM_hole}
	\end{subfigure}
	\hfill
	\begin{subfigure}[b]{0.495\textwidth}
		\centering
		\includegraphics[width=.55\textwidth]{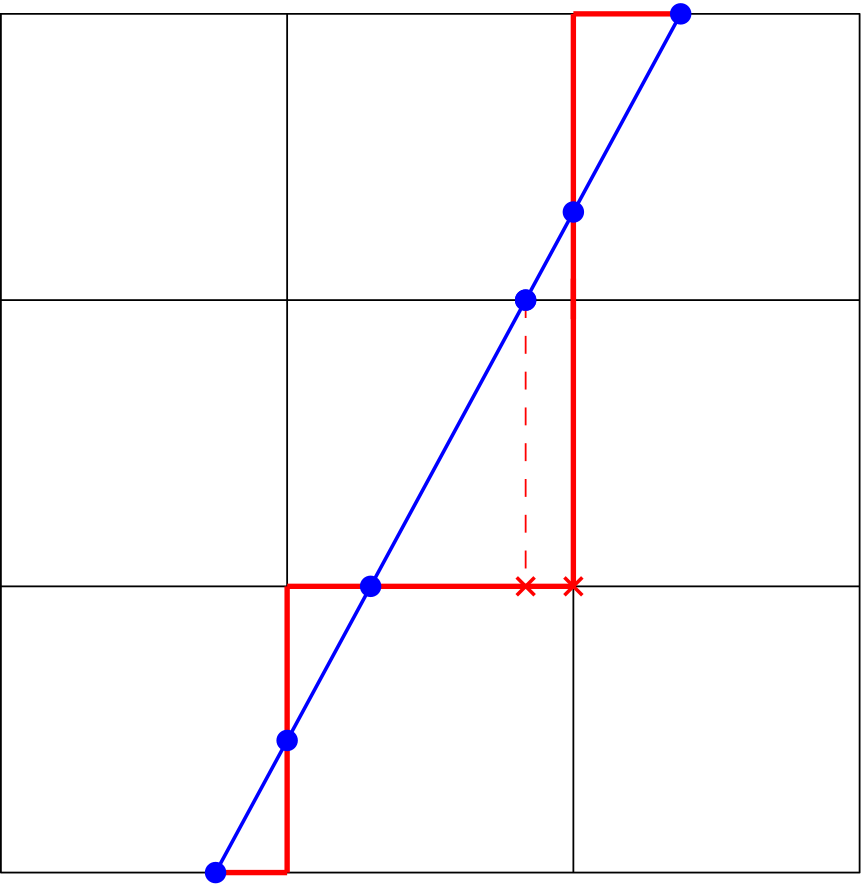}
		\caption{}
		\label{fig:pEDFM_no_hole}
	\end{subfigure}
	\caption{Fracture projection on matrix grid faces. (a) Only the cells neighbouring a cut cell are considered in the procedure. A hole is created in the fracture projection path. (b) Cells both neighbouring and not a cut cell are considered in the procedure.}
	\label{fig:pEDFM_hole_no_hole}
\end{figure}

As stated in \cite{Tene}, the pEDFM method requires the construction of a continuous projection path of the fracture cells on the matrix grid faces in order for it to work. 
Note that the construction of this path is, in general, not straighforward: indeed, as can be seen in Fig.~\ref{fig:pEDFM_hole_no_hole}, there are cases in which it is not sufficient to consider only the cells neighbouring a matrix cell cut by a fracture to obtain a continuous projection path of the fracture itself.
Not considering this possibility, as it happens for the original implementation of the pEDFM code in MRST, results in having holes in the fracture projection path in some cases, which results in leakage across the fracture.
If the fracture is impermeable, the holes greatly lower its barrier effect on the fluid flow leading to very inaccurate pressure solutions, as described in Section~\ref{ssec:imperm_oblique_frac}.

\begin{remark}
 There are special cases in which the straightforward algorithm implemented in MRST produces a correct result, i.e. the path is continuous even if one considers only neighbouring matrix cells. If we consider Cartesian grids, such cases correspond to: horizontal fractures, vertical fractures, and fractures forming an angle $\alpha = \arctan(h_y/h_x)$ with the horizontal direction, where $h_x$ and $h_y$ are the discretization steps in the $x$ and $y$ directions, respectively.
\end{remark}

The code has been fixed to obtain continuous fracture projection paths for any fracture inclination and grid refinement level, but it must be observed that this was possible thanks to the simplicity of the structure of the Cartesian grid. More complex grids would entail more difficulties in the implementation of the method.
\section{Comparison of pEDFM formulations in a 1D setting} \label{sec:1DpEDFM}

In this section we want to compare the implementation of pEDFM in MRST, denoted as ''MRST pEDFM'', and the  updated version of the code where the formulation proposed in~\cite{Jiang} is adopted, denoted as ''updated pEDFM''. The two formulations differ in the way some of the transmissibilities are defined, as explained in Appendix~\ref{sec:warn_pEDFM}. To this aim we consider, for the sake of simplicity, a 1D incompressible single-phase flow problem in presence of a single fracture with pEDFM. $\Omega = (0, 1) \subset \mathbb{R}$ and the fracture corresponds to a point $\gamma = \{ x_f \}$, with $x_f = 1/2$. The matrix domain is then given by $\Omega_m = \Omega \setminus \gamma$. Let us denote with $k_m$ and $k_f$ the porous matrix and fracture permeability, respectively, and consider a null left boundary pressure and a unit right boundary pressure, i.e. $p_0 = 0$ and $p_1 = 1$.

Far from being a proof, this comparison will however show that the MRST version cannot simulate well some impermeable fracture scenarios, especially if the fracture is thin. On the other hand, the updated version will prove to be accurate in any situation.

Let us now discretize the 1D domain $(0, 1)$ in $N$ subintervals of width $h = 1/N$. The nodes of the discretization are $x_i = i h$, with $i =  0, \ldots, N$, while the grid cells and their centroids are defined as $C_i = \left[ x_{i-1}, x_i \right]$ and $x_{C_i} = \left(x_{i-1} + x_i \right)/2$, respectively, with $i = 1, \ldots, N$.
$N$ is assumed to be an odd number so that the fracture does not coincide with a grid node.

The matrix-matrix transmissibility $T_{m_1m_2}$ and the modified matrix-matrix transmissibility $T_{m m_n}$, where $m_n$ is the non-neighbouring matrix cell, i.e. the matrix cell sharing with $m$ the grid point on which the fracture is projected, in this case are given by
\begin{equation*}
	T_{m_1m_2} = \frac{k_m}{h},
	\qquad
	T_{m m_n} = 0
\end{equation*}
for both pEDFM versions. Indeed, the two versions differ only in the expressions for the matrix-fracture and non-neighbouring matrix-fracture transmissibilities.
Note that $T_{m m_n}$ is always null in the one-dimensional case since the fracture projection always coincides with the entire interface shared by the cells $m$ and $m_n$, i.e. a point.

Particularizing to the 1D case the matrix-fracture and non-neighbouring matrix-fracture transmissibility expressions of the MRST pEDFM we obtain
\begin{equation} \label{eq:transm_MRST}
	T_{mf}^M = \frac{8}{h} \frac{k_m k_f}{k_m + k_f},
	\qquad
	T_{m_n f}^M = \frac{2}{h} \frac{k_m k_f}{k_m + k_f},
\end{equation}
where the superscript $M$ stands for ``MRST''.

Proceeding in a similar manner, from~\eqref{eq:Tmf_jiang},~\eqref{eq:Tmnf_jiang} and~\eqref{eq:half_transm_jiang} we obtain the following transmissibility formulae for the updated version
\begin{equation} \label{eq:transm_updated}
	T_{mf}^U = \frac{4 k_m k_f}{2 k_m d + k_f h},
	\qquad
	T_{m_n f}^U = \frac{2 k_m k_f}{k_m d + 2 k_f h},
\end{equation}
where the superscript $U$ stands for ``Updated''.

We now compare the results obtained using both versions with the corresponding analytical solution~\eqref{eq:1D_DFMM_L_p} of the 1D Lower-Dimensional Discrete Fracture-Matrix Model.
The fracture aperture and matrix permeability are set to $d = 10^{-4}$, $k_m = 1$, respectively, and we examine three different permeability contrasts $R_k \vcentcolon= k_f / k_m$ between porous matrix and fracture domains: $R_k = 10^8$, $R_k = 10^{-8}$, and $R_k = 10^{-4}$.

\begin{figure}[!t]
	\centering
	\begin{subfigure}[b]{0.328\textwidth}
		\centering
		\includegraphics[width=\textwidth]{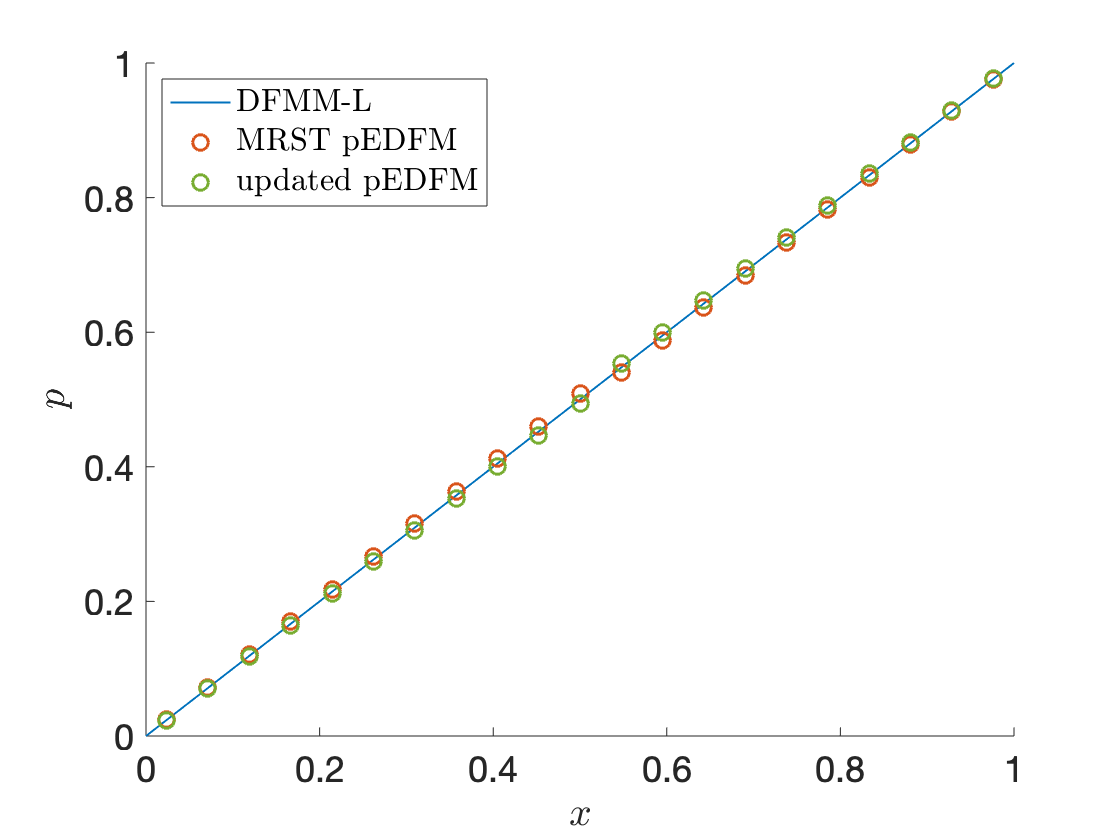}
		\caption{$R_k = 10^8$.}
		\label{fig:1D_pEDFM_1e+8}
	\end{subfigure}
	\hfill
	\begin{subfigure}[b]{0.328\textwidth}
		\centering
		\includegraphics[width=\textwidth]{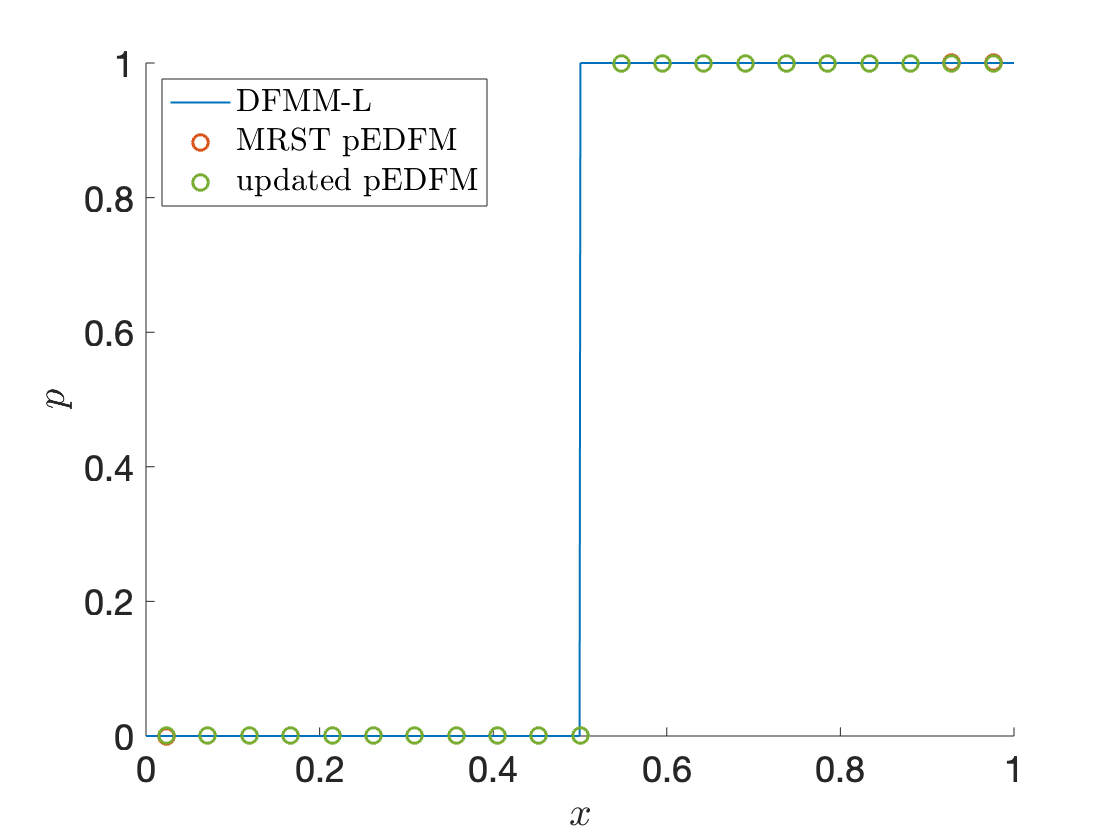}
		\caption{$R_k = 10^{-8}$.}
		\label{fig:1D_pEDFM_1e-8}
	\end{subfigure}
	\hfill
	\begin{subfigure}[b]{0.328\textwidth}
		\centering
		\includegraphics[width=\textwidth]{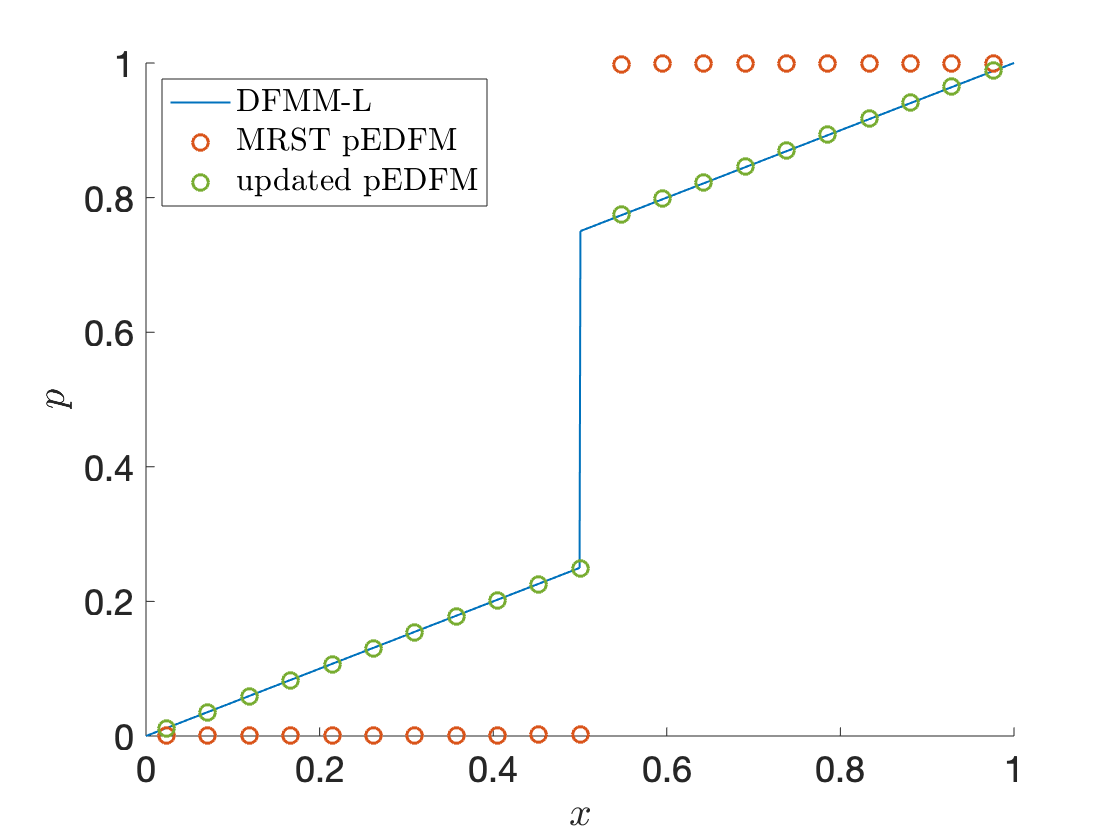}
		\caption{$R_k = 10^{-4}$.}
		\label{fig:1D_pEDFM_1e-4}
	\end{subfigure}
	\caption{MRST and updated pEDFM pressure solutions compared with the analytical solution of the 1D Lower-Dimensional Discrete Fracture-Matrix Model (DFMM-L) for different permeability contrast cases. The fracture aperture is set to $d = 10^{-4}$.}
	\label{fig:1D_pEDFM}
\end{figure}

Fig.~\ref{fig:1D_pEDFM} shows the numerical and analytical pressure solutions for the different permeability contrast cases, for $N = 21$.  
From Figs.~\ref{fig:1D_pEDFM_1e+8} and~\ref{fig:1D_pEDFM_1e-8} we note that both pEDFM versions provide results very close to the analytical ones in the highly permeable and impermeable cases.
However, in the intermediate impermeable case ($R_k = 10^{-4}$), depicted in Fig.~\ref{fig:1D_pEDFM_1e-4}, it is clear that the MRST pEDFM is not able to capture correctly the solution behaviour, contrarily to the updated version. In particular, the barrier effect of the fracture is overestimated, leading to a much higher pressure jump across the fracture than that of the analytical and updated pEDFM solutions.
This occurs because a simple harmonic averaging of the permeabilities is adopted in the matrix-fracture and non-neighbouring matrix-fracture transmissibility expressions of the MRST pEDFM, so that the fracture aperture is not accounted for.

To study the source of error, we determine the expressions of the numerical pressures for both the MRST and updated versions of pEDFM in a very coarse case, e.g. with $N = 5$ subintervals, solving exactly the linear system stemming from the finite volume discretization and then we compute the $L^1$ errors of the numerical solution with respect to the analytical one.
Although this is more an experiment rather than a formal proof, to make the comparison as fair as possible, the fracture position is shifted to $x_f = 0.6$ for the analytical solution, since in the pEDFM formulation the fracture is projected from $x_f = 0.5$ right onto $x = 0.6$, i.e. on the right interface of the matrix cell containing the fracture.
This is possible since the fracture position $x_f$ has no effect on the slope of the analytical linear pressure solution, as can be seen in~\eqref{eq:1D_DFMM_L_p}.

The linear system arising from the discretization can be written as
\begin{equation} \label{eq:1D_linsys}
	\mathbf{A} \widetilde{\mathbf{p}} = \mathbf{q},
\end{equation}
where $\widetilde{\mathbf{p}} = \left[\begin{array}{c;{2pt/2pt}c} \widetilde{\mathbf{p}}^\top_m & \widetilde{p}_f \end{array}\right]^\top \in \mathbb{R}^6$ is the unknown vector of pressures, $\widetilde{\mathbf{p}}_m \in \mathbb{R}^5$ being the vector of matrix pressures and $\widetilde{p}_f$ the fracture pressure. 
$\mathbf{A} \in \mathbb{R}^{6 \times 6}$ is the coefficient matrix and $\mathbf{q} \in \mathbb{R}^6$ the right hand side vector containing boundary contributions, and are given by
\begin{equation*}
	\mathbf{A} = 
	\left[
		\begin{array}{ccccc;{2pt/2pt}c}
			\frac{3 k_m}{h} & -\frac{k_m}{h} & 0 & 0 & 0 & 0\\
			-\frac{k_m}{h} & \frac{2 k_m}{h} &  -\frac{k_m}{h} & 0 & 0 & 0\\
			0 & -\frac{k_m}{h} & \frac{k_m}{h} + T_{mf} & 0 & 0 & -T_{mf}\\
			0 & 0 & 0 & \frac{k_m}{h} + T_{m_nf} & -\frac{k_m}{h} & -T_{m_nf}\\
			0 & 0 & 0 & -\frac{k_m}{h} & \frac{3 k_m}{h} & 0\\\hdashline[2pt/2pt]
			0 & 0 & -T_{mf} & -T_{m_nf} & 0 & T_{mf} + T_{m_nf}
		\end{array}
	\right],
	\qquad
	\mathbf{q} =
	\left[
		\begin{array}{c}
			\frac{2 k_m}{h} p_0\\
			0\\
			0\\
			0\\
			\frac{2 k_m}{h} p_1\\\hdashline[2pt/2pt]
			0
		\end{array}
	\right],
\end{equation*}
where the blocks contain the terms related to the porous matrix or pertaining to the fracture.

Setting $p_0 = 0$ and $p_1 = 1$ as boundary pressures, knowing the expressions for $T_{mf}$ and $T_{m_nf}$~\eqref{eq:transm_MRST}--\eqref{eq:transm_updated} for the MRST and updated pEDFM versions, respectively, and using the fact that $h = 1/5$ yields the following expressions for the matrix pressures once the linear system~\eqref{eq:1D_linsys} has been solved 

\begin{equation*}
	\widetilde{\mathbf{p}}_{m}^M = 
	\frac{1}{5 + 37 R_k} 
	\begin{bmatrix}
		4 R_k \\
		12 R_k \\
		20 R_k \\
		5 + 25 R_k  \\
		5 + 33 R_k
	\end{bmatrix},
	\qquad
	\widetilde{\mathbf{p}}_{m}^U = 
	\frac{1}{20 d + 21 R_k} 
	\begin{bmatrix}
		2 R_k \\
		6 R_k \\
		10 R_k \\
		20 d + 15 R_k \\
		20 d + 19 R_k
	\end{bmatrix},
\end{equation*}
where $\widetilde{\mathbf{p}}_{m}^M$ is the vector of matrix pressures of the MRST pEDFM and $\widetilde{\mathbf{p}}_{m}^U$ the one relative to the updated version.

Moreover, using~\eqref{eq:1D_DFMM_L_p}, we compute the vector of analytical pressures $\mathbf{p}_a$ at the grid cell centroids $x_{C_i}$:
\begin{equation*}
	\mathbf{p}_a =
	\frac{1}{10 (R_k + d)}
		\begin{bmatrix}
		R_k \\
		3 R_k \\
		5 R_k \\
		7 R_k + 10 d  \\
		9 R_k + 10 d
	\end{bmatrix}.
\end{equation*}

\begin{figure}[!t]
	\centering
	\includegraphics[width=.5\textwidth]{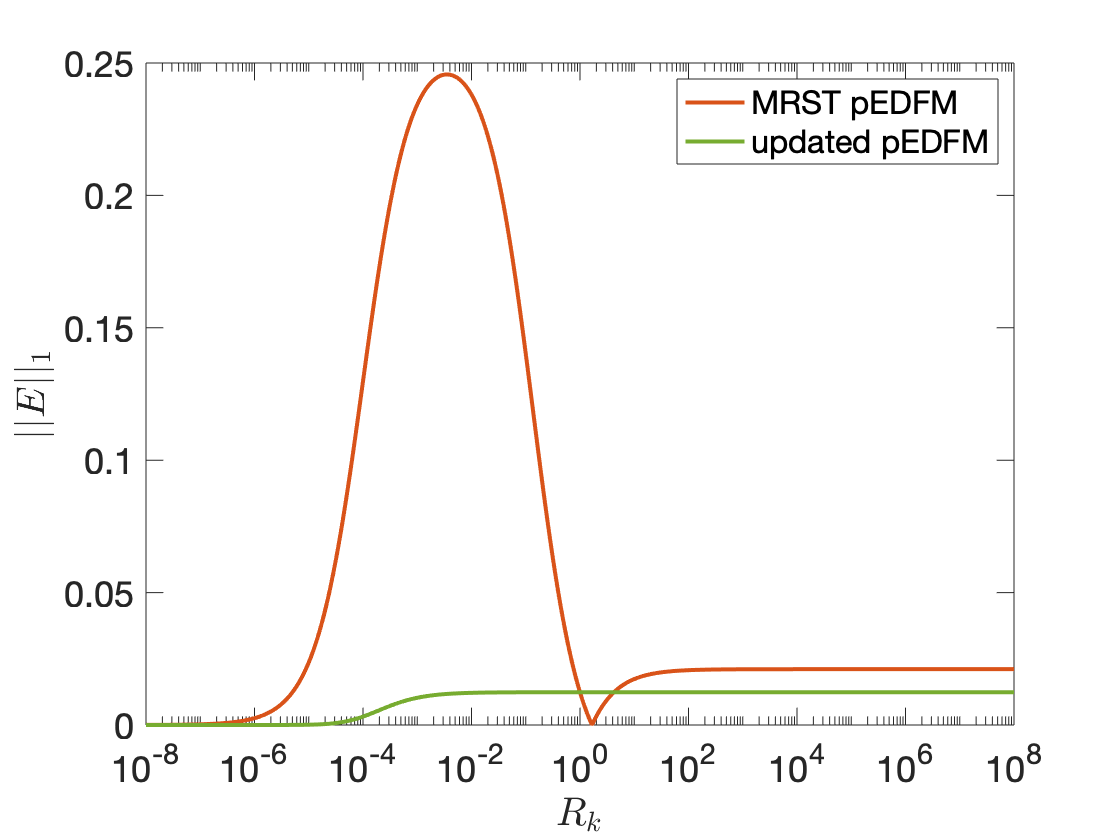}
	\caption{MRST and updated pEDFM $L^1$ errors for variable $R_k$. The fracture aperture is set to $d = 10^{-4}$.}
	\label{fig:1D_pEDFM_err}
\end{figure}

The $L^1$ error $|| E ||_1$ is defined as
\begin{equation*}
	|| E ||_1 = h \sum_{j=1}^{N} | E_j |
	\quad
	\textnormal{with }
	E_j = p_{a_j} - \widetilde{p}_{m_j},
\end{equation*}
where $p_{a_j}$ and $\widetilde{p}_{m_j}$ are the $j$-th elements of the pressure vectors $\mathbf{p}_a$ and $\widetilde{\mathbf{p}}_{m}$, respectively.
The $L^1$ errors associated to the MRST and updated versions of pEDFM, denoted as $|| E ||^M_1$ and $|| E ||^U_1$, are given respectively by
\begin{equation*}
	|| E ||^M_1 = \frac{13}{50} \frac{R_k |5 - 3 R_k - 40d|}{37 R^2_k + 5 R_k + 37 R_k d + 5 d},
	\qquad
	|| E ||^U_1 = \frac{13}{50} \frac{R^2_k}{21 R^2_k + 20 d^2 + 41 R_k d}. 
\end{equation*}

Fig.~\ref{fig:1D_pEDFM_err} depicts the trend of the errors with respect to the permeability ratio $R_k$ for a fixed value of the fracture aperture $d = 10^{-4}$.
We notice that when the fracture is much more impermeable than the porous matrix the errors tend to be small for both pEDFM versions, while it can be clearly seen that the error of the MRST version of pEDFM is much higher than that of the updated one from low to intermediate (impermeable) conductivity contrast values.

We stress that the error is computed for a fixed coarse grid. Refining the grid would then reduce the value of the aforementioned constant error. 

\begin{figure}[!t]
	\centering
	\begin{subfigure}[b]{0.495\textwidth}
		\centering
		\includegraphics[width=.9\textwidth]{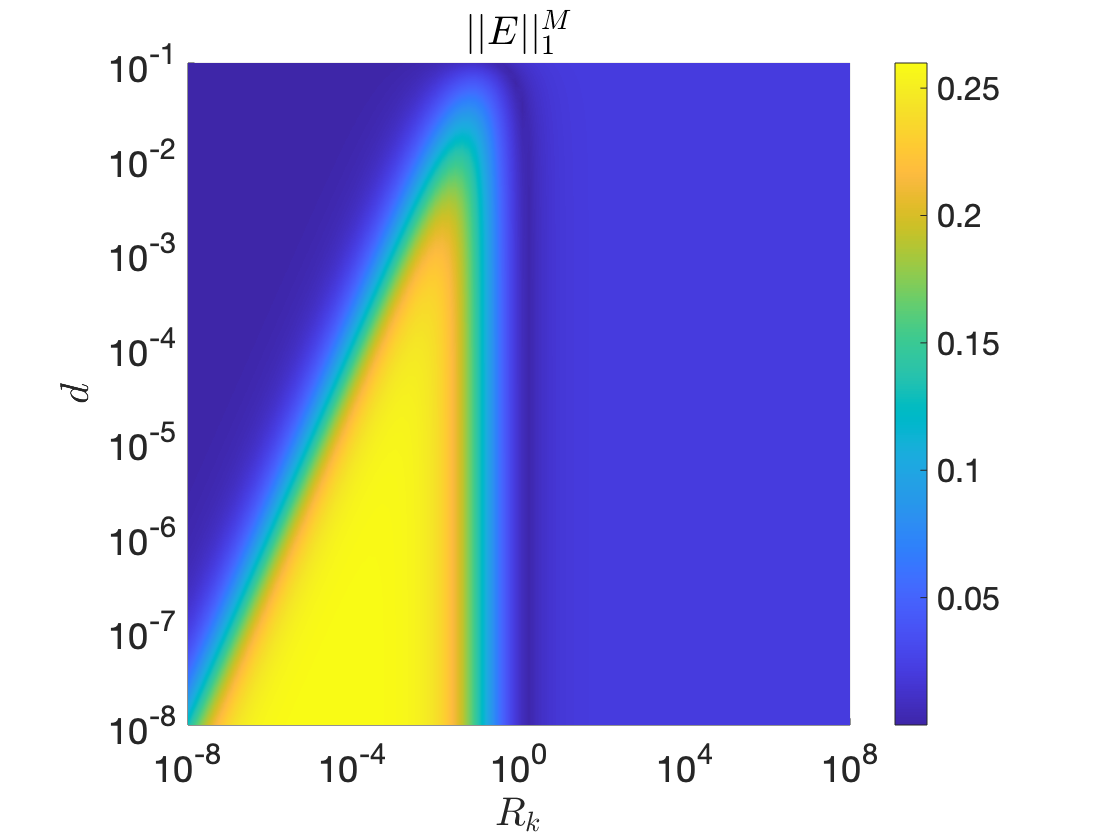}
		\caption{MRST pEDFM.}
		\label{fig:1D_pEDFM_err_o_3}
	\end{subfigure}
	\hfill
	\begin{subfigure}[b]{0.495\textwidth}
		\centering
		\includegraphics[width=.9\textwidth]{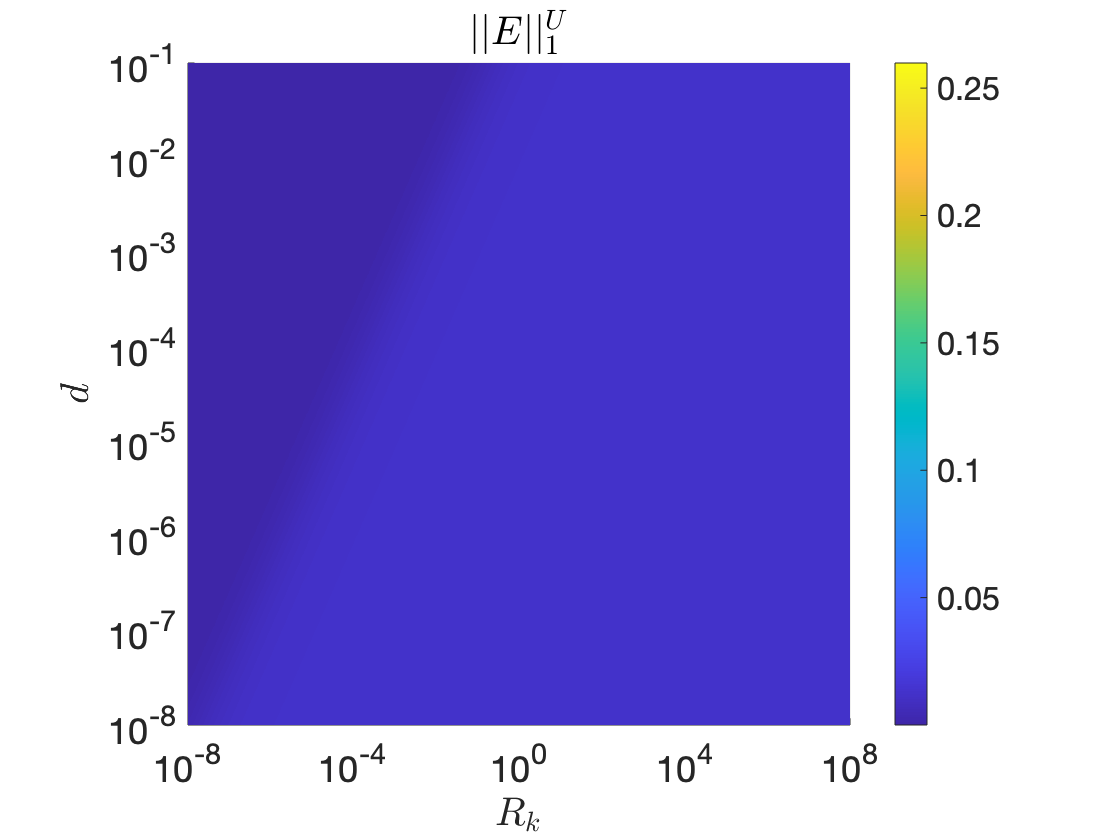}
		\caption{Updated pEDFM.}
		\label{fig:1D_pEDFM_err_i_3}
	\end{subfigure}
	\caption{$L^1$ error maps for variable $R_k$ and $d$.}
	\label{fig:1D_pEDFM_err_3}
\end{figure}

If we also consider the dependence of the error on the fracture aperture $d$ we obtain the 2D maps shown in Fig.~\ref{fig:1D_pEDFM_err_3} for the $L^1$ errors.
The interval of fracture aperture values is chosen to be $[10^{-8}, 10^{-1}]$ to cover cases of both thick and thin fractures with an upper bound equal to $h/2$, so that it is smaller than the grid size $h$.
The $L^1$ error map for the MRST pEDFM (Fig.~\ref{fig:1D_pEDFM_err_o_3}) shows that as $d$ decreases, the model is inaccurate for a wider range of $R_k$.
As $d$ approaches its upper bound value, instead, we obtain errors similar in magnitude to those obtained with the updated pEDFM.
The $L^1$ error for the updated pEDFM, shown in Fig.~\ref{fig:1D_pEDFM_err_i_3}, is always low, albeit depending on the fixed grid size. In particular, as the fracture gets thicker, the interval of impermeable $R_k$ values for which the error is very low widens starting from impermeable fractures.
 
The results obtained show that the adoption of the updated version of pEDFM is more appropriate than the MRST one, especially if one wants to correctly represent the presence of thin, impermeable fractures. 
Assuming that similar arguments also apply to finer 1D grids and multidimensional cases, the updated version of pEDFM should always be preferred over the MRST one.

\section*{Acknowledgments}
The authors acknowledge financial support of Eni S.p.A. and thank Hadi Hajibeygi for the fruitful discussions about pEDFM. 

\printbibliography
	
\end{document}